%% file: main.tex
\documentclass[fleqn,10pt]{article}
\usepackage{latexsym, graphicx, epsfig, amsmath, amssymb,amsfonts}
\usepackage{natbib,amsthm,version}
\usepackage{amsbsy,bm,multirow,enumerate}
\usepackage[titletoc,page]{appendix}
\usepackage[mathscr]{eucal}
\usepackage{mathtools}
\usepackage{color}
\usepackage[utf8]{inputenc}
\usepackage[english]{babel}
\usepackage{amsthm}
\usepackage{enumerate}
\usepackage[hidelinks]{hyperref}
\usepackage{url}
\usepackage{subfigure}
\usepackage[lined,boxed,linesnumbered,ruled]{algorithm2e}
\usepackage{float}

\newcommand{\mbs}[1]{\mathbf{#1}}
\newcommand{\mbb}[1]{\mathbb{#1}}

\newtheorem{theorem}{Theorem}[section]

\newtheorem{remark}{Remark}[section]

\theoremstyle{definition}

\title{
  Numerical Computation of Partial Differential Equations by
  Hidden-Layer Concatenated Extreme Learning Machine
} 
\author{
  Naxian Ni,\ \ Suchuan Dong\thanks{Author of correspondence. Email: sdong@purdue.edu}
  \\
  Center for Computational and Applied Mathematics \\
  Department of Mathematics \\
  Purdue University, USA 
 } 

\date{(\today)}

\begin{document}
\maketitle


\input Abstract


\vspace{0.05cm}



\input Introduction

\input Method

\input Test

\input Summary

\section*{Acknowledgement}
This work was partially supported by
NSF (DMS-2012415). 

\input Append

\bibliographystyle{plain}
\bibliography{elm,elm1,mypub,dnn,sem,obc,varpro}

\end{document}

%% file: Abstract.tex
\begin{abstract}

  The extreme learning machine (ELM) method can yield highly accurate
  solutions to linear/nonlinear partial differential equations (PDEs),
  but requires the last hidden layer of the neural network to be wide
  to achieve a high accuracy. If the last hidden layer
  is narrow, the accuracy of the existing ELM method
  will be poor, irrespective of the rest of the network configuration.
  In this paper we present a modified ELM method, termed HLConcELM (hidden-layer
  concatenated ELM), to overcome the above drawback of the conventional
  ELM method. The HLConcELM method can produce highly accurate solutions
  to linear/nonlinear PDEs when the last hidden layer of the network is narrow
  and when it is wide. The new method is based on a type of modified feedforward
  neural networks (FNN), termed HLConcFNN (hidden-layer concatenated FNN),
  which incorporates a logical concatenation of the hidden layers in
  the network and exposes all the hidden nodes to the output-layer nodes.
  HLConcFNNs have the interesting property that, given a
  network architecture, when additional hidden layers are appended to the network
  or when extra nodes are added to the existing hidden layers
  the representation
  capacity of the HLConcFNN associated with the new architecture is guaranteed
  to be not smaller than that of the original network architecture.
  Here representation capacity refers to the set of all functions that can be
  exactly represented by the neural network of a given architecture.
  We present ample benchmark tests with linear/nonlinear PDEs
  to demonstrate the computational accuracy and performance of the HLConcELM method
  and the superiority of this method to the conventional ELM from
  previous works.

\end{abstract}

%% file: Introduction.tex
\section{Introduction}
\label{sec:intro}

%
%
%

This work continues and extends our recent
studies~\cite{DongL2021,DongL2021bip,DongY2021}
of a type of random-weight neural networks,
the so-called extreme learning machines (ELMs)~\cite{HuangZS2006},
for scientific computing and in particular for computational
partial differential equations (PDEs). 
Specifically, we would like to address the following question:
\begin{itemize}
\item Can ELM achieve a high accuracy for solving linear/nonlinear PDEs
  on network architectures with a narrow last hidden layer?

\end{itemize}
For the existing ELM method~\cite{DongL2021,DongL2021bip,DongY2021}
(referred to as the conventional ELM hereafter), the answer to this
question is negative. The goal of this paper is to
introduce a modified method, referred to as the hidden-layer
concatenated ELM (HLConcELM), to overcome this drawback of
the conventional ELM and provide
a positive answer to the above question.


Exploiting randomization in neural networks has a long history~\cite{ScardapaneW2017}.
Turing's un-organized machine~\cite{Webster2012} and Rosenblatt's
perceptron~\cite{Rosenblatt1958} in the 1950s are early examples of randomized
neural networks. After a hiatus of several decades,
there has been a strong revival of
methods based on random-weight neural networks,
starting in the 1990s~\cite{SuganthanK2021}.
In recent years randomization based neural networks
have attracted a growing interest in a variety of
areas~\cite{ScardapaneW2017,FreireRB2020}.


Since it is enormously costly and hard to optimize the entire set of
adjustable parameters in the neural network, it seems advisable if one randomly
assigns and fixes a subset of the network's parameters so that
the ensuing optimization task of network training
can be simpler, and ideally linear, without severely compromising
the network's achievable approximation capability.
This philosophy underlies the randomization of neural networks.
When applied to feedforward or recurrent neural networks, randomization
leads to techniques such as the random vector functional link (RVFL)
networks~\cite{PaoT1992,PaoPS1994,IgelnikP1995}, the extreme
learning machine~\cite{HuangZS2004,HuangZS2006,HuangCS2006},
the echo-state network~\cite{JaegerLPS2007,LukoseviciusJ2009},
the no-propagation network~\cite{WidrowGKP2013}, and
the liquid state machine~\cite{MaasM2004}.
%
The random-weight neural networks (with a single hidden layer)
are universal function approximators. The universal approximation property
of such networks has been studied 
in~\cite{IgelnikP1995,LiCIP1997,HuangCS2006,NeedellNSS2020}.
The theoretical results of~\cite{IgelnikP1995,HuangCS2006,NeedellNSS2020}
establish that a single hidden-layer feedforward neural network
having random but fixed (not trained) hidden nodes
can approximate any continuous function to any desired
degree of accuracy, provided that the number of hidden nodes
is sufficiently large.
The expected rate of convergence in the approximation
of Lipschitz continuous functions is given
in~\cite{IgelnikP1995,RahimiR2008,NeedellNSS2020}.


ELM was originally developed in~\cite{HuangZS2004,HuangZS2006}
for single hidden-layer
feedforward neural networks 
for linear classification/regression problems.
It has since undergone a dramatic growth and found widespread applications
in a variety of areas (see e.g.~the reviews of~\cite{HuangHSY2015,Alabaetal2019}
and the references therein).
The method is based on two strategies: (i) randomly assigned but fixed (not trainable)
hidden-layer coefficients, and (ii) trainable linear output-layer
coefficients computed by a linear least squares method or by
using the pseudoinverse (Moore-Penrose inverse) of the coefficient
matrix~\cite{VermaM1994,PaoPS1994,BraakeS1995,GuoCS1995}. 


While ELM emerged nearly two decades ago,
the investigation of this technique for the numerical solution of
differential equations has appeared only quite recently,
alongside the proliferation of deep neural
network (DNN) based PDE solvers in the past few
years (see e.g.~\cite{Karniadakisetal2021,SirignanoS2018,RaissiPK2019,EY2018,WinovichRL2019,HeX2019,CyrGPPT2020,JagtapKK2020,WangL2020,WanW2020,DongN2021,LuMMK2021,TangWL2021,KrishnapriyanGZKM2021,WangYP2022},
among many others).
In~\cite{YangHL2018,Sunetal2019,LiuXWL2020,LiuHWC2021} 
the ELM technique has been used for solving linear ordinary or partial
differential equations (ODEs/PDEs) with single hidden-layer
feedforward neural networks, in which 
certain polynomials (e.g.~Legendre, Chebyshev, or Bernstein polynomials)
serve as the activation function.
%
In~\cite{PanghalK2020} the ELM algorithm is used for solving linear
ODEs and PDEs on neural networks with a single hidden layer, in which
the Moore-Penrose inverse of the coefficient matrix has been used.
%
In~\cite{DwivediS2020} a physics-informed ELM method is proposed
for solving linear PDEs  by combining
the physics-informed neural network and the ELM idea.
The neural network consists of a single hidden layer, and the Moore-Penrose
inverse is employed to solve the resultant linear system.
Interestingly, the authors set the number of hidden nodes
to be equal to the total number of conditions in the problem.
A solution strategy based on the normal equation associated
with the linear system is studied in~\cite{DwivediS2022}.

The ELM approach is extended to the numerical solution of nonlinear
PDEs in~\cite{DongL2021} on local or global feedforward neural networks with
a single or multiple hidden layers. A nonlinear least squares method
with perturbations (NLLSQ-perturb) and a Newton-linear least squares
(Newton-LLSQ) method are developed for solving the resultant
nonlinear algebraic system for the output-layer
coefficients of the ELM neural network.
The NLLSQ-perturb algorithm therein takes advantage of the nonlinear
least squares routine from the scipy library, which implements
a Gauss-Newton type method combined with a trust-region strategy.
A block time marching (BTM) scheme is proposed in~\cite{DongL2021}
for long-time dynamic simulations of linear and nonlinear PDE problems,
in which the temporal dimension (if large) is divided into a number of windows
(called time blocks) and the PDE problem is solved on the time blocks
individually and successively.
More importantly, a systematic comparison of the accuracy and the
computational cost (network training time) between the ELM method
and two state-of-the-art deep neural network (DNN) based
PDE solvers, the deep Galerkin method (DGM)~\cite{SirignanoS2018}
and the physics-informed neural network (PINN) method~\cite{RaissiPK2019},
has been conducted in~\cite{DongL2021}, as well as a systematic comparison between
ELM and the classical finite element method (FEM).
The comparisons show that the ELM method far outperforms DGM and PINN
in terms of the accuracy and the computational cost, and that
ELM is on par with the classical FEM in computational performance
and outperforms the FEM as the problem size becomes larger.
In~\cite{DongL2021} the hidden-layer coefficients are set
to uniform random values generated on $[-R_m,R_m]$,
where $R_m$ is a user-prescribed constant.
The results of~\cite{DongL2021} show that
the $R_m$ value has a strong influence on the numerical accuracy of
the ELM results
and that the best accuracy is associated with
a range of moderate $R_m$ values for a given problem.
This is consistent with the observation for classification
problems~\cite{ZhangS2016}.

A number of further developments of the ELM technique for solving linear and nonlinear
PDEs appeared recently; see 
e.g.~\cite{DongL2021bip,CalabroFS2021,GalarisFCSS2021,DongY2021,FabianiCRS2021},
among others.
In order to address the influence of random initialization of the hidden-layer coefficients
on the ELM accuracy,
a modified batch intrinsic plasticity (modBIP) method is developed
in~\cite{DongL2021bip} for pre-training the random coefficients in the ELM
network. This method, together with ELM, is applied to a number of linear and
nonlinear PDEs. The accuracy of the combined modBIP/ELM method has been shown to be
insensitive to the random initializations of the hidden-layer coefficients.
In~\cite{CalabroFS2021} the authors present a method for
solving one-dimensional linear elliptic
PDEs  based on ELM with single hidden-layer feedforward neural networks
and the sigmoid activation function. The random parameters in the activation function
are set based on the location of the domain of interest and the function
derivative information.
In~\cite{GalarisFCSS2021} the authors present a method based on randomized
neural networks with a single hidden layer for solving stiff ODEs.
The time integration
therein appears similar to the block time marching
strategy~\cite{DongL2021} but with an adaptation
on the time block sizes. It is observed that the presented
method is advantageous over
the stiff ODE solvers from MatLab.
Noting the influence of the maximum random-coefficient magnitude (i.e.~the $R_m$ constant) on
the ELM accuracy as shown by~\cite{DongL2021},
in~\cite{DongY2021} we have presented a method for computing the optimal $R_m$
in ELM based on the differential evolution algorithm, as well as an improved
implementation for computing the differential operators of the last hidden-layer
data. These improvements significantly enhance
the ELM computational performance and dramatically reduce its network training
time as compared with that of~\cite{DongL2021}.
The improved ELM method is compared systematically with the traditional
second-order and high-order
finite element methods for solving a number of linear and nonlinear PDEs in~\cite{DongY2021}.
The improved ELM far outperforms the second-order FEM.
For smaller problem sizes it is comparable to
the high-order FEM in performance, and for larger
problem sizes the improved ELM outperforms the high-order FEM markedly.
Here, by ``outperform'' we mean that one method achieves a better accuracy under
the same computational cost or incurs a lower computational cost to achieve
the same accuracy.
In~\cite{FabianiCRS2021} an ELM method is presented for the numerical
solution of stationary nonlinear PDEs based on the sigmoid and radial
basis activation functions. The authors 
observe that the ELM method exhibits a better accuracy than the finite difference
and the finite element method.
%
Another recent development related to ELM is~\cite{DongY2022},
in which a method based on the variable projection strategy
is proposed for solving linear and nonlinear PDEs with artificial neural networks.
For linear PDEs, the neural-network representation of the PDE solution leads to a
separable nonlinear least squares problem, which
is then reformulated to eliminate the output-layer coefficients, leading to
a reduced problem about the hidden-layer coefficients only. The reduced problem
is solved first by the nonlinear least squares method to determine the hidden-layer
coefficients, and the output-layer coefficients are then computed by
the linear least squares method~\cite{DongY2022}.
For nonlinear PDEs, the problem is first linearized by the Newton's method with
a particular linearization form, and the linearized system is solved by
the variable projection framework together with neural networks.
The ELM method can be considered as a special case of the variable projection,
i.e.~with zero iteration when solving the reduced problem for
the hidden-layer coefficients~\cite{DongY2022}. It is shown in~\cite{DongY2022} that
the variable projection method exhibits an accuracy significantly superior to the ELM
under identical conditions and network configurations.


As has been shown in previous works~\cite{DongL2021,DongL2021bip,DongY2021},
ELM can produce highly accurate results for solving linear and
nonlinear PDEs. For smooth field solutions
the ELM errors decrease exponentially as the number of training data points
or the number of training parameters in the neural network increases, and
the errors can reach a level close to the machine zero when
the number of degrees of freedom becomes large.
To achieve a high accuracy, however, the existing ELM method requires the
number of nodes in the last hidden layer of the neural network
to be sufficiently large~\cite{DongL2021}.
Therefore, the ELM network usually has a wide hidden layer in the case
of a shallow neural network, or a wide
last hidden layer in the case of deeper neural networks.
If the last hidden layer  contains
only a small number of nodes, the
results computed with the existing (conventional) ELM
will tend to be poor in accuracy,
regardless of the configuration with the rest of the network.

In this paper, we focus on feedforward neural networks (FNNs) with multiple hidden layers,
and present a modified ELM method (termed HLConcELM) for solving linear/nonlinear
PDEs to overcome
the above drawback associated with the conventional ELM method.
The HLConcELM method can produce accurate solutions to
linear/nonlinear PDEs when the last hidden layer of the network is narrow,
and when the last hidden layer is wide.

The new method is based on a type of modified feedforward neural networks,
referred to as hidden-layer concatenated FNN (or HLConcFNN) herein, which incorporates
a logical concatenation of  the hidden layers so that all the hidden nodes
are exposed to and connected with
the nodes in the output layer (see Section~\ref{sec:method} for details).
The HLConcFNNs have the interesting property that,
given a network architecture, when additional hidden layers are appended
to the neural network or when extra nodes are added to the existing hidden layers,
the representation capacity of the HLConcFNN associated with
the new architecture is guaranteed
to be not smaller than that associated with the original network architecture.
Here by representation compacity we refer to the set of all functions that
can be exactly represented by the neural network (see Section~\ref{sec:method}
the definition).
In contrast, conventional FNNs do not have a parallel property when additional
hidden layers are appended to the network.

The HLConcELM is attained by setting (and fixing) the weight/bias coefficients in
the hidden layers of the HLConcFNN to random values
and allowing the connection coefficients
between all the hidden nodes and the output nodes to be adjustable (trainable).
More specifically, given a network architecture with $L$ hidden layers,
we set the weight/bias coefficients of the $l$-th ($1\leqslant l\leqslant L$)
hidden layer to uniform random values generated on the interval $[-R_l,R_l]$,
where $R_l$ is a constant. The vector of $R_l$ constants (referred to as
the hidden magnitude vector herein), $\mbs R=(R_1,R_2,\dots,R_L)$,
influences the accuracy of HLConcELM, and in this paper we determine
the optimal $\mbs R$ using the method from~\cite{DongY2021} based on
the differential evolution algorithm.
HLConcELMs partly inherit the non-decreasing representation capacity property
of HLConcFNNs. For example,
given a network architecture, when extra hidden layers are appended to
the network, the representation capacity of the HLConcELM associated with
the new architecture will not be smaller than that associated with the original
architecture, provided that the random hidden-layer coefficients
for the new architecture are assigned in an appropriate fashion.
On the other hand,
when extra nodes are added to the existing hidden layers, HLConcELMs
in general do not have a non-decreasing property with regard to its representation
capacity, because of the randomly assigned hidden-layer coefficients.


The exploration of neural-network architecture has been actively pursued in
machine learning research, and the connectivity patterns are
the focus of a number of research efforts.
Neural networks incorporating shortcut connections (concatenations)
between the input nodes, the hidden
nodes, and the output nodes are explored
in e.g.~\cite{HuangLMW2018,WilamowskiY2010,CortesGKMY2016,KatuwalST2019}
(among others). The hidden-layer concatenated neural network adopted
in the current paper can be considered in spirit as a simplification of
the connection patterns in
the DenseNet~\cite{HuangLMW2018} architecture, and it is similar to the
deep RVFL architecture of~\cite{KatuwalST2019} but without the
connection between the input nodes and the output nodes.
We note that these previous works are for image and data classification problems,
while the current work focuses on scientific computing and in particular the
numerical solutions of partial differential equations.


We present extensive numerical experiments with several
linear and nonlinear PDEs  to test the performance
of the HLConcELM method and compare  the current method
with the conventional ELM method.
These benchmark tests demonstrate unequivocally that HLConcELM
can achieve highly accurate results when the last hidden layer
in the neural network is narrow or wide, and that it is much superior in accuracy to
the conventional ELM method.
The implementation of the current method is in Python and employs the Tensorflow
(www.tensorflow.org), Keras (keras.io), and the scipy libraries.
All the benchmark tests are performed on a MAC computer (3.2GHz Intel Core i5 CPU,
24GB memory) in the authors' institution.


The contribution of this work lies in two aspects. The main contribution
lies in the HLConcELM method
for solving linear and nonlinear PDEs. The other aspect is with regard to
the non-decreasing representation capacity of HLConcFNNs when additional
hidden layers are appended to an existing network architecture.
To the best the authors' knowledge, this property seems unknown to the community
so far. Bringing this property into collective consciousness can be
another contribution of this paper.

The rest of this paper is organized as follows.
In Section~\ref{sec:method} we discuss the structures of HLConcFNNs and
HLConcELMs, as well as their non-decreasing representation capacity property
when additional hidden layers are appended to an existing architecture.
We then develop the algorithm for solving linear and nonlinear PDEs
employing the HLConcELM architecture.
In Section~\ref{sec:tests} we present extensive benchmark examples
to test the current HLConcELM method
and compare  its performance with that of
the conventional ELM.
Section~\ref{sec:summary} concludes the presentation with several further comments
about the presented method.
In the Appendix we include constructive proofs to some
theorems from Section~\ref{sec:method}
concerning the representation capacity of HLConcFNNs.


%% file: Method.tex
\section{Hidden-Layer Concatenated Extreme Learning Machine}
\label{sec:method}

%
%
%

%

\subsection{Conventional ELM and Drawback}

\begin{figure}
  \centerline{
    \includegraphics[height=1.2in]{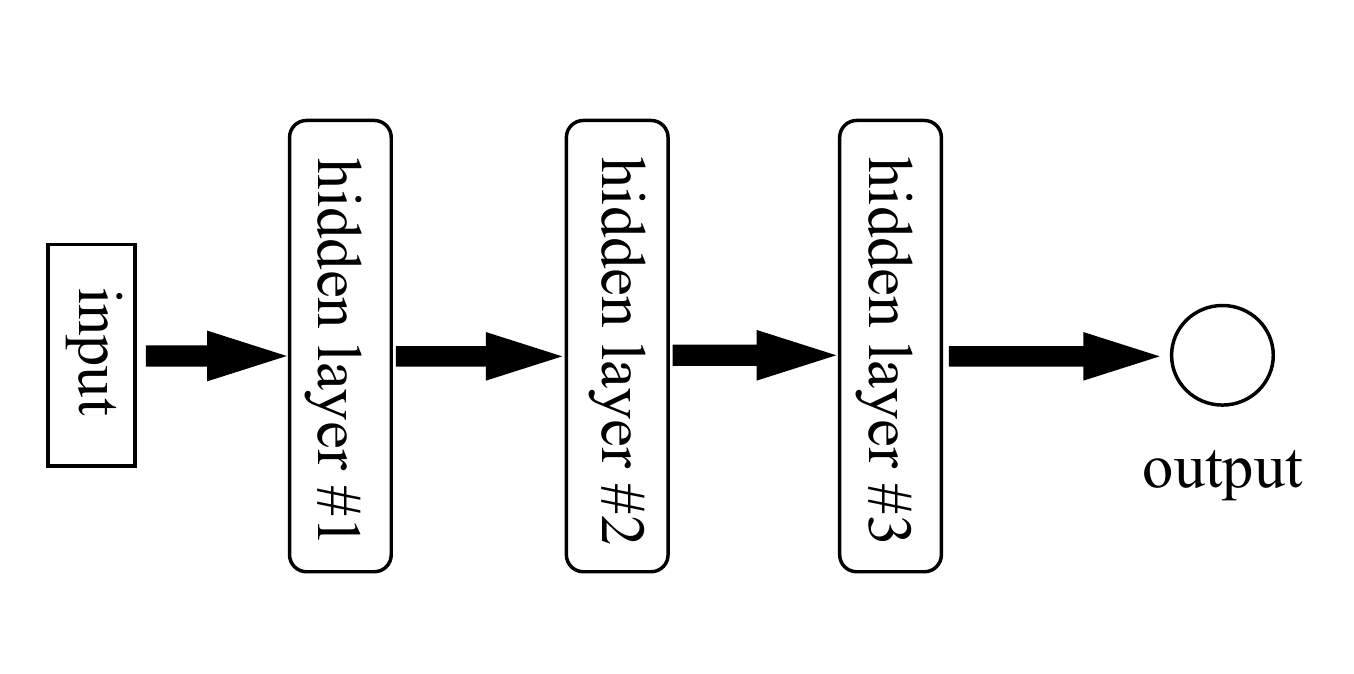}(a)\qquad
    \includegraphics[height=1.2in]{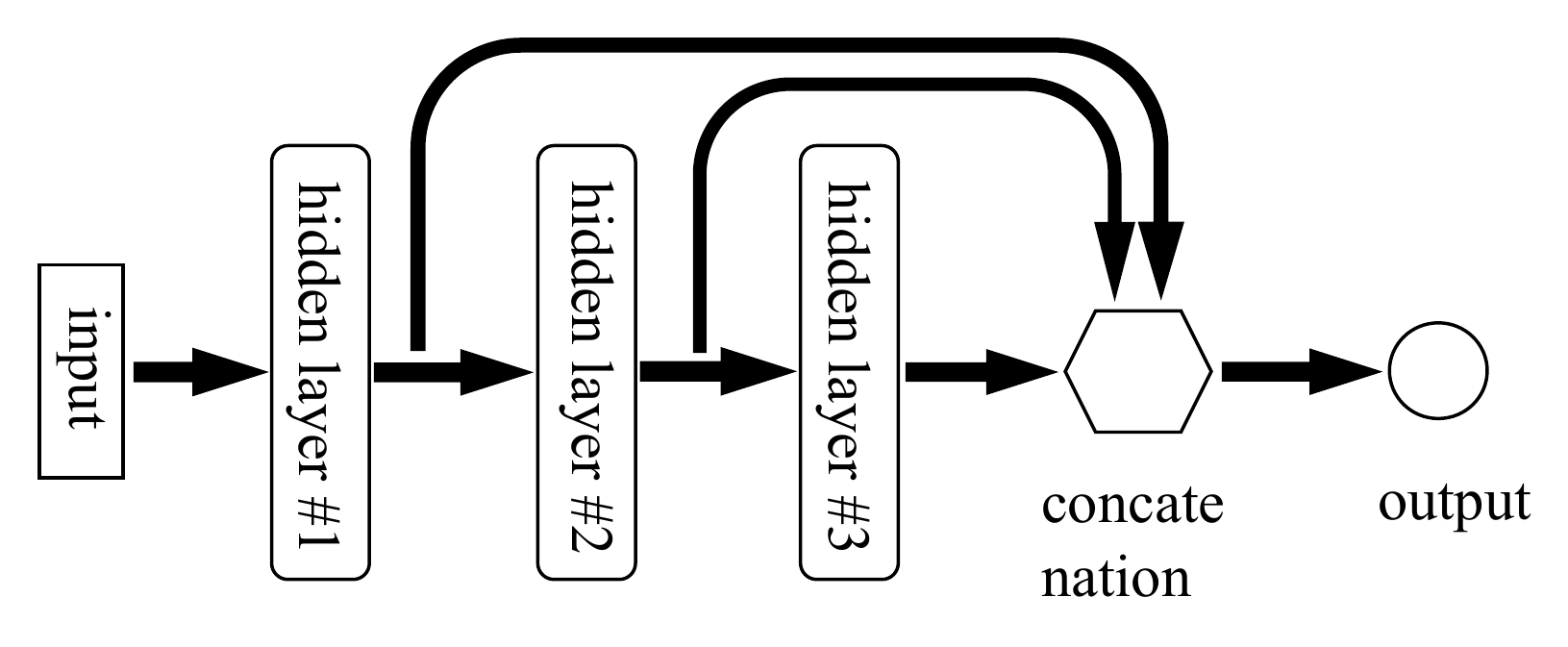}(b)
  }
  \caption{Illustration of neural network structure (with $3$ hidden layers):
    (a) conventional FNN, and
    (b) hidden-layer concatenated FNN (HLConcFNN).
    In HLConcFNN all the hidden nodes are exposed
    to the output nodes, while in conventional FNN only
    the last hidden-layer nodes are exposed to the output nodes.
  }
  \label{fg_1}
\end{figure}

The ELM method with feedforward
neural networks (FNN) for solving linear/nonlinear PDEs
has been described in e.g.~\cite{DongL2021,DongL2021bip,DongY2021}.
Figure~\ref{fg_1}(a) illustrates such a network containing three
hidden layers. From layer to layer, the arrow in the sketch
represents the usual FNN logic,
an affine transform followed by a function composition
with an activation function~\cite{GoodfellowBC2016}.
For ELM we require that no activation function
be applied to the output layer and that
the output layer contain no bias.
So the output layer is linear and has zero bias.
This requirement is adopted throughout this paper.

As discussed in~\cite{DongL2021},
we pre-set the weight/bias coefficients in all
the hidden layers to random values and fix these values once
they are set. Only the weight
coefficients of the output layer are trainable. The hidden-layer
coefficients in the neural network are not trainable with ELM.

To solve a given linear or nonlinear PDE with ELM,
we first enforce the PDE and the associated boundary/initial conditions
on a set of collocation points in the domain or on the
appropriate domain boundaries. This gives rise to
a linear least squares problem for linear PDEs, or a
nonlinear least squares problem for nonlinear PDEs, about the
output-layer coefficients (trainable parameters)
of the neural network~\cite{DongL2021}.
We solve this
least squares problem for the output-layer coefficients
by a linear least squares method for linear PDEs and
by a nonlinear least squares method for nonlinear PDEs~\cite{DongL2021}.

ELM can produce highly accurate solutions
to PDEs. In particular, for smooth solutions its errors decrease
exponentially as the number of collocation points or
the number of trainable parameters (the number of
nodes in the last hidden layer) increases~\cite{DongL2021,DongL2021bip}.
In addition, it has a low computational cost
(network training time)~\cite{DongL2021,DongY2021}.


%

Hereafter we refer to a vector or a list of positive integers as
an architectural vector (denoted by $\mbs M$),
\begin{equation} \label{eq_1}
  \text{architectural vector:}\quad
  \mbs M = [M_0, M_1, \dots, M_{L-1}, M_L]
\end{equation}
where $(L+1)$ is the dimension of the vector with $L\geqslant 2$,
and $M_i$ ($0\leqslant i\leqslant L$) are positive integers.
We associate a given  $\mbs M$
to the architecture
of an FNN with $(L+1)$ layers,
where $M_i$ ($0\leqslant i\leqslant L$) is the number of nodes
in the $i$-th layer. The layer $0$ and the layer $L$ represent
the input and the output layers, respectively. The layers in between
are the hidden layers.

\begin{figure}
  \centerline{
    \includegraphics[width=1.8in]{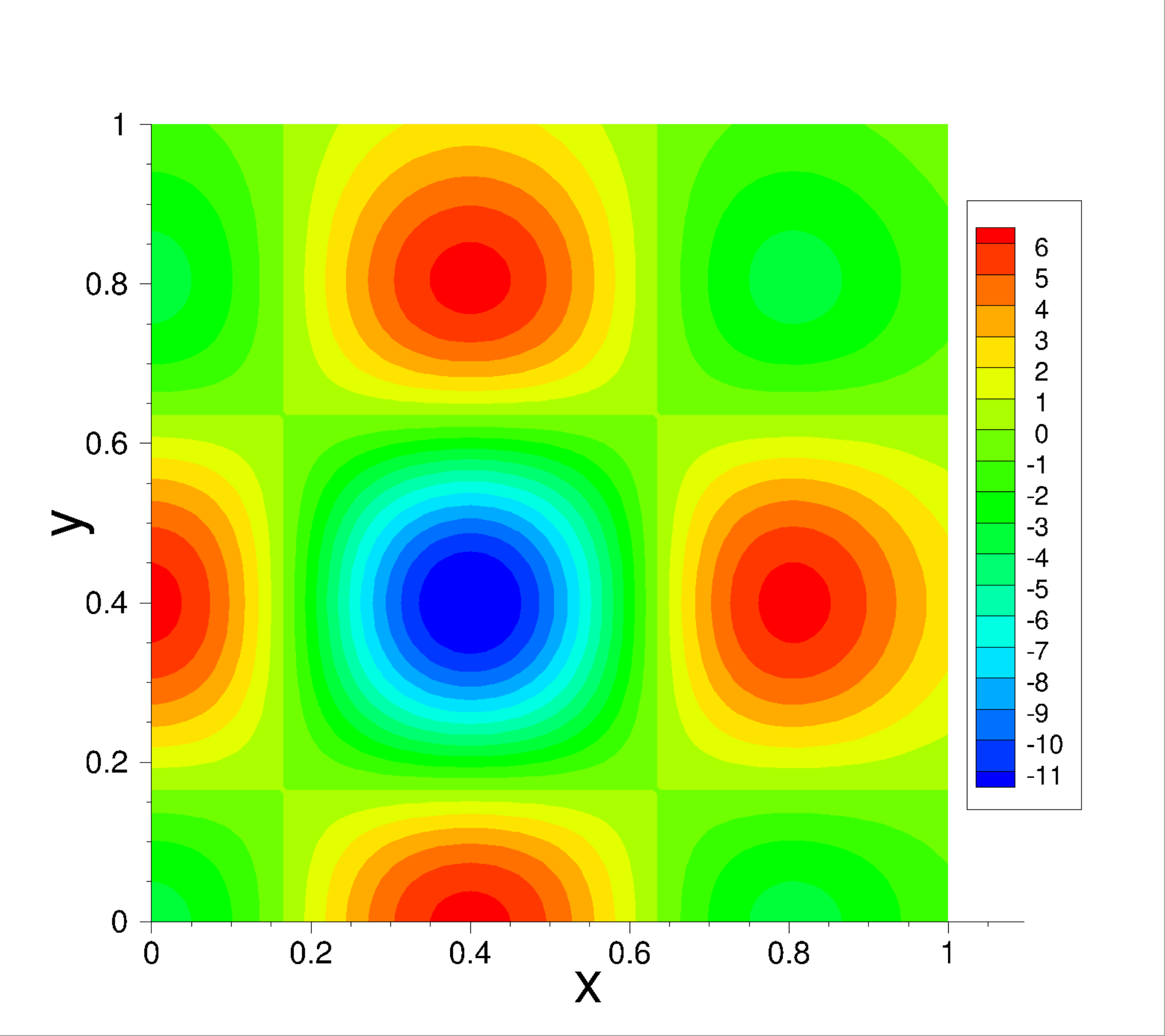}(a)
    \includegraphics[width=1.8in]{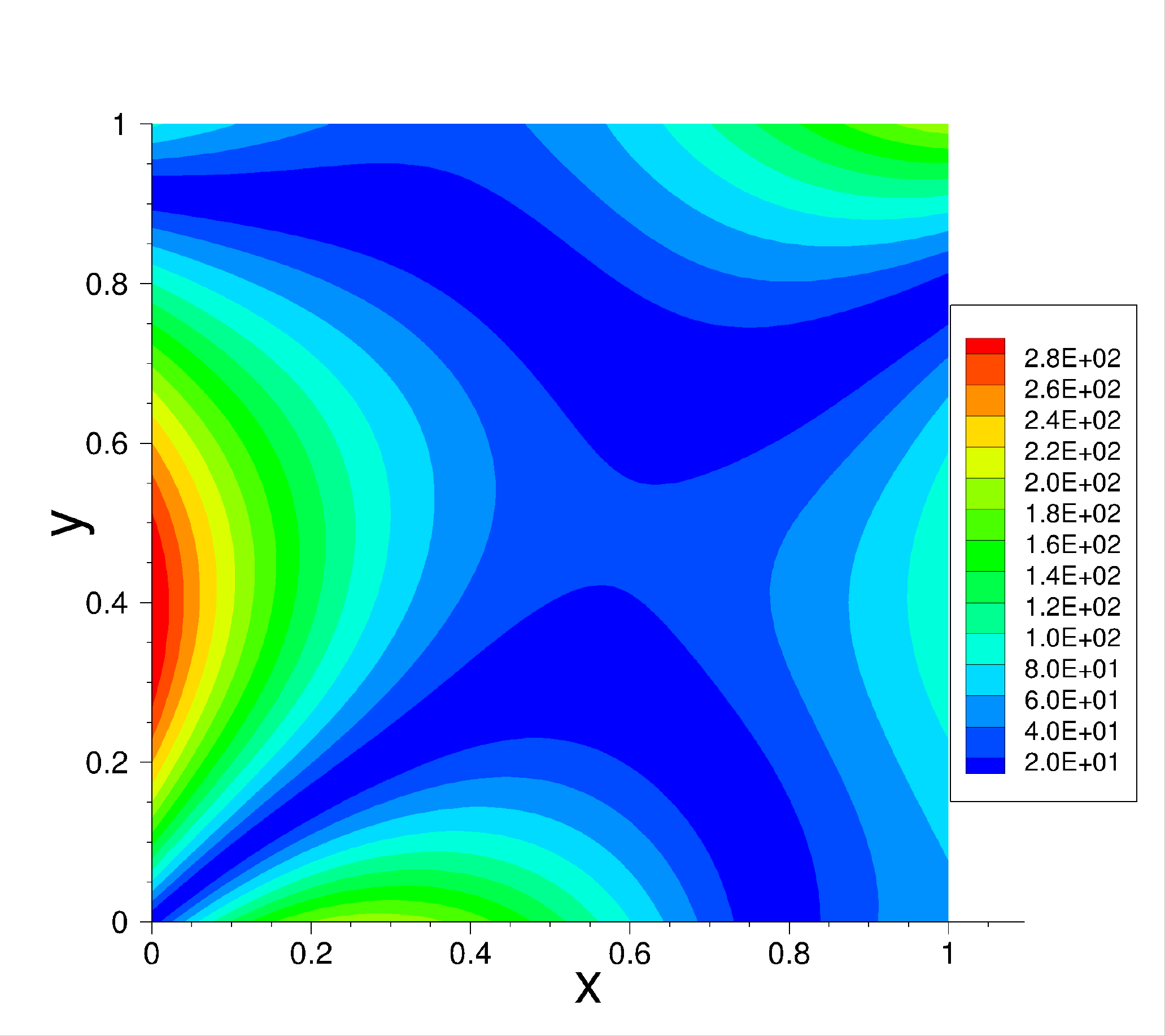}(b)
    \includegraphics[width=1.8in]{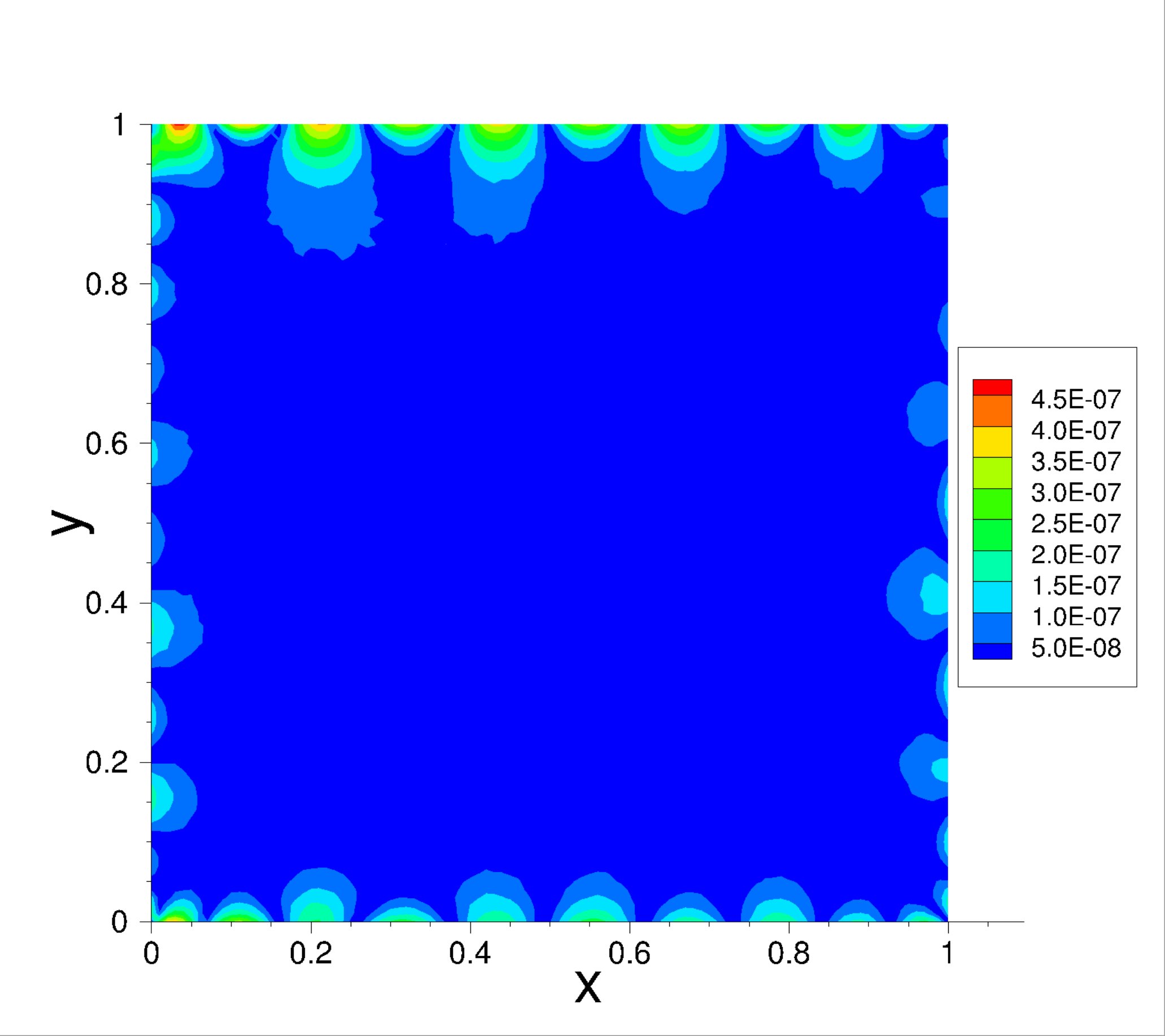}(c)
  }
  \caption{Illustration of error characteristics of conventional ELM (Poisson equation):
    Distributions of (a) the exact solution, (b) the ELM error obtained
    with the architecture $\mbs M_1=[2, 300, 30, 1]$, and (c) the ELM error obtained
    with the architecture $\mbs M_2=[2, 300, 400, 1]$.
  }
  \label{fg_2}
\end{figure}

Despite its high accuracy and attractive computational performance,
certain aspect of the ELM
method is less appealing and remains to be improved. One particular aspect in
this regard concerns the size of the last hidden layer of the ELM network.
ELM  requires the last hidden layer of the neural network
to be wide in order to achieve a high accuracy,
irrespective of the sizes of the rest of the network architecture.
If the last hidden layer contains only
a small number of nodes,
the ELM accuracy will be poor even though the preceding
hidden layers can be wide enough. This point is illustrated by
Figure~\ref{fg_2}, which shows the ELM results for
solving the two-dimensional (2D) Poisson equation on
a unit square. Figure~\ref{fg_2}(a) shows the distribution of
the exact solution. Figures~\ref{fg_2}(b) and (c) show
the ELM error distributions  obtained using
two network architectures given by
$[2, 300, 30, 1]$ and $[2, 300, 400, 1]$, respectively,
under otherwise identical conditions.
Both neural networks have the Gaussian
activation function $\sigma(x)=e^{-x^2}$ for all the hidden nodes,
and are trained on a uniform set of $21\times 21$ collocation points.
The only difference between them is the size of the last hidden layer.
With $400$ nodes in the last hidden layer the ELM solution is
highly accurate, with the maximum error on the order $10^{-7}$ in
the domain. With $30$ nodes in the last hidden layer,
on the other hand, the ELM solution exhibits no accuracy at all,
with the maximum error on the order of $10^2$,
despite the fact that the first hidden layer is fairly large (containing
$300$ nodes).
With the existing ELM method, only the last hidden-layer nodes
directly contribute to the the output of the neural network,
while the nodes
in the preceding hidden layers do not directly affect the network output.
In other words, with the existing ELM,
all the degrees of freedom provided by the nodes in
the preceding hidden layers are to some extent ``wasted''.

Can one achieve a high accuracy even if the last hidden layer is narrow
in the ELM network?
Can we  take advantage of the degrees of freedom provided by
the hidden nodes in the preceding hidden layers with ELM?
%
These are the questions we are interested in and would like to
address in the current work.


The above drawback of the existing ELM method, which will be referred to
as the conventional ELM hereafter, motivates the developments
in what follows. We present a modified ELM method to address this issue
and discuss how to use the modified method for numerical simulations of
PDEs.

\subsection{Modifying ELM Neural Network with Hidden-Layer Concatenation}

To address the aforementioned drawback, we consider a type of modified
FNNs for ELM computation.
The idea of the modified network is illustrated in Figure~\ref{fg_1}(b)
using three hidden layers as an example.

The main strategy here is to expose {\em all} the hidden nodes
in the neural network to the output-layer nodes. Starting with a standard FNN,
 we incorporate a logical concatenation layer between the last hidden layer and
the output layer.
This logical layer concatenates the output fields of all the hidden nodes,
from the first to the last
hidden layers, in the original network architecture.
From the logical concatenation layer to the output layer a usual affine
transform, together with possibly a function composition with an activation function, 
is performed to attain the output fields of the overall neural network.
Note that the logical concatenation layer involves no parameters.

Hereafter we refer to this type of modified neural networks
as the hidden-layer concatenated FNN (HLConcFNN), and
the original FNN as the base neural network.
Thanks to the logical concatenation, in HLConcFNN all the hidden nodes in
the base network architecture are connected with the output nodes.

One can also include the input fields in the logical concatenation layer. Numerical
experiments show that, however, there is no advantage in terms of
the accuracy when the input fields are included. In the current paper
we do not include the input fields in the concatenation.


Let us next use a real-valued function of $d$ ($d\geqslant 1$) variables,
$u(\mbs x)$ ($\mbs x \in\Omega\subset \mbb{R}^d$),
represented by a HLConcFNN
to illustrate some of its properties.
Consider a HLConcFNN, whose base architecture is given by
 $\mbs M$ in~\eqref{eq_1}, where $M_0=d$ and $M_L=1$.
The $d$ input nodes represent the $d$  components of
$\mbs x=(x_1,x_2,\dots,x_d)$, and the single output node represents the
function $u(\mbs x)$. Let $\sigma: \mbb R\rightarrow\mbb R$ denote
the activation function for all the hidden nodes.
As stated before, we require that no activation function be applied
to the output node and that it contain no bias.

Let $\bm \Phi^{(i)}(\mbs x)=\left(\phi^{(i)}_1,\dots,\phi^{(i)}_{M_i}\right)$,
$1\leqslant i\leqslant L-1$,
denote the $M_i$ output fields of the $i$-th hidden layer.
The logical concatenation layer contains a total of
$N_c(\mbs M) = \sum_{i=1}^{L-1}M_i$ logical nodes. 
Then we have the following expansion relation,
\begin{equation}\label{eq_2}
  u(\mbs x) = \sum_{i=1}^{L-1}\sum_{j=1}^{M_i}\beta_{ij}\phi_j^{(i)}(\mbs x)
  =\sum_{i=1}^{L-1} \bm\Phi^{(i)}(\mbs x)\bm\beta_i^T
  =\bm\Phi(\mbs x)\bm\beta^T,
\end{equation}
where $\beta_{ij}$ ($1\leqslant i\leqslant L-1$, $1\leqslant j\leqslant M_i$)
denotes the weight coefficient of the output layer, i.e.~the connection coefficient
between the output node and the $j$-th hidden node in the $i$-th hidden layer, and
\begin{equation}\label{eq_3}
  \left\{
  \begin{split}
    &
    \bm\beta_i=(\beta_{i1},\beta_{i2},\dots,\beta_{iM_i}), \quad
    \bm\beta = (\bm\beta_1, \bm\beta_2,\dots,\bm\beta_{L-1})
    =(\beta_{11},\dots,\beta_{1M_1},\beta_{21},\dots,\beta_{L-1,M_{L-1}}),
    \\
    &
    \bm\Phi = \left(\bm\Phi^{(1)}, \bm\Phi^{(2)},\dots,\bm\Phi^{(L-1)}\right)
    =\left(\phi^{(1)}_1,\dots,\phi^{(1)}_{M_1},\phi^{(2)}_1,\dots,\phi^{(L-1)}_{M_{L-1}}\right).
  \end{split}
  \right.
\end{equation}

The logic from layer $(i-1)$ to layer $i$, for $1\leqslant i\leqslant L-1$,
represents an affine transform followed by a function composition with
the activation function,
\begin{equation}\label{eq_4}
  \left\{
  \begin{split}
    &
    \phi_j^{(i)}(\mbs x) = \sigma\left(\sum_{k=1}^{M_{i-1}}\phi^{(i-1)}_k(\mbs x)w_{kj}^{(i)} + b_j^{(i)}\right),
    \quad 1\leqslant j\leqslant M_i, \\
    &
    \bm\Phi^{(i)}(\mbs x) = \sigma\left(\bm\Phi^{(i-1)}(\mbs x)\mbs W^{(i)} + \mbs b^{(i)}\right).
  \end{split}
  \right.
\end{equation}
In the above equation, the constants
$w_{kj}^{(i)}$ ($1\leqslant k\leqslant M_{i-1}$, $1\leqslant j\leqslant M_i$)
are the weights and $b_j^{(i)}$ ($1\leqslant j\leqslant M_i$) are
the biases of layer $i$, and
\begin{equation}\label{eq_5}
\mbs W^{(i)}=\begin{bmatrix} w_{kj}^{(i)} \end{bmatrix}_{M_{i-1}\times M_i}, \quad
\mbs b^{(i)}=\left(b_1^{(i)},b_2^{(i)},\dots,b_{M_i}^{(i)}\right), \quad
\mbs \Phi^{(0)}(\mbs x)=\left(\phi_1^{(0)},\phi_2^{(0)},\dots,\phi_{M_0}^{(0)}\right)=\mbs x.
\end{equation}
%
%
Define
\begin{equation}\label{eq_6}
  \left\{
  \begin{split}
    &
    \bm\theta^{(i)} = \text{\bf{flatten}}\left[\mbs W^{(i)}, \mbs b^{(i)}  \right],
    \quad 1\leqslant i\leqslant L-1,
    \\
    &
    \bm\theta = (\bm\theta^{(1)}, \bm\theta^{(2)}, \dots,\bm\theta^{(L-1)} )
    =(\theta_1,\theta_2,\dots,\theta_{N_h}),
  \end{split}
  \right.
\end{equation}
where ``{\bf flatten}'' denotes the operation to reshape and combine a list of matrices
or vectors into a single vector, and $N_h(\mbs M)=\sum_{i=1}^{L-1}(M_{i-1}+1)M_i$.
Here $\bm\theta^{(i)}$ denotes the vector of weight/bias coefficients
of layer $i$ for $1\leqslant i\leqslant L-1$, $\bm\theta$
denotes the vector of weight/bias coefficients in all the hidden layers,
and $N_h$ is the total number of hidden-layer coefficients
in the neural network.

The output field of the neural network depends on $(\bm\theta,\bm\beta)$,
and the output fields of each hidden layer depend on $\bm\theta$.
To make these dependencies more explicit, we re-write equation~\eqref{eq_2} into
\begin{equation}\label{eq_7}
  u(\bm\theta,\bm\beta,\mbs x)
  =\sum_{i=1}^{L-1}\sum_{j=1}^{M_i}\beta_{ij}\phi^{(i)}_j(\bm\theta,\mbs x)
  = \sum_{i=1}^{L-1}\bm\Phi^{(i)}(\bm\theta,\mbs x)\bm\beta_i^T
  = \bm\Phi(\bm\theta,\mbs x)\bm\beta^T,
\end{equation}
where $\bm\Phi$, $\bm\Phi^{(i)}$, $\bm\beta_i$ and $\bm\beta$ are
defined in~\eqref{eq_3}.

A hidden-layer concatenated FNN is characterized by
the architectural vector of the base network
and the activation function.
Given the architectural vector $\mbs M$
and an activation function $\sigma$, let HLConcFNN$(\mbs M,\sigma)$ denote
the associated hidden-layer concatenated neural network.
For a given domain $\Omega\subset\mbb R^d$, an architectural
vector $\mbs M=(M_0,M_1,\dots,M_L)$ with $M_0=d$ and $M_L=1$,
and an activation function $\sigma(\cdot)$, we define
\begin{equation}\label{eq_8}
  U(\Omega,\mbs M,\sigma) = \left\{
  u(\bm\theta,\bm\beta,\mbs x)\ |\ u(\bm\theta,\bm\beta,\mbs x)\
  \text{is output of HLConcFNN}(\mbs M,\sigma),\
  \mbs x\in\Omega,\ \bm\theta\in\mbb R^{N_h},\
  \bm\beta\in\mbb R^{N_c}
  \right\}
\end{equation}
as the collection of all the possible output fields
of this  HLConcFNN$(\mbs M,\sigma)$.
Note that $U(\Omega,\mbs M,\sigma)$ denotes the set of all functions that
can be exactly represented by this HLConcFNN$(\mbs M,\sigma)$ on $\Omega$.
Hereafter we refer to $U(\Omega,\mbs M,\sigma)$ as the representation
capacity of the HLConcFNN$(\mbs M,\sigma)$ for the domain $\Omega$.

\begin{remark}\label{rem_m1}
  It should be noted that $U(\Omega,\mbs M,\sigma)$ as defined by~\eqref{eq_8}
  is a set, not a linear
  space, for the simple fact that it is not closed under addition because
  of the nonlinear parameters $\bm\theta$.
\end{remark}

The HLConcFNNs have an interesting property.
If one appends extra hidden layers to the network architecture,
or adds nodes to any of the existing hidden layers,
the representation capacity of the resultant HLConcFNN is at least as
large as that of the original one.
On the other hand, conventional FNNs lack such a property when additional hidden
layers are appended to the architecture.
Specifically, we have the following results.

\begin{theorem}\label{thm_1}
  Given an architectural vector $\mbs M_1=(m_0,m_1,\dots,m_{L-1},m_L)$
  with $m_L=1$, define a new vector $\mbs M_2=(m_0,m_1,\dots,m_{L-1},n,m_L)$, where
  $n\geqslant 1$ is an integer.
  For a given domain $\Omega\subset\mbb R^{m_0}$ and an activation function
  $\sigma(\cdot)$, the following relation holds
  \begin{equation}\label{eq_9}
    U(\Omega,\mbs M_1,\sigma) \subseteq U(\Omega,\mbs M_2,\sigma),
  \end{equation}
  where $U$ is defined by~\eqref{eq_8}.

\end{theorem}

\begin{theorem}\label{thm_2}
  Given an architectural vector $\mbs M_1=(m_0,m_1,\dots,m_L)$
  with $m_L=1$, define a new vector $\mbs M_2=(m_0,m_1,\dots,m_{s-1},m_s+1,m_{s+1},\dots,m_L)$
  for some $s$ ($1\leqslant s\leqslant L-1$).
  For a given domain $\Omega\subset\mbb R^{m_0}$ and an activation function
  $\sigma(\cdot)$, the following relation holds
  \begin{equation}\label{eq_10}
    U(\Omega,\mbs M_1,\sigma) \subseteq U(\Omega,\mbs M_2,\sigma),
  \end{equation}
  where $U$ is defined by~\eqref{eq_8}.
  
\end{theorem}
\noindent These properties can be shown to be true by simple constructions.
These constructions are straightforward and border on being trivial.
Risking on the side of naivety, we still include
the constructive proofs for these two theorems
in an Appendix of this paper for the benefit of a skeptical reader.

It should be noted that for conventional FNNs the relation given by~\eqref{eq_10}
is  true, but the relation given by~\eqref{eq_9} does not hold.
Relation~\eqref{eq_9} is true for HLConcFNNs thanks to
the concatenation of hidden layers in such networks.

Suppose we start with a base neural network architecture $\mbs M_0$
and generate a sequence of architectures $\mbs M_i$ ($i\geqslant 1$),
with each one obtained either by adding extra nodes to the existing
hidden layers of or
by appending additional hidden layers to the previous architecture.
Then based on the above two theorems the HLConcFNNs associated
with this sequence of architectures 
exhibit a hierarchical structure, in the sense that the representation capacities
of this sequence of HLConcFNNs do not decrease, namely
\begin{equation}
  U(\Omega,\mbs M_0,\sigma) \subseteq
  U(\Omega,\mbs M_1,\sigma) \subseteq \cdots \subseteq
  U(\Omega,\mbs M_n,\sigma) \subseteq \cdots.
\end{equation}
If the activation function $\sigma(\cdot)$ is nonlinear,
the representation capacities of this sequence of HLConcFNNs should 
strictly increase.

\begin{remark}\label{rem_0}
  If the number of nodes in the output layer of the HLConcFNN
  is more than one, the relations given by Theorems~\ref{thm_1} and
  \ref{thm_2} about the representation capacities equally hold.

\end{remark}


\paragraph{Hidden-Layer Concatenated Extreme Learning Machine (HLConcELM)}

Let us next combine the hidden-layer concatenated feedforward neural networks
with the idea of ELM.
We adopt HLConcFNNs as the neural network for the ELM computation.
We pre-set (and fix) all the weight/bias coefficients in the hidden layers (i.e.~$\bm\theta$)
of the HLConcFNN to random values, and train/compute the output-layer
coefficients (i.e.~$\bm\beta$) by a linear or nonlinear least squares method.
We will refer to the resultant method as the
hidden-layer concatenated extreme learning
machine (HLConcELM).

Given an architectural vector $\mbs M$, an activation function $\sigma(\cdot)$,
and the randomly assigned values for the hidden-layer coefficients $\bm\theta$,
let HLConcELM($\mbs M,\sigma,\bm\theta$) denote the
associated hidden-layer concatenated ELM. For a given domain
$\Omega\subset\mbb R^d$, a vector $\mbs M=(M_0,M_1,\dots,M_L)$ with $M_0=d$ and $M_L=1$,
and given $\bm\theta\in\mathbb{R}^{N_h}$ and $\sigma$, we define
\begin{equation}\label{eq_12}
  U(\Omega,\mbs M,\sigma,\bm\theta) = \left\{
  u(\bm\theta,\bm\beta,\mbs x)\ |\ u(\bm\theta,\bm\beta,\mbs x)\
  \text{is the output of HLConcFNN}(\mbs M,\sigma),\
  \mbs x\in\Omega,\ \bm\beta\in\mbb R^{N_c}
  \right\}
\end{equation}
as the set of all the possible output fields of HLConcELM($\mbs M,\sigma,\bm\theta$)
on $\Omega$,
where $N_c$ denotes the total number of the output-layer coefficients.
Hereafter we refer to $U(\Omega,\mbs M,\sigma,\bm\theta)$ as the representation
capacity of the HLConcELM($\mbs M,\sigma,\bm\theta$).
Note that $U(\Omega,\mbs M,\sigma,\bm\theta)$ forms a linear space.

Analogous to Theorem~\ref{thm_1}, when one appends hidden layers
to a given network architecture,
the representation capacity of the HLConcELM associated with the resultant
architecture will be at least as large as that associated with the original one,
on condition that the random hidden-layer coefficients of the new HLConcELM
are set appropriately. On the other hand,
if one adds extra nodes to a hidden layer (other than the last one)
of a given architecture, there is no  analogous result to Theorem~\ref{thm_2}
for HLConcELM, because the hidden-layer coefficients in ELM are
randomly set. Specifically, we have the following result.
\begin{theorem}\label{thm_3}
  Given an architectural vector $\mbs M_1=(m_0,m_1,\dots,m_{L-1},m_L)$
  with $m_L=1$, define a new vector $\mbs M_2=(m_0,m_1,\dots,m_{L-1},n,m_L)$, where
  $n\geqslant 1$ is an integer.
  Let $\bm\theta\in\mbb R^{N_{h1}}$ and $\bm\vartheta\in\mbb R^{N_{h2}}$
  denote two random vectors, with the relation
  $\bm\vartheta[1:N_{h1}]=\bm\theta[1:N_{h1}]$, where
  $N_{h1}=\sum_{i=1}^{L-1}(m_{i-1}+1)m_i$ and
  $N_{h2}=N_{h1}+(m_{L-1}+1)n$.
  For a given domain $\Omega\subset\mbb R^{m_0}$ and an activation function
  $\sigma(\cdot)$, the following relation holds
  \begin{equation}\label{eq_13}
    U(\Omega,\mbs M_1,\sigma,\bm\theta) \subseteq U(\Omega,\mbs M_2,\sigma,\bm\vartheta),
  \end{equation}
  where $U$ is defined by~\eqref{eq_12}.

\end{theorem}
\noindent By $\bm\vartheta[1:N_{h1}]=\bm\theta[1:N_{h1}]$ we mean that
the first $N_{h1}$ entries of $\bm\vartheta$ and $\bm\theta$ are the same.
Because of this condition, the random bases for $U(\Omega,\mbs M_2,\sigma,\bm\vartheta)$
would contain those bases for $U(\Omega,\mbs M_1,\sigma,\bm\theta)$,
giving rise to the relation~\eqref{eq_13}.
It should be noted that conventional ELMs lack a comparable property
as expressed by the relation~\eqref{eq_13}


In the current paper
we set the random hidden-layer coefficients $\bm\theta$ in HLConcELM
in the following fashion.
Given an architectural vector $\mbs M=(m_0,m_1,\dots,m_{L-1},m_L)$,
let $\bm\xi\in\mbb R^{N_h}$ be a random vector generated
on the interval $[-1,1]$ from a uniform distribution,
where $N_h=\sum_{i=1}^{L-1}(m_{i-1}+1)m_i$.
Once generated, $\bm\xi$ will be
fixed throughout the computation for the given architecture $\mbs M$.
We next partition $\bm\xi$ into $(L-1)$ sub-vectors,
$\bm\xi=(\bm\xi_1,\bm\xi_2,\dots,\bm\xi_{L-1})$, with
$\bm\xi_{i}$ having a dimension $(m_{i-1}+1)m_i$ for $1\leqslant i\leqslant L-1$.
Let $\bm R=(R_1,R_2,\dots,R_{L-1})$ denote $(L-1)$ constants.
We then set $\bm\theta$ in HLConcELM for the given architecture $\mbs M$ to
\begin{equation}\label{eq_14}
  \bm\theta(\mbs M, \mbs R, \bm\xi) =
  \text{\bf flatten}\left[R_1\bm\xi_1,R_2\bm\xi_2,\dots,R_{L-1}\bm\xi_{L-1}\right],
\end{equation}
where ``{\bf flatten}'' concatenates the list of vectors into
a single vector.


Hereafter we refer to the above vector $\mbs R=(R_1,\dots,R_{L-1})$
as the hidden magnitude vector for the network architecture $\mbs M$.
When assigning random hidden-layer coefficients as described above,
we have essentially set the weight/bias coefficients in the $i$-th hidden layer
to uniform random values generated on the interval $[-|R_i|,|R_i|]$,
where $R_i$ is the $i$-th component of $\mbs R$,
for $1\leqslant i\leqslant L-1$.
The constant $|R_i|$ denotes the maximum magnitude of
the random coefficients for the $i$-th hidden layer.

The constants $R_i$ ($1\leqslant i\leqslant L-1$) are
the hyperparameters of the HLConcELM. The idea of generating
random coefficients for different hidden layers with different maximum
magnitudes is first studied in~\cite{DongY2021} for conventional
feedforward neural networks, and a method based on
the differential evolution algorithm is developed therein
for computing the optimal values of those magnitudes.
In the current work, for a given PDE problem,
we use the method of~\cite{DongY2021} to compute
the optimal (or near-optimal) hidden magnitude vector $\mbs R^*$,
and employ $\mbs R=\mbs R^*$ in HLConcELM
for the simulations.

Hereafter we use HLConcELM($\mbs M,\sigma,\mbs R,\bm\xi$)
to denote the hidden-layer concatenated ELM
characterized by the architectural vector $\mbs M$,
the activation function $\sigma(\cdot)$,
the randomly-assigned but fixed vector $\bm\xi$ on $[-1,1]$,
and the hidden magnitude vector $\mbs R$.
According to Theorem~\ref{thm_3}, when additional hidden layers
are appended to a given HLConcELM($\mbs M,\sigma,\mbs R,\bm\xi$),
the representation capacity of the resultant HLConcELM
will not be smaller than that of the original one,
if the vectors $\mbs R$ and $\bm\xi$ of the resultant network
are set appropriately.


\subsection{Solving linear/nonlinear PDEs with Hidden-Layer Concatenated ELM}
\label{sec:pde}


We next discuss how to use the hidden-layer concatenated
ELM for the numerical solution of PDEs.
Consider a domain $\Omega\subset\mbb R^d$ and 
the following boundary value problem on this domain,
\begin{subequations}\label{eq_15}
\begin{align}
  &
  \mathcal{L}u + F(u) = f(\mbs x),  \label{eq_15a} \\
  &
  Bu + G(u) = g(\mbs x), \quad \mbs x\in\partial\Omega. \label{eq_15b}
\end{align}
\end{subequations}
In these equations
$u(\mbs x)$ is the field function to be solved for.
$\mathcal{L}$ is a linear differential operator. $F(u)$ is a nonlinear
operator acting on $u$ and also possibly its derivatives.
Equation~\eqref{eq_15b} represents the boundary conditions, where
$B$ is a linear differential or algebraic operator.
The boundary condition may
possibly contain some nonlinear operator $G(u)$ acting on $u$
and also possibly its derivatives.
If  both $F(u)$ and $G(u)$ are absent
the problem becomes linear.
We assume that this problem is well posed.

In addition, we assume that $\mathcal{L}$ may possibly
include time derivatives (e.g.~$\frac{\partial}{\partial t}$,
$\frac{\partial^2}{\partial t^2}$).
In this case the problem~\eqref{eq_15} becomes time-dependent, and
we will treat the time $t$ in the same way as the spatial coordinate $\mbs x$
and consider $t$ as the last dimension in $d$ dimensions.
We require that the equation~\eqref{eq_15b} should include appropriate
initial condition(s) for such a case. So
the problem~\eqref{eq_15} may refer to time-dependent cases, which will not
be distinguished in the following discussions.

We represent the solution field $u(\mbs x)$ by a hidden-layer
concatenated ELM from the previous subsection. 
Consider a network architecture given by $\mbs M=(m_0,m_1,\dots,m_L)$,
where $m_0=d$ and $m_L=1$, and an activation function $\sigma(\cdot)$.
We use the HLConcFNN($\mbs M,\sigma$) to represent
the solution field $u(\mbs x)$ (see Figure~\ref{fg_1}(b)). Here the $d$
input nodes  represent $\mbs x$, and the single output node
represents $u(\mbs x)$.
The activation function $\sigma$ is applied to all the hidden nodes
in the network. As noted before, we require that
the output layer contain no activation function 
and have zero bias.
With a given hidden magnitude vector $\mbs R=(R_1,R_2,\dots,R_{L-1})$ and a
randomly generated vector $\bm\xi\in\mbb R^{N_h}$  on $[-1,1]$,
where $N_h=\sum_{i=1}^{L-1}(m_{i-1}+1)m_i$, we
set and fix the random hidden-layer coefficients according to
equation~\eqref{eq_14}.

Under these settings,
the output field of the neural network is given by equation~\eqref{eq_2}.
Substituting this expression for $u(\mbs x)$ into the system~\eqref{eq_15},
we have
\begin{subequations}
  \begin{align}
    &
    \sum_{i=1}^{L-1}\sum_{j=1}^{m_i}\beta_{ij}\left[\mathcal{L}\phi_j^{(i)}(\mbs x) \right]
    + F\left(\sum_{i=1}^{L-1}\sum_{j=1}^{m_i}\beta_{ij}\phi_j^{(i)}(\mbs x) \right) = f(\mbs x)
    \label{eq_16a} \\
    &
    \sum_{i=1}^{L-1}\sum_{j=1}^{m_i}\beta_{ij}\left[B\phi_j^{(i)}(\mbs x) \right]
    + G\left(\sum_{i=1}^{L-1}\sum_{j=1}^{m_i}\beta_{ij}\phi_j^{(i)}(\mbs x) \right) = g(\mbs x),
    \quad \text{on}\ \partial\Omega,
    \label{eq_16b}
  \end{align}
\end{subequations}
where $\phi_j^{(i)}(\mbs x)$ ($1\leqslant i\leqslant L-1$, $1\leqslant j\leqslant m_i$)
denotes the output field of the $j$-th node in the $i$-th hidden layer,
and $\beta_{ij}$ ($1\leqslant i\leqslant L-1$, $1\leqslant j\leqslant m_i$)
are the weight coefficients in the output layer of the HLConcELM.
It should be noted that, since the hidden-layer coefficients are randomly set but fixed,
$\phi_j^{(i)}(\mbs x)$ are random but fixed functions.
The coefficients $\beta_{ij}$ are the trainable parameters in HLConcELM.


We next choose a set of $Q$ ($Q\geqslant 1$) points on $\Omega$,
referred to as the collocation points,
which can be regular grid points,
random points,  or chosen based on
some other distribution. Among these points we assume
that $Q_b$ ($1\leqslant Q_b\leqslant Q-1$) points
reside on the boundary $\partial\Omega$
and the rest are from the interior of $\Omega$.
Let $\mbb{X}$ denote the set of all the collocation points, and $\mbb{X}_b$
denote the set of the boundary collocation points.

We enforce the equation~\eqref{eq_16a} on all the
collocation points from $\mbb{X}$,
and enforce the
equation~\eqref{eq_16b} on all the boundary collocation points
from $\mbb{X}_b$.
This leads to
\begin{subequations}\label{eq_17}
  \begin{align}
    &
    \sum_{i=1}^{L-1}\sum_{j=1}^{m_i}\beta_{ij}\left[\mathcal{L}\phi_j^{(i)}(\mbs x_p) \right]
    + F\left(\sum_{i=1}^{L-1}\sum_{j=1}^{m_i}\beta_{ij}\phi_j^{(i)}(\mbs x_p) \right) = f(\mbs x_p),
    \quad \mbs x_p\in\mbb X, \quad
    1\leqslant p\leqslant Q;
    \label{eq_17a} \\
    &
    \sum_{i=1}^{L-1}\sum_{j=1}^{m_i}\beta_{ij}\left[B\phi_j^{(i)}(\mbs x_q) \right]
    + G\left(\sum_{i=1}^{L-1}\sum_{j=1}^{m_i}\beta_{ij}\phi_j^{(i)}(\mbs x_q) \right) = g(\mbs x_q),
    \quad \mbs x_q\in\mbb X_b, \quad 1\leqslant q\leqslant Q_b.
    \label{eq_17b}
  \end{align}
\end{subequations}
This is a system of $(Q+Q_b)$ nonlinear algebraic equations about
$N_c=\sum_{i=1}^{L-1}m_i$ unknowns, $\beta_{ij}$.
The differential operators involved in these equations, such as
$\mathcal{L}\phi_j^{(i)}(\mbs x_p)$, $B\phi_j^{(i)}(\mbs x_q)$,
$F(u(\mbs x_p))$ and $G(u(\mbs x_q))$ where
$u(\mbs x_p)=\sum_{i=1}^{L-1}\sum_{j=1}^{m_i}\beta_{ij}\phi_j^{(i)}(\mbs x_p)$,
can be computed by auto-differentiations of the neural network.

The system~\eqref{eq_17} is a rectangular system, in which
the number of equations and the number of unknowns are not the same.
We seek a least squares solution to this system.
This is a nonlinear least squares problem, and it can be solved
by the Gauss-Newton method together with the trust region
strategy~\cite{NocedalW2006}.
Several quality implementations of the Gauss-Newton method are available
from the scientific libraries.
In this work we employ the Gauss-Newton implementation together with
a trust region reflective algorithm~\cite{BranchCL1999,ByrdSS1988}
from the scipy package in Python (scipy.optimize.least\_squares)
to solve this problem. We refer to the method implemented
in this scipy routine as
the nonlinear least squares method in this paper.

The nonlinear least squares method requires
two procedures for solving the system~\eqref{eq_17},
one for computing the residual of this system and the other
for computing the Jacobian matrix for a given arbitrary $\beta_{ij}$
($1\leqslant i\leqslant L-1$, $1\leqslant j\leqslant m_i$).
For a given arbitrary $\bm\beta$ (see~\eqref{eq_3}),
the residual $\mbs r(\bm\beta)\in\mbb R^{Q+Q_b}$ is given by
\begin{equation}\label{eq_18}
  \left\{
  \begin{split}
    &
  \mbs r(\bm\beta) = \left(\mbs r_1(\bm\beta),\mbs r_2(\bm\beta)\right), \quad
  \mbs r_1(\bm\beta) = (r_{11},r_{12},\dots,r_{1Q}), \quad
  \mbs r_2(\bm\beta) = (r_{21}, r_{22},\dots, r_{2Q_b}); \\
  &
  r_{1p} = \sum_{i=1}^{L-1}\sum_{j=1}^{m_i}\beta_{ij}\left[\mathcal{L}\phi_j^{(i)}(\mbs x_p) \right]
    + F\left(\sum_{i=1}^{L-1}\sum_{j=1}^{m_i}\beta_{ij}\phi_j^{(i)}(\mbs x_p) \right) - f(\mbs x_p),
    \quad \mbs x_p\in\mbb X, \quad
    1\leqslant p\leqslant Q; \\
    &
    r_{2q} = \sum_{i=1}^{L-1}\sum_{j=1}^{m_i}\beta_{ij}\left[B\phi_j^{(i)}(\mbs x_q) \right]
    + G\left(\sum_{i=1}^{L-1}\sum_{j=1}^{m_i}\beta_{ij}\phi_j^{(i)}(\mbs x_q) \right) - g(\mbs x_q),
    \quad \mbs x_q\in\mbb X_b, \quad 1\leqslant q\leqslant Q_b.
  \end{split}
  \right.
\end{equation}
The Jacobian matrix $\frac{\partial\mbs r}{\partial\bm\beta}\in\mbb R^{(Q+Q_b)\times N_c}$
is given by,
\begin{equation}\label{eq_19}
  \left\{
  \begin{split}
    &
    \frac{\partial\mbs r}{\partial\bm\beta}=\begin{bmatrix}
    \frac{\partial\mbs r_1}{\partial\bm\beta} \\ \frac{\partial\mbs r_2}{\partial\bm\beta}
    \end{bmatrix}_{(Q+Q_b)\times N_c}, \quad
    \frac{\partial\mbs r_1}{\partial\bm\beta} = \begin{bmatrix}
      \frac{\partial r_{1p}}{\partial\beta_{ij}}
    \end{bmatrix}_{Q\times N_c}, \quad
    \frac{\partial\mbs r_2}{\partial\bm\beta} = \begin{bmatrix}
      \frac{\partial r_{2q}}{\partial\beta_{ij}}
    \end{bmatrix}_{Q_b\times N_c}; \\
    &
    \frac{\partial r_{1p}}{\partial\beta_{ij}} = \mathcal{L}\phi_j^{(i)}(\mbs x_{p})
    + F^{\prime}(u(\mbs x_p))\phi_j^{(i)}(\mbs x_p), \quad
    \mbs x_p\in\mbb X, \quad 1\leqslant p\leqslant Q, \quad
    1\leqslant i\leqslant L-1,\ \ 1\leqslant j\leqslant m_i;
    \\
    &
    \frac{\partial r_{2q}}{\partial\beta_{ij}} = B\phi_j^{(i)}(\mbs x_q)
    + G'(u(\mbs x_q))\phi_j^{(i)}(\mbs x_q), \quad
    \mbs x_q\in\mbb X_b, \quad 1\leqslant p\leqslant Q_b, \quad
    1\leqslant i\leqslant L-1,\ \ 1\leqslant j\leqslant m_i;
  \end{split}
  \right.
\end{equation}
where $u(\mbs x_p)$ is computed based on equation~\eqref{eq_2}.
The $F'(u)$ and $G'(u)$ terms denote the derivatives with respect to $u$, and
may represent the effect of an operator. For example,
the nonlinear function $F(u) = u\frac{\partial u}{\partial x}$ (as
in the Burgers' equation) leads to
$F'(u)\phi = \frac{\partial u}{\partial x}\phi + u\frac{\partial\phi}{\partial x}$.


Therefore,
to solve the problem~\eqref{eq_15} by HLConcELM,
the input training data (denoted by $\mbs X$) to the neural
network is a $Q\times d$ matrix,
consisting of the coordinates
of all the collocation points,
$
\mbs X= \begin{bmatrix}\mbs x_p  \end{bmatrix}_{Q\times d}
$
(for all $\mbs x_p\in\mbb X$).
The output data (denoted by $\mbs U$) of the neural
network is a $Q\times m_L$ matrix,
representing the field solution $u(\mbs x)$ on the collocation
points,
$
\mbs U = \begin{bmatrix} u(\mbs x_p)  \end{bmatrix}_{Q\times m_L}.
$
The output data of the logical concatenation layer
(denoted by $\bm\Psi$) of the HLConcELM
is a $Q\times N_c$ matrix given by
$
\bm\Psi = \begin{bmatrix}\bm\Phi(\bm\theta,\mbs x_p)  \end{bmatrix}_{Q\times N_c}.
$
It represents the output fields of the all the hidden nodes
on all the collocation points.
Here $N_c$ denotes the total number of hidden nodes in the network, and
$\bm\theta$ denotes the random hidden-layer coefficients given by~\eqref{eq_14}.
The relation~\eqref{eq_7} is translated into,
in terms of the neural-network data,
\begin{equation}\label{eq_20}
  \mbs U = \bm\Psi\bm\beta^T,
\end{equation}
where $\bm\beta$ denotes the output-layer coefficients given
by~\eqref{eq_3}.

\begin{remark}\label{rem_1}
  The output data of the logical concatenation layer
  $\bm\Psi$ can be computed by a forward
  evaluation of the neural network (for up to
  the logical concatenation layer) on the input data $\mbs X$.
  In our implementation we have created a Keras sub-model
  with the input layer as its input and the logical
  concatenation layer as its output.
  By evaluating this sub-model on the input data we can attain
  the output data for all the hidden nodes on the collocation points.
  The first and higher derivatives 
  of $\bm\Psi$ with respect to $\mbs X$ are computed by a
  forward-mode auto-differentiation, implemented by the ``ForwardAccumulator''
  in the Tensorflow library. This forward-mode auto-differentiation
  is crucial for the computational performance, because the total number of
  hidden nodes ($N_c$) in HLConcELM is typically much larger
  than the number of input nodes ($d$).
  The differential operators on the output fields of the hidden nodes
  involved in~\eqref{eq_18} and~\eqref{eq_19}, such as
  $\mathcal{L}\phi_j^{(i)}(\mbs x_p)$ ($\mbs x_p\in\mbb X$),
  $B\phi_j^{(i)}(\mbs x_q)$ ($\mbs x_q\in\mbb X_b$) and
  $F'(u(\mbs x_p))\phi_j^{(i)}(\mbs x_p)$, can be computed
  based on or extracted from $\bm\Psi$ and its derivatives with
  respect to $\mbs X$.
  Once $\bm\Psi$ is attained, for a given $\bm\beta$,
  the output data of the neural network can be computed by~\eqref{eq_20},
  which provides the $u(\mbs x_p)$
  ($\mbs x_p\in\mbb X$ or $\mbs x_p\in\mbb X_b$) for computing
  the terms $F'(u(\mbs x_p))$ and $G'(u(\mbs x_q))$ in~\eqref{eq_18}
  and~\eqref{eq_19}.

\end{remark}

\begin{remark}\label{rem_2}
  If the boundary value problem~\eqref{eq_15} is linear, i.e.~in
  the absence of the terms $F(u)$ and $G(u)$,
  the resultant system~\eqref{eq_17} is a linear algebraic
  system of $(Q+Q_b)$ equations about $N_c$ unknowns of
  the parameters $\beta_{ij}$.
  In this case we use the linear least squares method
  to solve this system to compute a least squares solution
  for $\beta_{ij}$. In our implementation we employ
  the linear least squares routine from scipy (scipy.linalg.lstsq),
  which in turn employs the linear
  least squares implementation from the LAPACK
  library.
  
\end{remark}

\begin{remark}\label{rem_3}
  If the problem~\eqref{eq_15} is time-dependent, for longer-time
  or long-time simulations, we employ the block marching
  scheme from~\cite{DongL2021} together with the HLConcELM
  for its computation. The temporal dimension, which can be
  potentially large in this case, is first divided into
  a number of windows (referred as time blocks), so that each
  time block is of a moderate size. The problem on each time block
  is solved by HLConcELM individually and successively.
  After one time block is computed, the solution evaluated at
  the last time instant, possibly together with its derivatives,
  is used as the initial condition for computing the time block
  that follows. We refer the reader to~\cite{DongL2021}
  for more detailed discussions of the block time marching scheme.
  
\end{remark}

\begin{remark}\label{rem_4}
  HLConcFNNs can be used together with the locELM
  (local extreme learning machine)
  method~\cite{DongL2021}
  and domain decomposition for solving PDE problems.
  In this case, we employ a HLConcELM 
  for the local neural network on each sub-domain, and the algorithm
  for computing the PDE solution is essentially the same. The only difference
  lies in that in the system~\eqref{eq_17} one needs to additionally include the
  $C^k$ continuity conditions on those collocation points that
  reside on the sub-domain boundaries.
  The residuals in~\eqref{eq_18} and
  the Jacobian matrix in~\eqref{eq_19} need to be modified accordingly to account for
  these additional equations from the $C^k$ continuity conditions.
  We refer the reader to~\cite{DongL2021} for detailed
  discussions of these aspects.
  For convenience of presentation,
  hereafter we will refer to the locELM method based on HLConcFNNs
  as the locHLConcELM method (local hidden-layer
  concatenated ELM).
  
\end{remark}

\begin{remark}\label{rem_5}
  For a given problem, the optimal or near-optimal value $\mbs R^*$
  for the hidden magnitude 
  vector $\mbs R$
  can be computed by the method from~\cite{DongY2021} based on the
  differential evolution algorithm.
  For all the test problems in Section~\ref{sec:tests}, we employ
  $\mbs R=\mbs R^*$ computed based on the method of~\cite{DongY2021}
  in the HLConcELM simulations. 

\end{remark}


%% file: Test.tex
\section{Numerical Benchmarks}
\label{sec:tests}

In this section we employ several benchmark problems
in two dimensions (2D) or in one spatial dimension (1D) plus time
to test the performance of the HLConcELM method
for solving linear and nonlinear PDEs.
We show that this method can produce highly accurate 
results  when the network  architecture has
a narrow last hidden layer.
In contrast, the conventional ELM
method in this case utterly loses accuracy.


The HLConcELM method is implemented in Python based on
the Tensorflow and Keras libraries. The linear and nonlinear least squares
methods employed in HLConcELM are based on the implementations in the scipy package
(scipy.linalg.lstsq and scipy.optimize.least\_squares), as discussed before.
The differential operators on the hidden-layer data
(see equations~\eqref{eq_17a}--\eqref{eq_17b}) are computed by a forward-mode
auto-differentiation, as stated in Remark~\ref{rem_1}.
In all the numerical tests of this section
we employ the Gaussian activation function
$\sigma(x)=e^{-x^2}$ for all the hidden nodes,
while the output layer is linear and has zero bias.


The ELM errors reported in the following subsections are computed
as follows. We have considered regular rectangular
domains for simplicity in the current paper.
For a given architecture we train the HLConcELM
network on $Q=Q_1\times Q_1$ uniform collocation points (i.e.~regular grid points)
by the linear or nonlinear least squares method,
with $Q_1$ uniform points in each direction of the 2D domain
or the spatial-temporal domain.
After the network is trained, we evaluate the neural network
on a finer set of $Q_{eval}=Q_2\times Q_2$ uniform grid points, with $Q_2$
much larger than $Q_1$, to attain the HLConcELM solution data.
We evaluate the exact solution  to the problem, if available, on
the same set of $Q_{eval}$ grid points. Then we compare the HLConcELM
solution data and the exact solution data on the $Q_2\times Q_2$ grid points
to compute the maximum ($l^{\infty}$) and root-mean-squares (rms, or ${l}^2$)
errors. We refer to the errors computed above as the HLConcELM errors
associated with the given network architecture and the
$Q=Q_1\times Q_1$ training collocation points.
When $Q_1$ is varied in a range for the convergence tests, we have made sure
that $Q_2$ is much larger than the largest $Q_1$ in the prescribed range.
When the block time marching scheme is used for longer-time simulations
together with HLConcELM (see Remark~\ref{rem_3}), the $Q=Q_1\times Q_1$
and $Q_{eval}=Q_2\times Q_2$ points above refer to the points
in each time block.
When the locHLConcELM method together with domain decomposition is used
to solve a problem (see Remark~\ref{rem_4}),
the $Q$ and $Q_{eval}$ points refer to the points
in each sub-domain.
In the current paper we employ a fixed $Q_{eval}=101\times 101$ (i.e.~$Q_2=101$)
when evaluating the neural network and computing the HLConcELM
errors for all the test problems in this section.


As in our previous works~\cite{DongL2021,DongL2021bip},
we employ a fixed seed for the random number generator in the Tensorflow
library in order to make the reported numerical results herein
exactly reproducible. While the seed value is different for the test problems
in different subsections, it has been fixed to a particular value
for the numerical tests within each subsection.
Specifically, the seed to the random number generator is $10$ in
Section~\ref{sec:poisson}, $50$ in Section~\ref{sec:nonhelm}, and
$100$ in Sections~\ref{sec:advec}, \ref{sec:burger} and~\ref{sec:kdv}.


In comparisons with the conventional ELM method~\cite{DongL2021} in
the following subsections, all the hidden-layer coefficients in
conventional ELM are assigned (and fixed) to uniform random values
generated on the interval $[-R_m,R_m]$, with $R_m=R_{m0}$,
where $R_{m0}$ is the optimal $R_m$ computed by the method of~\cite{DongY2021}
based on the differential evolution algorithm.


\subsection{Variable-Coefficient Poisson Equation}
\label{sec:poisson}

\begin{figure}
  \centerline{
    \includegraphics[width=2.in]{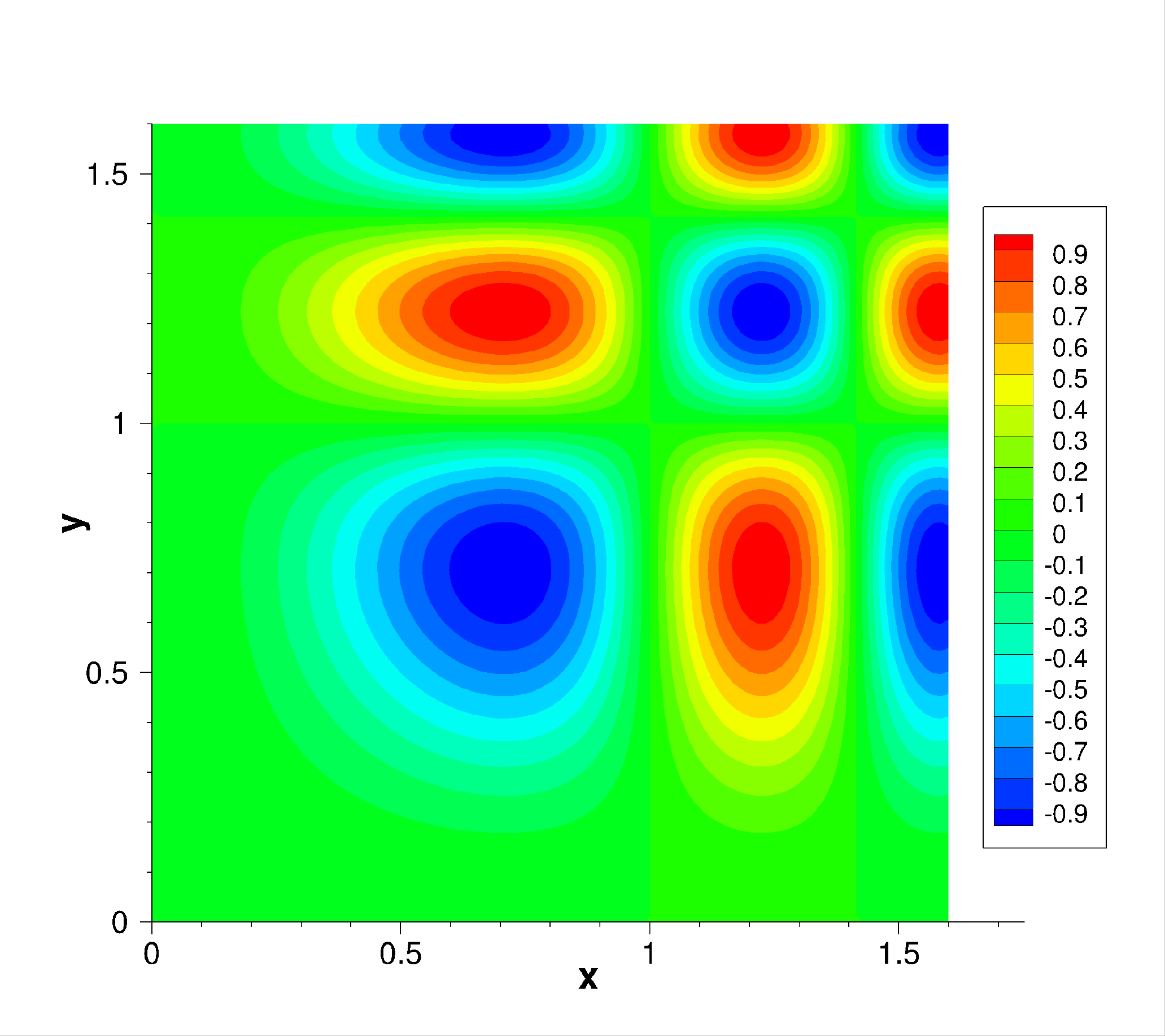}(a)
    \includegraphics[width=2.in]{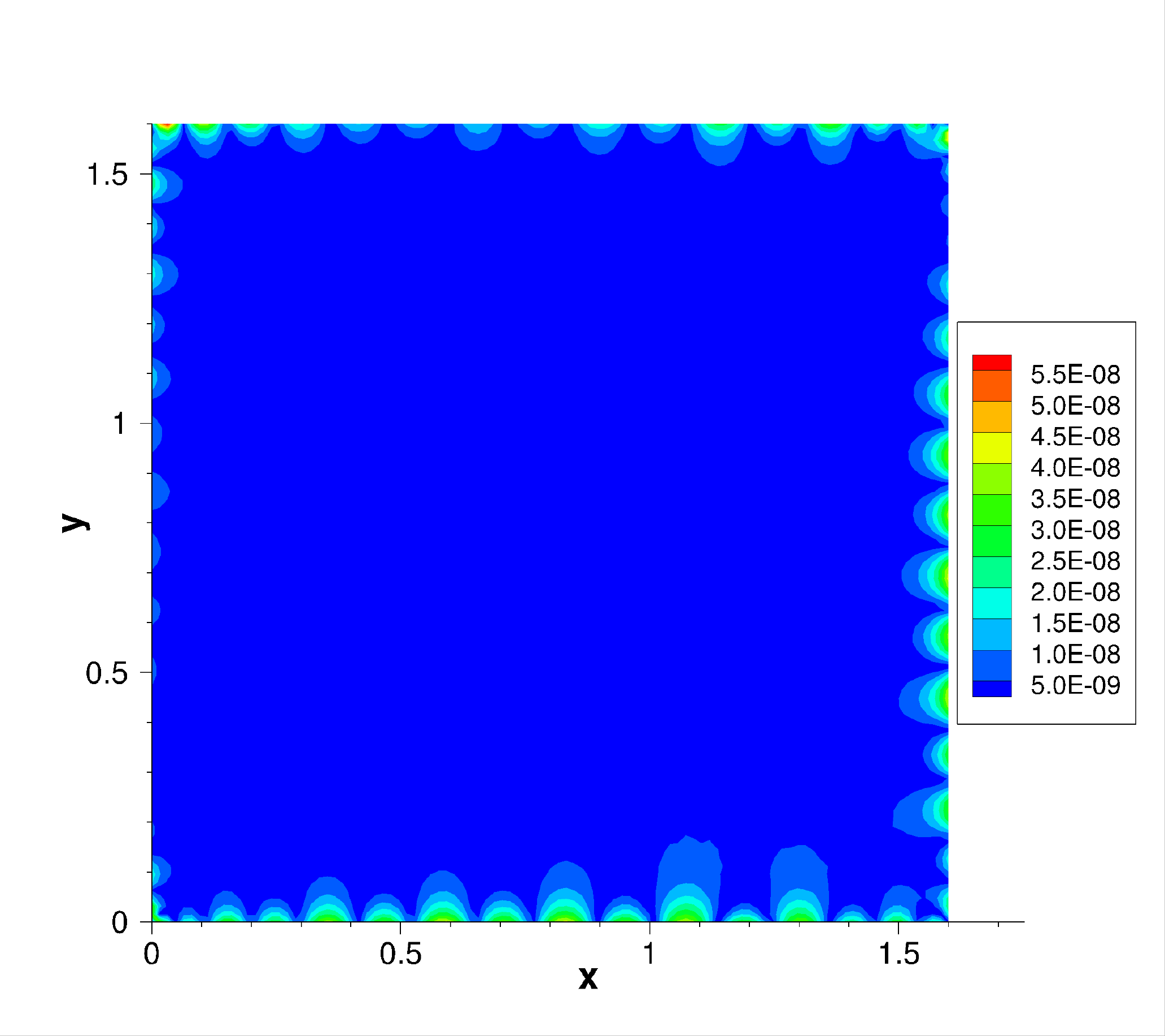}(b)
    \includegraphics[width=2.in]{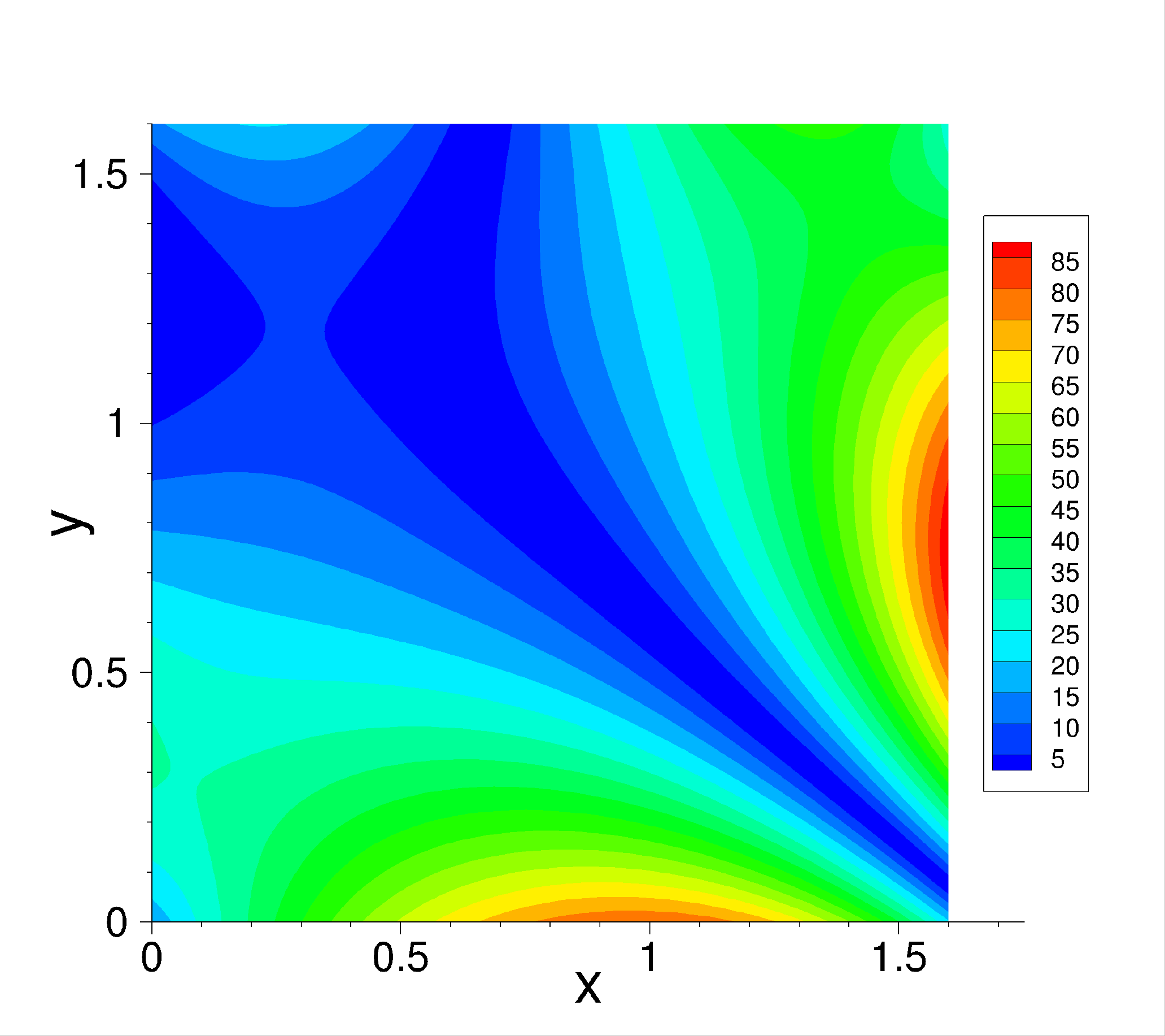}(c)
  }
  \caption{Variable-coefficient Poisson equation:
    distributions of (a) the exact solution,
    (b) the absolute error of the HLConcELM solution, and (c) the absolute
    error of the conventional ELM solution.
    In (b,c), network architecture $\mbs M=[2, 800, 50, 1]$,
    Gaussian activation function,
    $Q=35\times 35$ uniform collocation points.
    $\mbs R=(3.0, 0.005)$ for HLConcELM in (b). $R_m=R_{m0}=0.35$
    for conventional ELM in (c), where $R_{m0}$ is the optimal
    $R_m$ computed using the method of~\cite{DongY2021}.
  }
  \label{fg_3}
\end{figure}

The first numerical test involves  the 2D Poisson equation with
a variable coefficient field. Consider the 2D domain
$\Omega = [0,1.6]\times [0,1.6]$ and the the following boundary value
problem on $\Omega$,
\begin{subequations}\label{eq_21}
  \begin{align}
    &
    \frac{\partial}{\partial x}\left(a(x,y)\frac{\partial u}{\partial x}  \right)
    + \frac{\partial}{\partial y}\left(a(x,y)\frac{\partial u}{\partial y}  \right)
    = f(x,y), \label{eq_21a} \\
    &
    u(0,y) = g_1(y), \quad
    u(1.6,y) = g_2(y), \quad
    u(x,0) = h_1(x), \quad
    u(x,1.6) = h_2(x), \label{eq_21b}
  \end{align}
\end{subequations}
where $(x,y)$ are the spatial coordinates,
$u(x,y)$ is the field function to be solved for, $f(x,y)$ is a prescribed
source term, $a(x,y)$ is the coefficient field given by
$a(x,y) = 2 + \sin(x+y)$, and $g_1$, $g_2$, $h_1$ and $h_2$
are prescribed boundary distributions.
We choose the source term $f$ and the boundary data $g_i$ and $h_i$ ($i=1,2$)
appropriately such that the following function satisfies the system~\eqref{eq_21},
\begin{equation}\label{eq_22}
  u(x,y) = -\sin(\pi x^2)\sin(\pi y^2).
\end{equation}
The distribution of this exact solution in the $xy$ plane
is illustrated by Figure~\ref{fg_3}(a).


We employ the HLConcELM method from Section~\ref{sec:method} to solve
the system~\eqref{eq_21}. Let the vector $\mbs M=[2, m_1,\dots,m_{L-1},1]$
denote the architecture of
the HLConcELM, where the two input nodes
represent the coordinates $x$ and $y$ and the single output node represents
the solution $u$. We employ the Gaussian activation function for all the hidden
nodes, as stated at the beginning of Section~\ref{sec:tests}.
The output layer is linear  and has
no bias. The number of hidden layers and the number of hidden nodes
are varied, and the specific
architure will be given below when discussing the results.

We employ a uniform set of $Q=Q_1\times Q_1$ grid points on the domain $\Omega$,
with $Q_1$ uniform points on each side of the boundary,
as the collocation points
for training the neural network. $Q_1$ is varied in the tests.
As discussed earlier, after the neural network is trained,
we evaluate the neural network on
another finer set of $Q_{eval}=Q_2\times Q_2$, with $Q_2=101$,
uniform grid points on $\Omega$ 
and compute the HLConcELM errors.

Figures \ref{fg_3}(b) and (c) show a comparison of the absolute-error distributions
in the $xy$ plane of the HLConcELM solution and the conventional ELM
solution obtained using a neural network with a narrow last hidden layer.
Note that in conventional ELM the usual feedforward neural network
has been employed (see Figure~\ref{fg_1}(a)).
For both HLConcELM and conventional ELM
we employ here a neural network with the  architecture
$\mbs M=[2, 800, 50, 1]$ and a uniform set of $Q=35\times 35$ collocation
points for the network training. With HLConcELM,
for setting the random hidden-layer coefficients, we employ
a hidden magnitude vector $\mbs R=(3.0,0.005)$, which is close to
the optimum $\mbs R^*$ obtained based on the method of~\cite{DongY2021}.
With conventional ELM, we set the hidden-layer coefficients
to uniform random values generated on the interval $[-R_m,R_m]$
with $R_m=R_{m0}$, where $R_{m0}=0.35$ is the optimal $R_m$ obtained
using the method of~\cite{DongY2021} for this case.
Because the last hidden layer 
is quite narrow (with $50$ nodes), we observe that
the result of the conventional ELM
exhibits no accuracy, with a maximum error around $85$ in the domain.
In contrast, the HLConcELM method produces highly accurate
results, with the maximum error on the order of $10^{-8}$ in the domain.

\begin{figure}
  \centerline{
    \includegraphics[width=2.in]{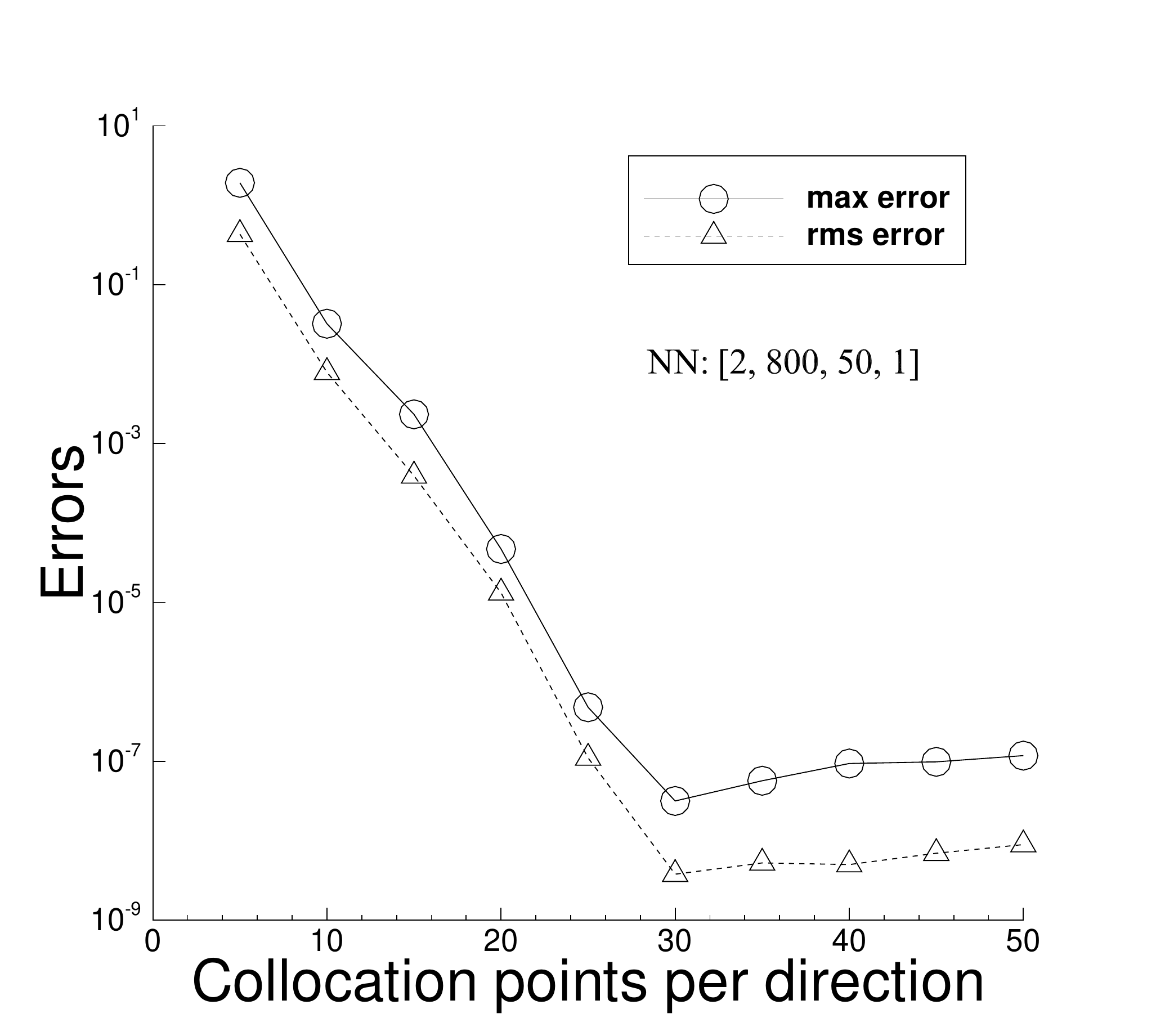}(a)
    \includegraphics[width=2.in]{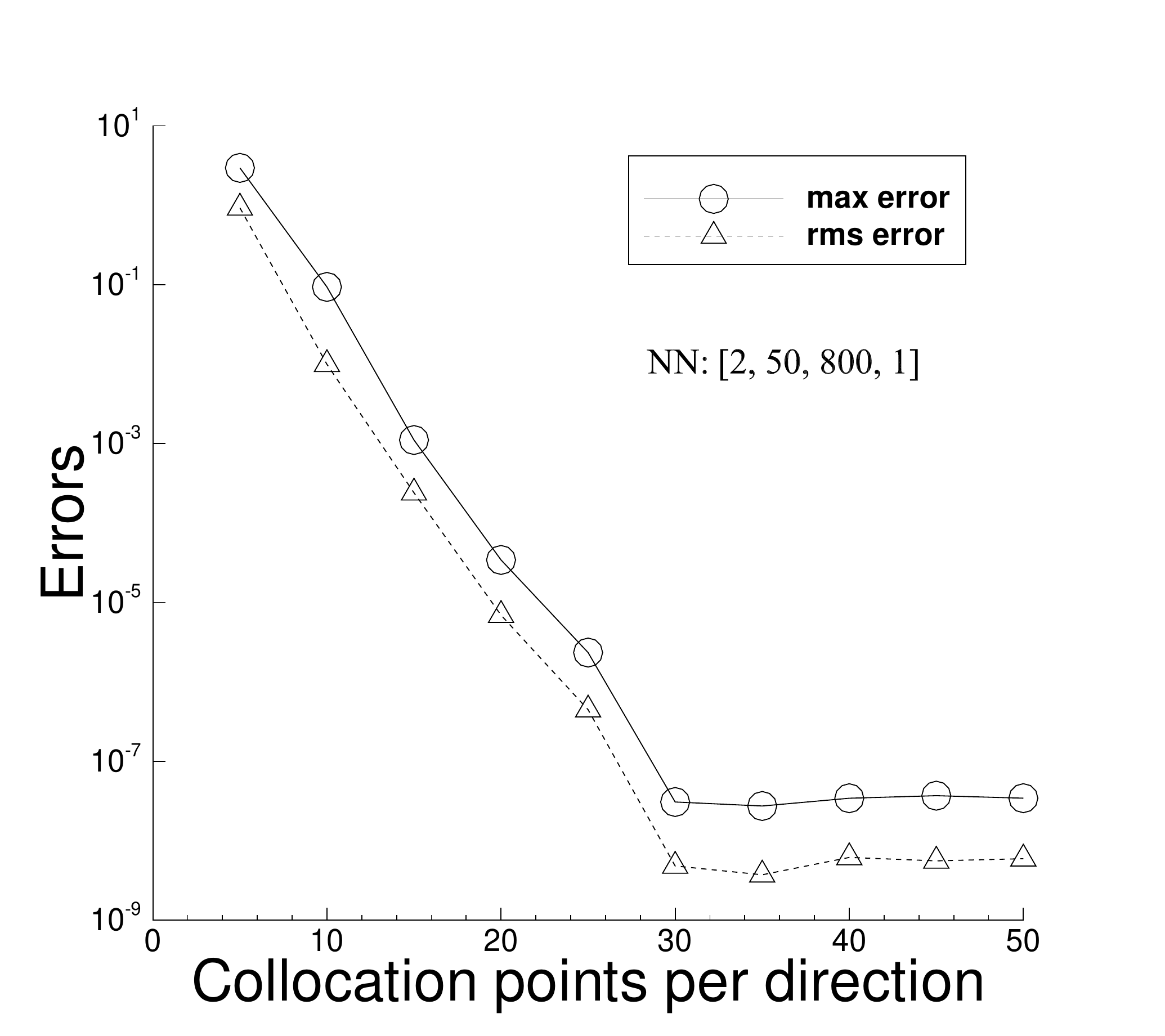}(b)
  }
  \caption{Variable-coefficient Poisson equation:
    the maximum/rms errors of the HLConcELM solution
    versus the number of collocation points per direction
    obtained with a network architecture of (a) $[2, 800, 50, 1]$ and
    (b) $[2, 50, 800, 1]$.
    $\mbs R=(3.0, 0.005)$ in (a), and $\mbs R=(0.68, 0.82)$ in (b).
  }
  \label{fg_4}
\end{figure}

Figure~\ref{fg_4} illustrates the convergence behavior of the HLConcELM
solution with respect to the number of collocation points in
the network training.
Two neural networks are considered, with the architectures given by
$\mbs M_1=[2, 800, 50, 1]$ and $\mbs M_2=[2, 50, 800, 1]$, respectively.
We vary the number of collocation points per direction (i.e.~$Q_1$)
systematically between $Q_1=5$ and $Q_1=50$, and record the corresponding
HLConcELM errors.
Figures~\ref{fg_4}(a) and (b) show the maximum/rms errors of HLConcELM
as a function of $Q_1$ for the two neural networks.
For the network $\mbs M_1$ we employ a hidden magnitude vector
$\mbs R=(3.0, 0.005)$, and for the network $\mbs M_2$ we employ
a hidden magnitude vector $\mbs R=(0.68,0.82)$. These $\mbs R$ values
are obtained using the method of~\cite{DongY2021}.
The results indicate that the HLConcELM errors decrease approximately exponentially
with increasing number of collocation points (when $Q_1\leqslant 30$). 
The errors stagnate as $Q_1$ increases further, because of the fixed network
size. Note that  the last hidden layer of
the network $\mbs M_1$ is quite narrow ($50$ nodes), while that of
the network $\mbs M_2$ is quite wide ($800$ nodes).
The HLConcELM method produces accurate results with both types of
neural networks.

\begin{figure}
  \centerline{
    \includegraphics[width=2.in]{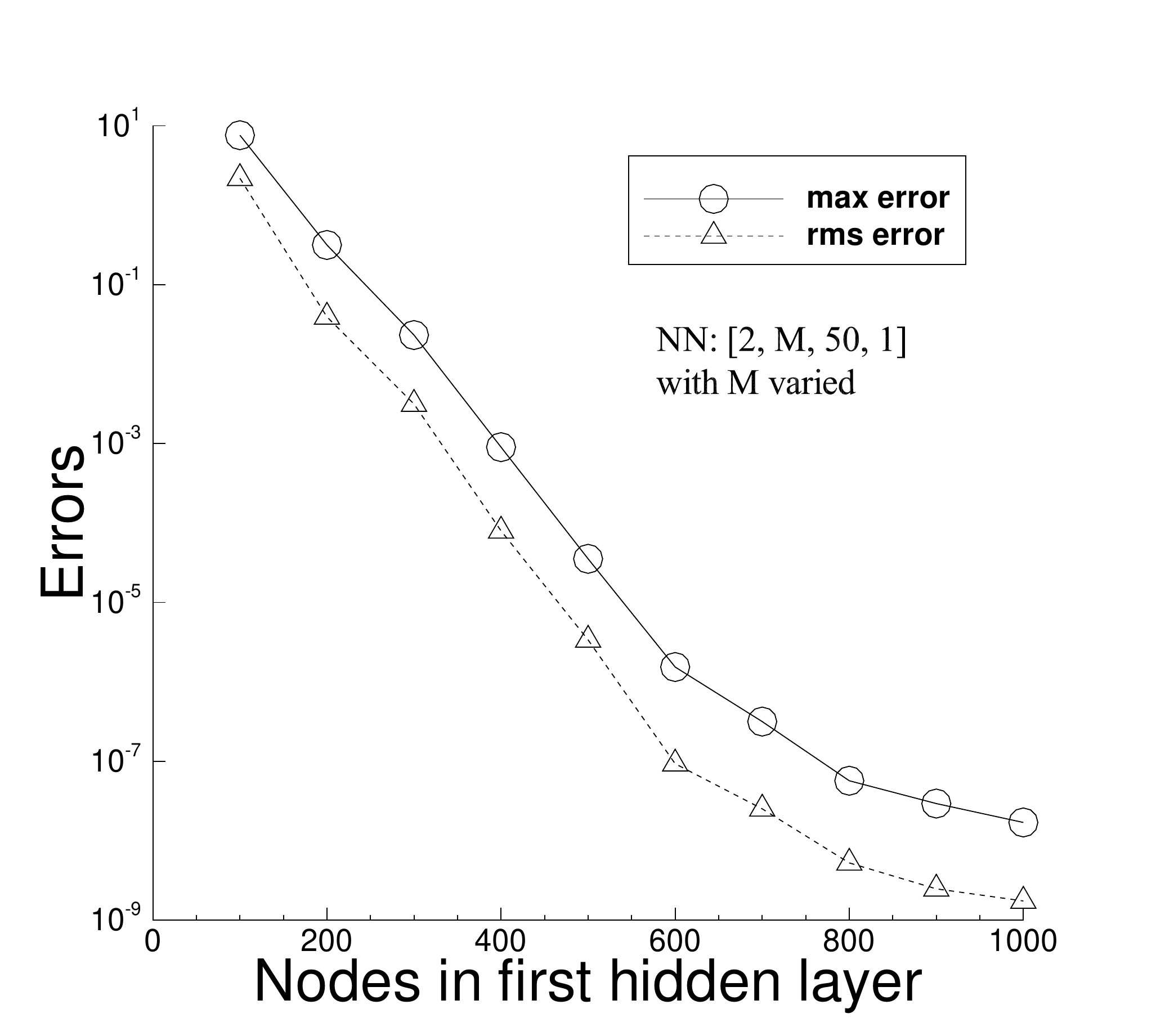}(a)
    \includegraphics[width=2.in]{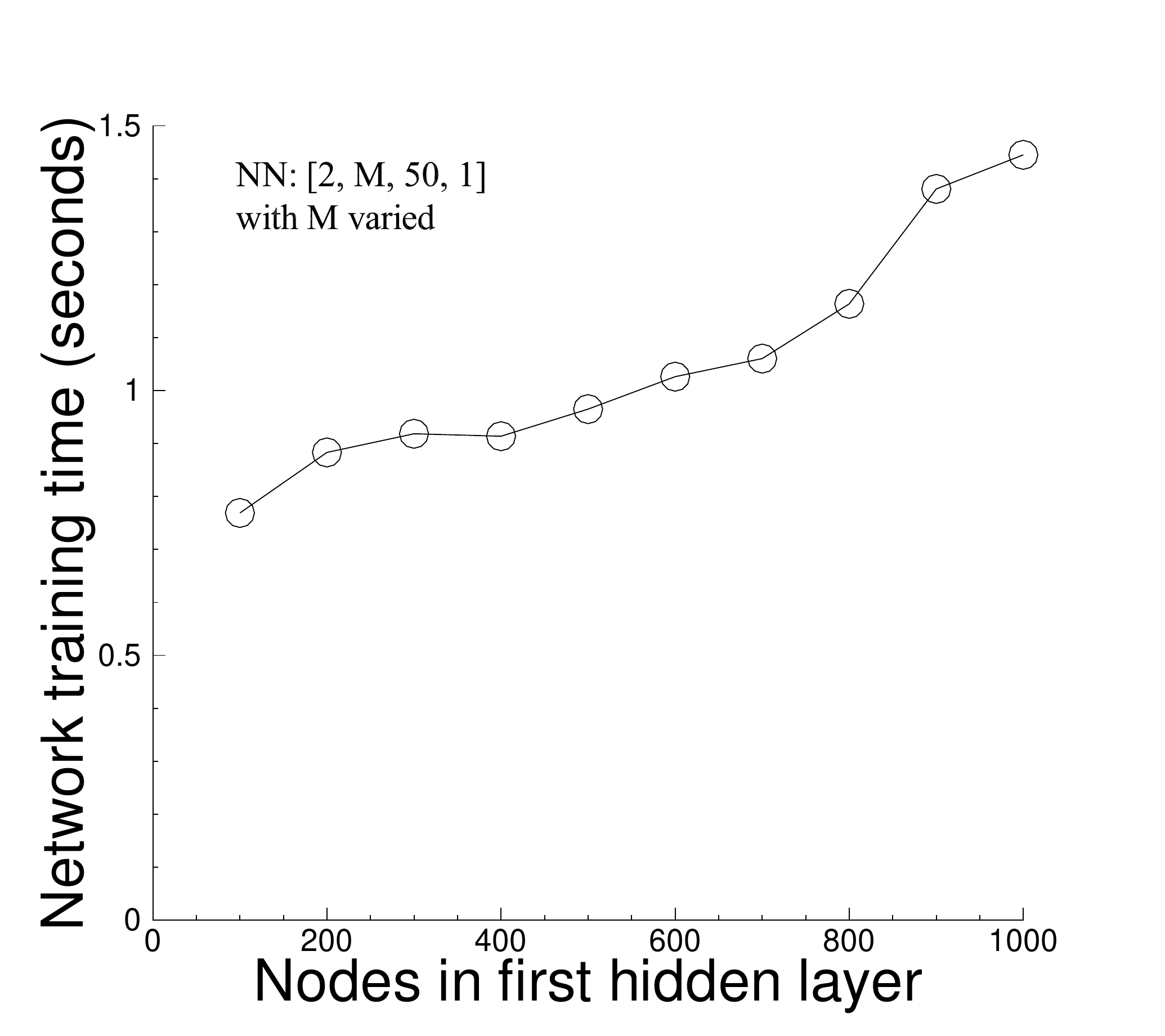}(b)
  }
  \centerline{
    \includegraphics[width=2.in]{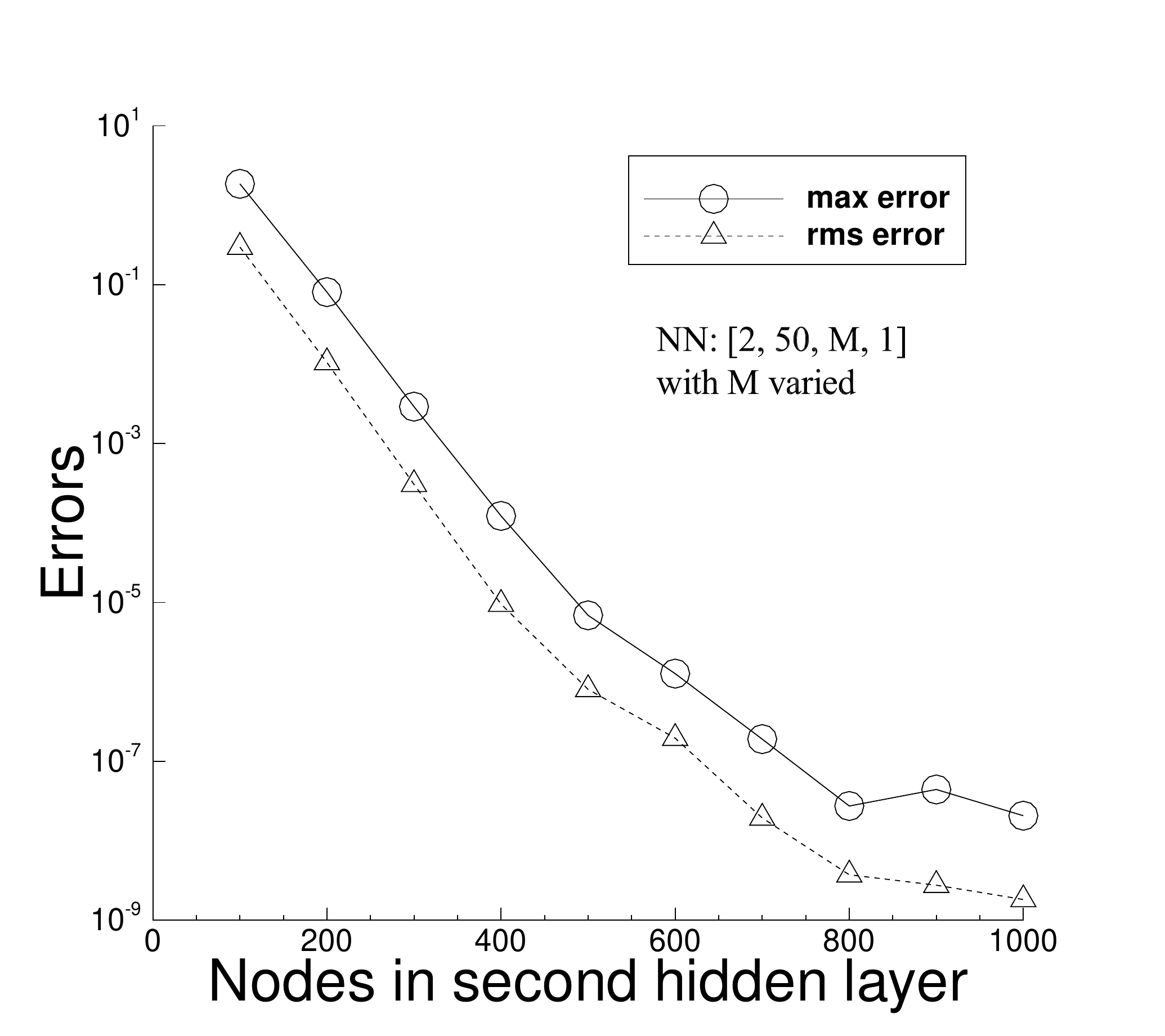}(c)
    \includegraphics[width=2.in]{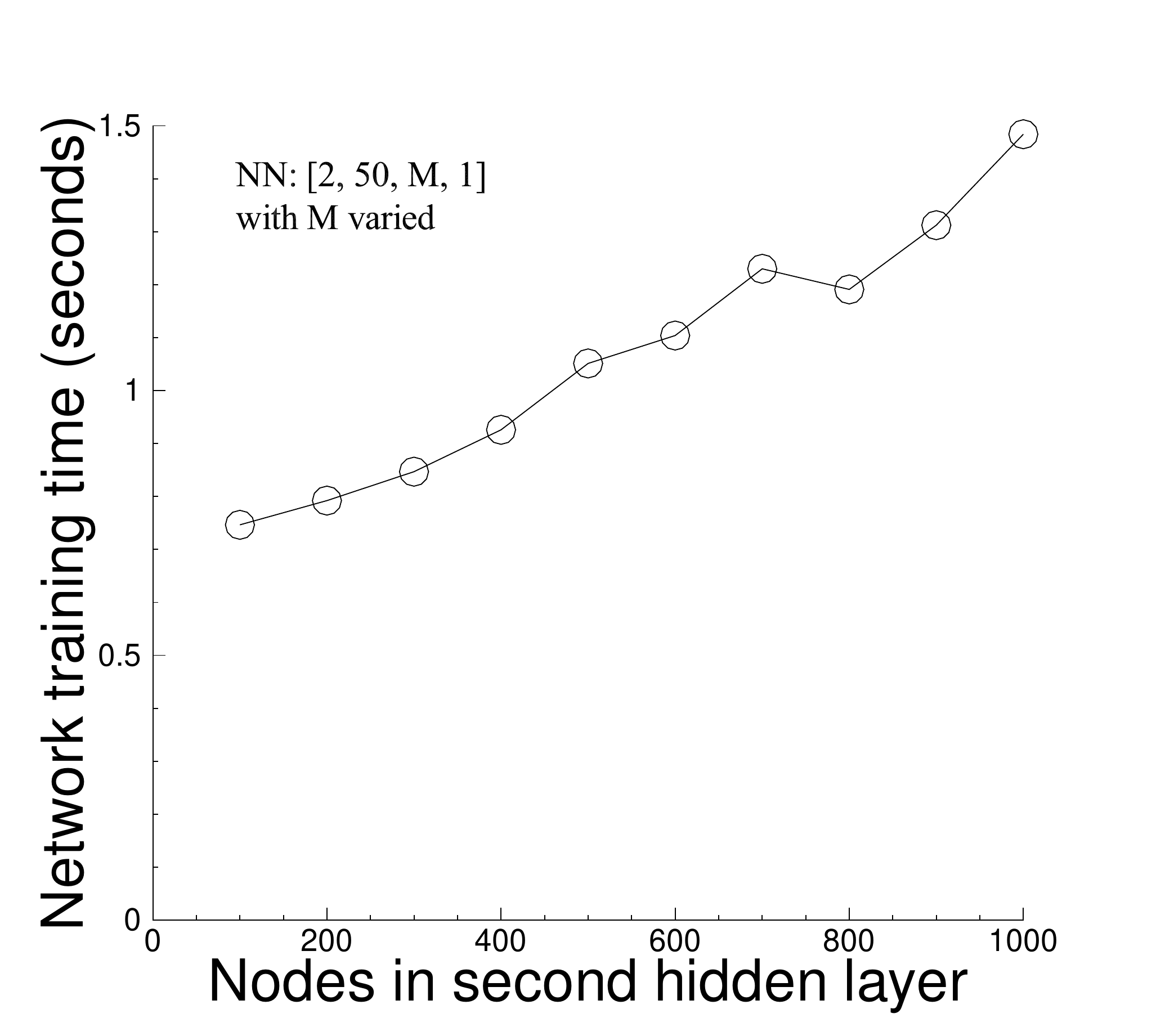}(d)
  }
  \caption{Variable-coefficient
    Poisson equation: (a) The maximum/rms errors and (b) the network
    training time of HLConcELM versus
    the number of nodes in the first hidden layer for a network architecture
    $[2, M, 50, 1]$ (varying $M$).
    (c) The maximum/rms errors and (d) the network training time
    of HLConcELM versus the number of nodes in
    the second hidden layer for a network architecture $[2, 50, M, 1]$ (varying $M$).
    $Q=35\times 35$ uniform collocation points in (a,b,c,d).
    $\mbs R=(3.0, 0.005)$ in (a,b), and $\mbs R=(0.68,0.82)$ in (c,d).
  }
  \label{fg_5}
\end{figure}

Figure~\ref{fg_5} illustrates the convergence behavior, as well as the
network training time, of the HLConcELM method with respect to
the number of nodes in the neural network.
We consider two groups of neural networks, with the architectures
given by $\mbs M_1=[2, M, 50, 1]$ and $\mbs M_2=[2, 50, M, 1]$, respectively,
where $M$ is varied systematically between $M=100$ and $M=1000$.
For all the test cases, we employ a fixed uniform set of $Q=35\times 35$
collocation points to train the neural network.
For generating the hidden-layer coefficients, we use 
a hidden magnitude vector $\mbs R=(3.0,0.005)$ with
the first group of networks $\mbs M_1$, and a vector $\mbs R=(0.68,0.82)$
with the second group of networks $\mbs M_2$.
Figures~\ref{fg_5}(a) and (c) depict the maximum/rms errors of HLConcELM
as a function of $M$ for these two groups of neural networks.
Figures~\ref{fg_5}(b) and (d) depict the corresponding wall time it takes to
train these neural networks with HLConcELM.
It can be observed that the HLConcELM errors decrease approximately
exponentially with increasing $M$ (before saturation).
When  $M$ becomes large the HLConcELM results are
highly accurate. The network training time of the HLConcELM method
increases approximately linearly with increasing $M$. In the range of $M$
values tested here, it takes around a second to train
the neural network to attain the HLConcELM results.


\begin{table}[tb]
  \centering
  \begin{tabular}{l|l|cc|cc}
    \hline
    network & collocation & current & HLConcELM & conventional & ELM  \\ \cline{3-6}
    architecture & points & max error & rms error & max error & rms error \\ \hline
    $[2,800,50,1]$ & $5\times 5$ & $1.91E+0$ & $4.31E-1$ & $3.06E+1$ & $6.34E+0$  \\
    & $10\times 10$ & $3.22E-2$ & $7.88E-3$ & $5.48E+1$ & $1.78E+1$  \\
    & $15\times 15$ & $2.33E-3$ & $3.92E-4$ & $6.24E+1$ & $2.11E+1$  \\
    & $20\times 20$ & $4.70E-5$ & $1.32E-5$ & $6.97E+1$ & $2.42E+1$  \\
    & $25\times 25$ & $4.78E-7$ & $1.10E-7$ & $7.67E+1$ & $2.70E+1$  \\
    & $30\times 30$ & $3.17E-8$ & $3.79E-9$ & $8.31E+1$ & $2.96E+1$  \\ \hline
    $[2,50,800,1]$ & $5\times 5$ & $2.95E+0$ & $9.26E-1$ & $2.48E+0$ & $8.85E-1$  \\
    & $10\times 10$ & $9.35E-2$ & $9.93E-3$ & $1.32E-1$ & $1.50E-2$ \\
    & $15\times 15$ & $1.11E-3$ & $2.40E-4$ & $6.52E-3$ & $1.00E-3$  \\
    & $20\times 20$ & $3.42E-5$ & $6.91E-6$ & $7.63E-5$ & $1.33E-5$ \\
    & $25\times 25$ & $2.34E-6$ & $4.45E-7$ & $1.83E-6$ & $4.14E-7$ \\
    & $30\times 30$ & $3.07E-8$ & $4.81E-9$ & $9.87E-8$ & $2.02E-8$ \\
    \hline
  \end{tabular}
  \caption{Variable-coefficient
    Poisson equation: Comparison of the maximum/rms errors
    computed using the current HLConcELM method and the
    conventional ELM method.
    The HLConcELM data  in this table correspond to a portion of those in
    Figure~\ref{fg_4}(a) for the network $[2, 800, 50, 1]$ and to
    those in Figure~\ref{fg_4}(b) for the network $[2, 50, 800, 1]$.
    For conventional ELM, the random hidden-layer coefficients
    are assigned to uniform random values generated on $[-R_m,R_m]$
    with $R_m=R_{m0}$. Here $R_{m0}$ is the optimal $R_m$ obtained
    using the method of~\cite{DongY2021}, with $R_{m0}=0.35$ for the network
    $[2,800,50,1]$ and $R_{m0}=0.75$ for the network $[2, 50, 800, 1]$.
  }
  \label{tab_1}
\end{table}

Table~\ref{tab_1} provides a comparison of the HLConcELM accuracy
and the conventional ELM accuracy
for solving the variable-coefficient Poisson equation
on two network architectures,
$\mbs M_1=[2, 800, 50, 1]$ and $\mbs M_2=[2, 50, 800, 1]$.
The network $\mbs M_1$ contains a relatively
small number of nodes in its last hidden layer,
and the conventional ELM would not perform well.
The network $\mbs M_2$ contains a large number of nodes in
its last hidden layer,
and the conventional ELM should perform quite well.
We consider a sequence of uniform collocation points,
ranging from $Q=5\times 5$ to $Q=30\times 30$.
Table~\ref{tab_1} lists the maximum/rms errors of the HLConcELM
solution and the conventional ELM solution corresponding to each set of collocation
points. 
The data indicate that the conventional ELM exhibits no accuracy
with the network $\mbs M_1$, and exhibits exponentially increasing
accuracy with increasing collocation points on the network $\mbs M_2$.
On the other hand, the current HLConcELM method exhibits exponentially
increasing accuracy with increasing collocation points on
both networks $\mbs M_1$ and $\mbs M_2$.

\begin{figure}
  \centerline{
    \includegraphics[width=2.in]{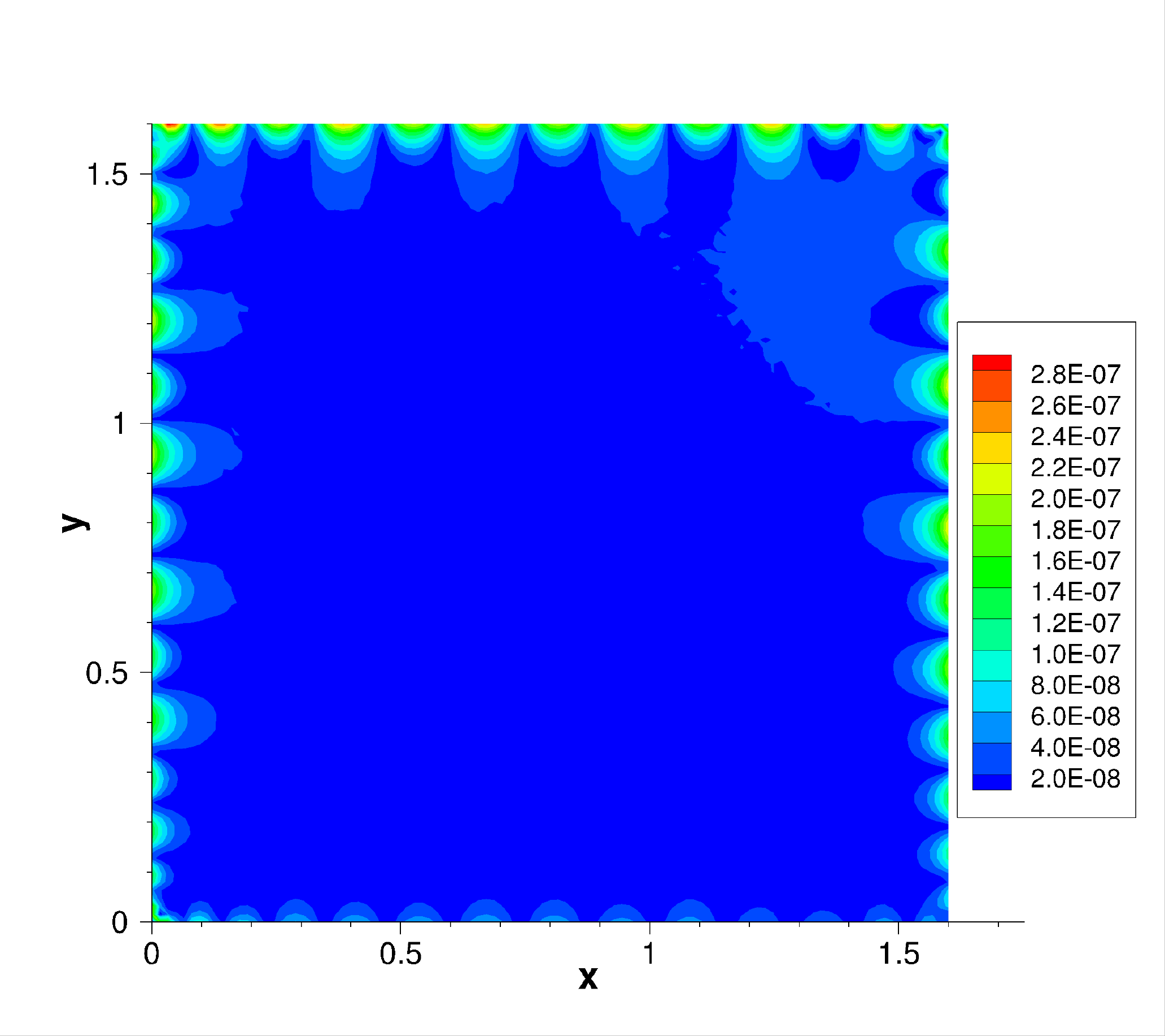}(a)
    \includegraphics[width=2.in]{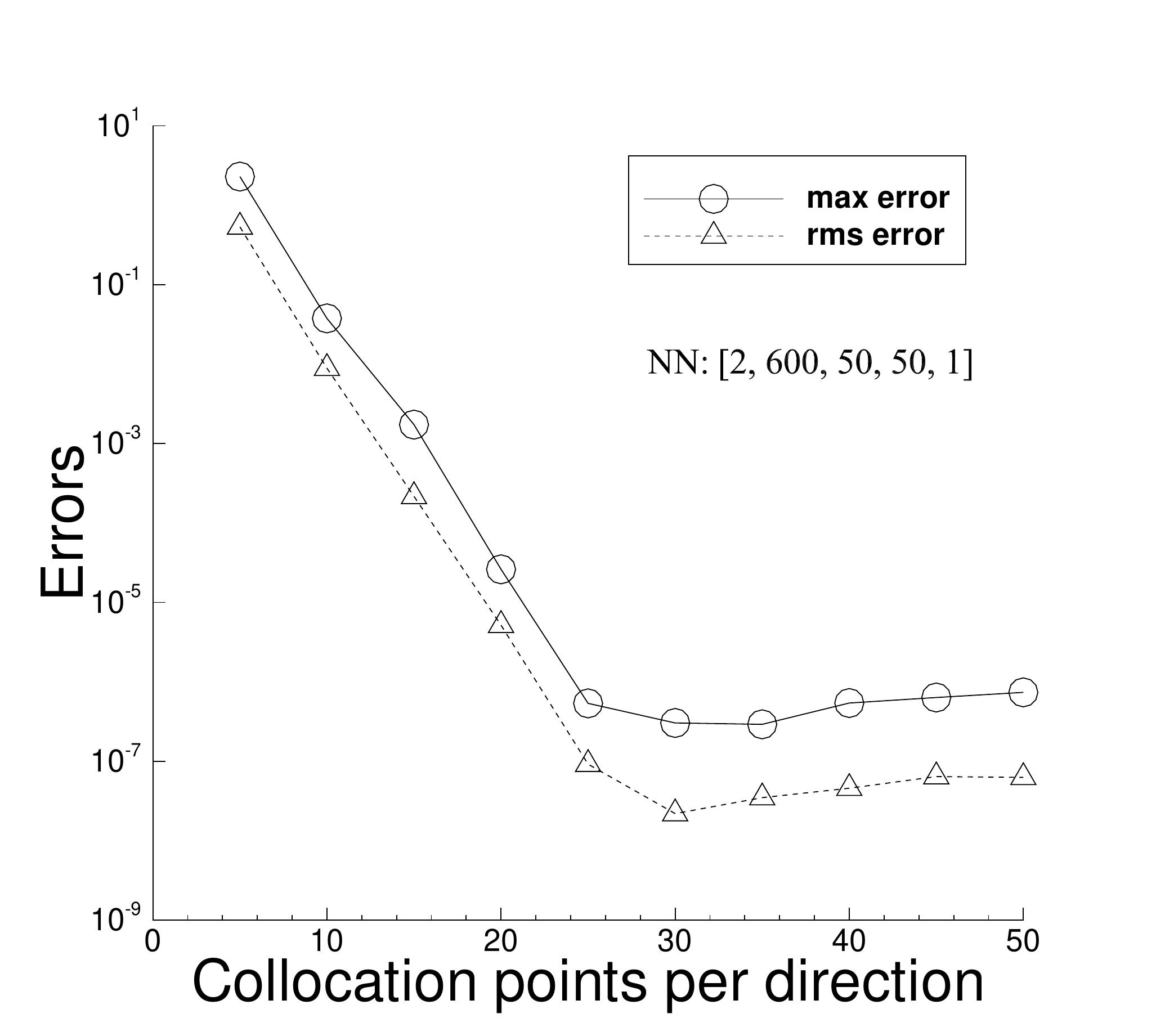}(b)
    \includegraphics[width=2.in]{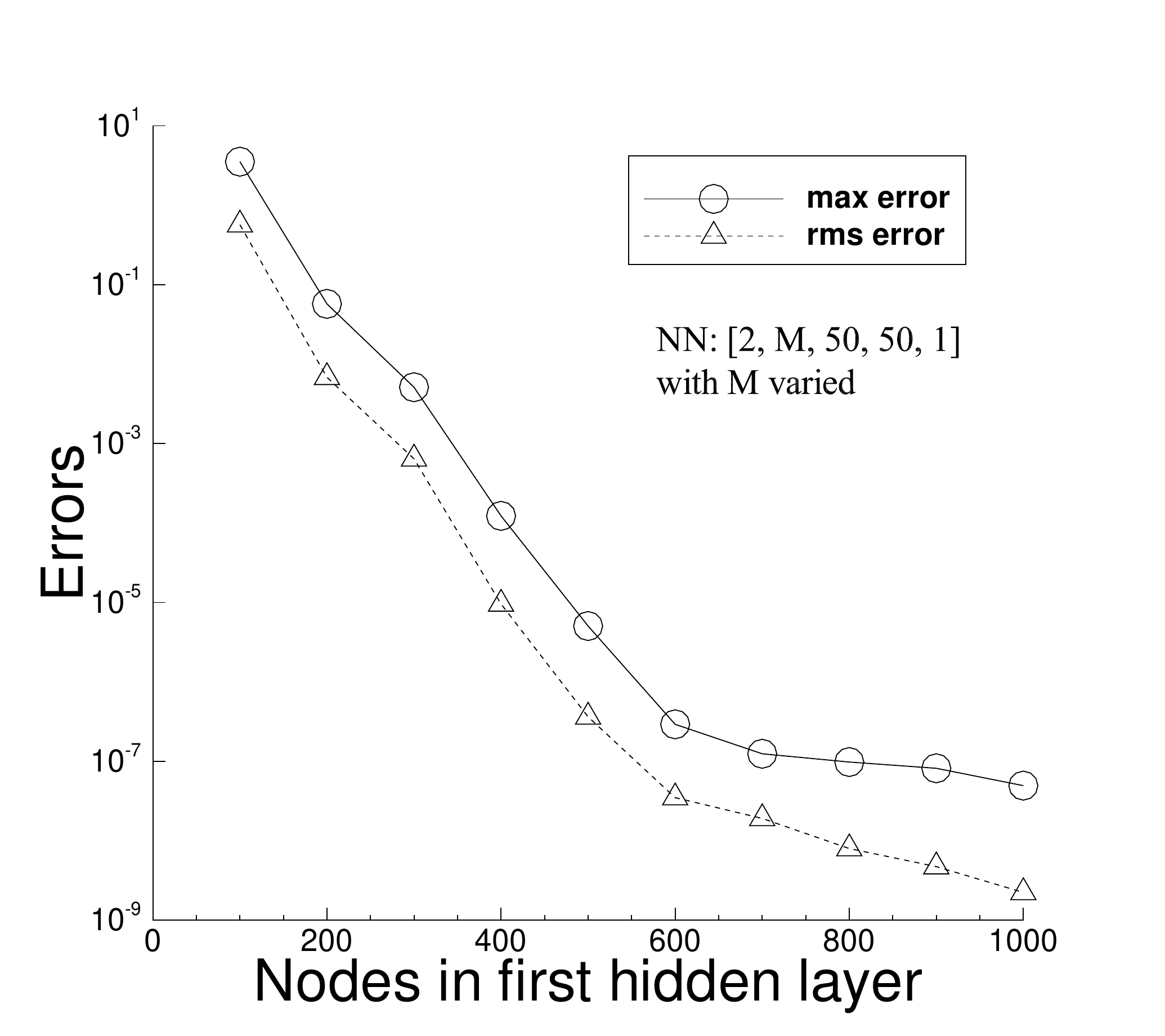}(c)
  }
  \caption{Variable-coefficient
    Poisson equation (3 hidden layers in NN): (a) Distribution of
    the absolute error of the HLConcELM solution.
    The maximum/rms errors of HLConcELM versus (b)
    the number of collocation points per direction, and (c) the number
    of the nodes in the first hidden layer of the neural network.
    Neural network architecture $\mbs M=[2, M, 50, 50, 1]$,
    $Q=Q_1\times Q_1$ uniform collocation points.
    $M=600$ in (a,b) and is varied in (c).
    $Q_1=35$ in (a,c) and is varied in (b).
    $\mbs R=(2.6, 0.005, 0.8)$ in (a,b,c).
  }
  \label{fg_6}
\end{figure}

Figure~\ref{fg_6} is an illustration of the HLConcELM results obtained
on neural networks with three hidden layers.
Here we consider a network architecture $\mbs M=[2, M, 50, 50, 1]$,
with $M$ either fixed at $M=600$ or varied systematically between
$M=100$ and $M=1000$.
The set of collocation points (uniform) is either fixed at $Q=35\times 35$
or varied systematically between $Q=5\times 5$ and $Q=50\times 50$.
We employ a fixed hidden magnitude vector $\mbs R=(2.6,0.005,0.8)$,
obtained using the method of~\cite{DongY2021}.
Figure~\ref{fg_6}(a) shows the HLConcELM error distribution
corresponding to $\mbs M=[2, 600, 50, 50, 1]$ and
$Q=35\times 35$, indicating a quite high accuracy, with the maximum
error in the domain on the order $10^{-7}$.
Figures~\ref{fg_6}(b) and (c) demonstrate the exponential convergence
(before saturation)
of the HLConcELM errors with respect to the collocation points $Q_1$
and the number of nodes $M$, respectively. These results show
that the current HLConcELM method can produce highly accurate results
on neural networks
with multiple hidden layers and a narrow last hidden layer.

\subsection{Advection Equation}
\label{sec:advec}

\begin{figure}
  \centering
  \includegraphics[width=5in]{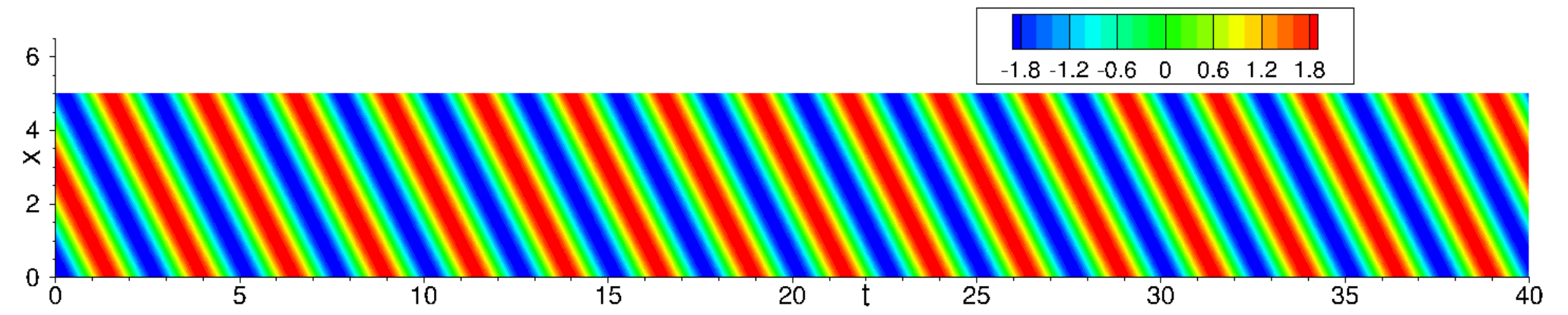}(a)
  \includegraphics[width=5in]{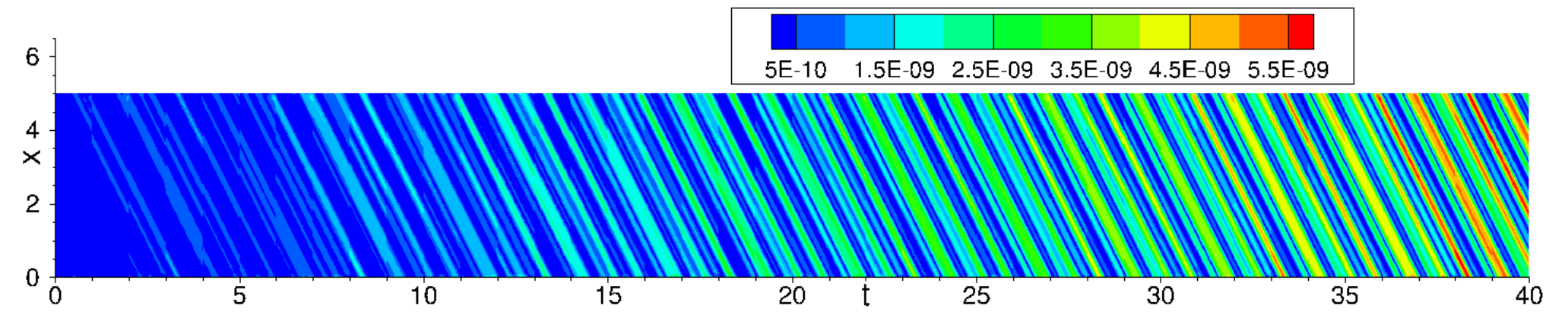}(b)
  \includegraphics[width=5in]{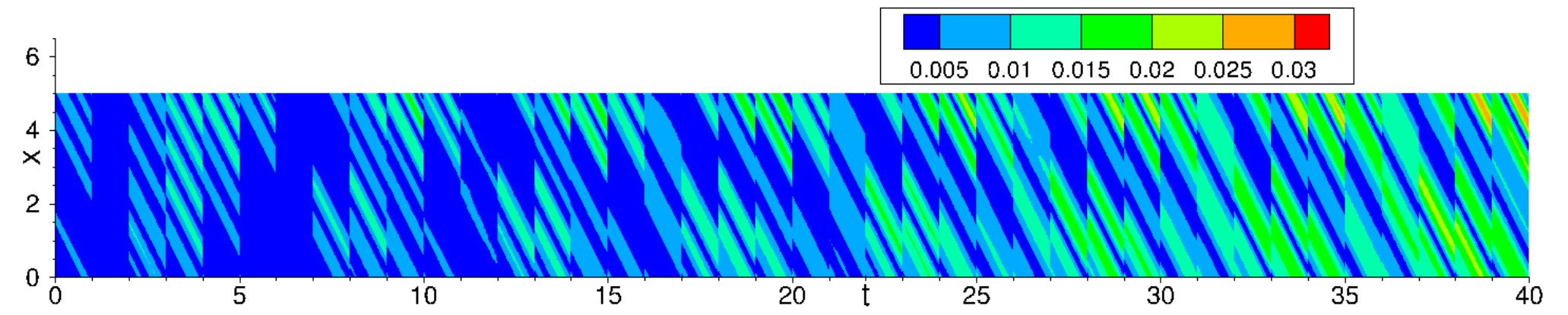}(c)
  \caption{Advection equation: Distributions in the spatial-temporal domain
    of (a) the exact solution,
    (b) the absolute error of the HLConcELM solution, and (c) the absolute error
    of the conventional ELM solution.
    In (b,c), network architecture [2, 500, 50, 1],
    $40$ uniform time blocks,
    $Q=35\times 35$ uniform collocation points per time block.
    $\mbs R=(3.0,1.0)$ in (b) for HLConcELM.
    $R_m=R_{m0}=0.065$ in (c) for conventional ELM.
  }
  \label{fg_7}
\end{figure}

In the next example we employ the 1D advection equation (plus time)
to test the  HLConcELM method.
Consider the spatial-temporal domain,
$(x,t)\in\Omega=[0,5]\times [0,40]$, and the following initial/boundary
value problem on $\Omega$,
\begin{subequations}\label{eq_23}
  \begin{align}
    &
    \frac{\partial u}{\partial t} - 2\frac{\partial u}{\partial x} = 0, \\
    &
    u(0,t) = u(5,t), \\
    &
    u(x,0) = 20\tanh\left(\frac{1}{10}\cos\left(\frac{2\pi}{5}\left(x -3 \right) \right) \right).
  \end{align}
\end{subequations}
In the above equations $u(x,t)$ is the field function to be solved for,
and we impose the periodic boundary condition in the spatial direction.
This system has the following exact solution,
\begin{equation}\label{eq_24}
  u(x,t) = 20\tanh\left(\frac{1}{10}\cos\left(\frac{2\pi}{5}\left(x + 2t -3 \right) \right) \right).
\end{equation}
The distribution of this solution on the spatial-temporal domain
is illustrated in Figure~\ref{fg_7}(a).

To solve the system~\eqref{eq_23}, we employ the HLConcELM method combined
with the block time marching scheme (see Remark~\ref{rem_3} and~\cite{DongL2021}).
We divide the domain $\Omega$ into $40$ uniform time blocks in time.
For computing each time block with HLConcELM,
we employ a network architecture
$\mbs M=[2, m_1,\dots,m_{L-1},1]$, where the two input nodes represent $x$ and $t$
and the single output node represents $u(x,t)$.
Let $Q=Q_1\times Q_1$ denote the uniform set of collocation points for each
time block ($Q_1$ grid points in both $x$ and $t$ directions), where $Q_1$ is varied
in the tests.
As discussed before,
upon completion of training, the neural network is evaluated
on a uniform set of $Q_{eval}=101\times 101$ grid points on each time block
and the corresponding errors are computed.
The maximum and rms errors reported below refer to
the errors of the HLConcELM solution on
the entire domain $\Omega$ (over $40$ time blocks).


Figures~\ref{fg_7}(b) and (c) illustrate the absolute-error distributions
on $\Omega$ of the HLConcELM solution and the conventional ELM solution, respectively.
For both methods, we employ $40$ time blocks in block time marching,
a neural network architecture $\mbs M=[2, 500, 50, 1]$ with the Gaussian activation
function, and a set of $Q=35\times 35$ uniform collocation points per time block.
For HLConcELM we employ $\mbs R=(3.0,1.0)$, which is computed by the method of~\cite{DongY2021}.
For conventional ELM we employ $R_m=R_{m0}=0.065$, which is also obtained
by the method of~\cite{DongY2021}, for generating
the random hidden-layer coefficients.
Because the number of nodes in the last hidden layer is quite small,
the conventional ELM exhibits a low accuracy, with the maximum error on
the order of $10^{-2}$ in the domain.
On the other hand, the HLConcELM method produces a highly accurate
solution, with the maximum error on the order of $10^{-9}$ in the domain.

\begin{figure}
  \centerline{
    \includegraphics[width=2.in]{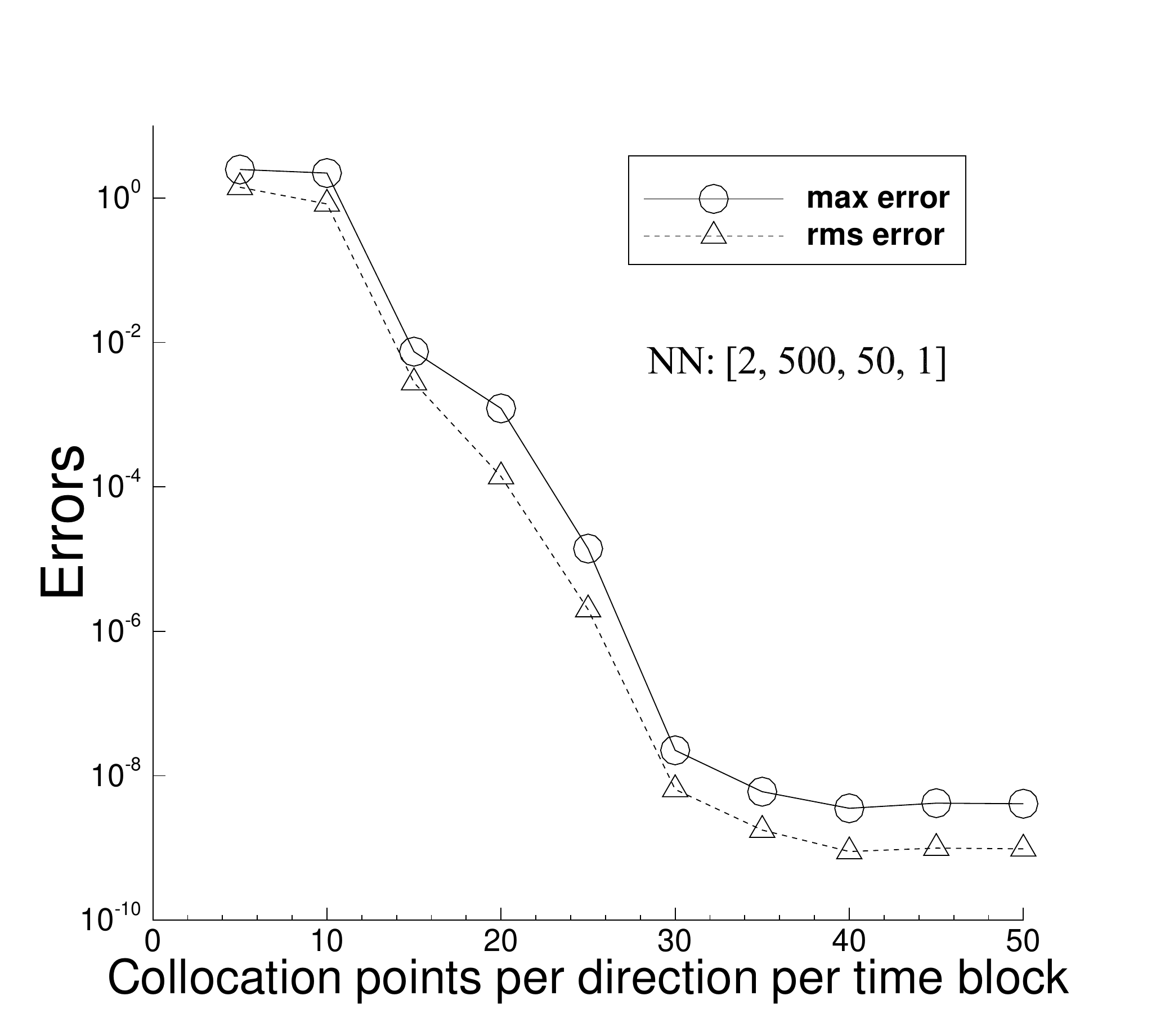}(a)
    \includegraphics[width=2.in]{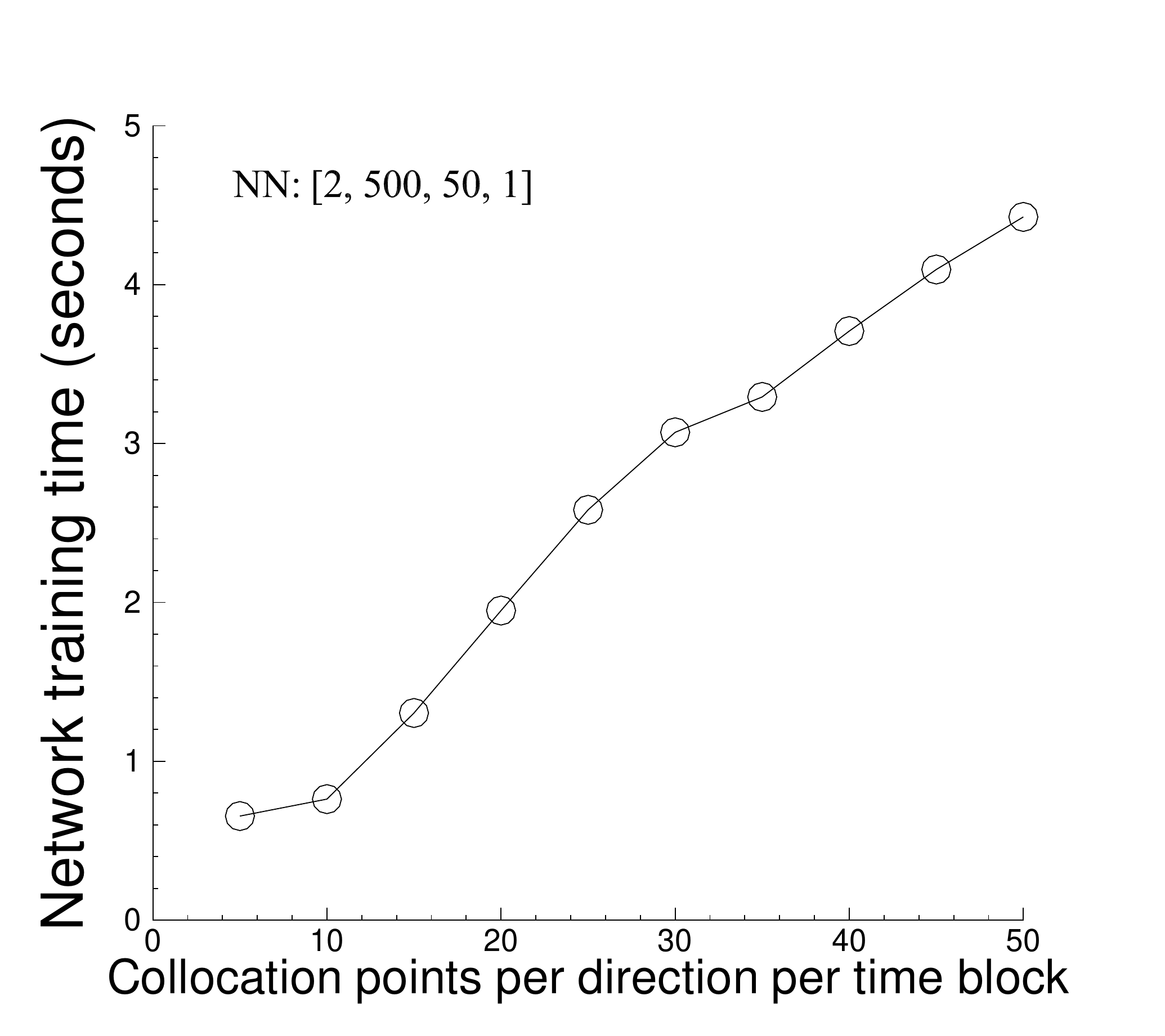}(b)
  }
  \centerline{
    \includegraphics[width=2.in]{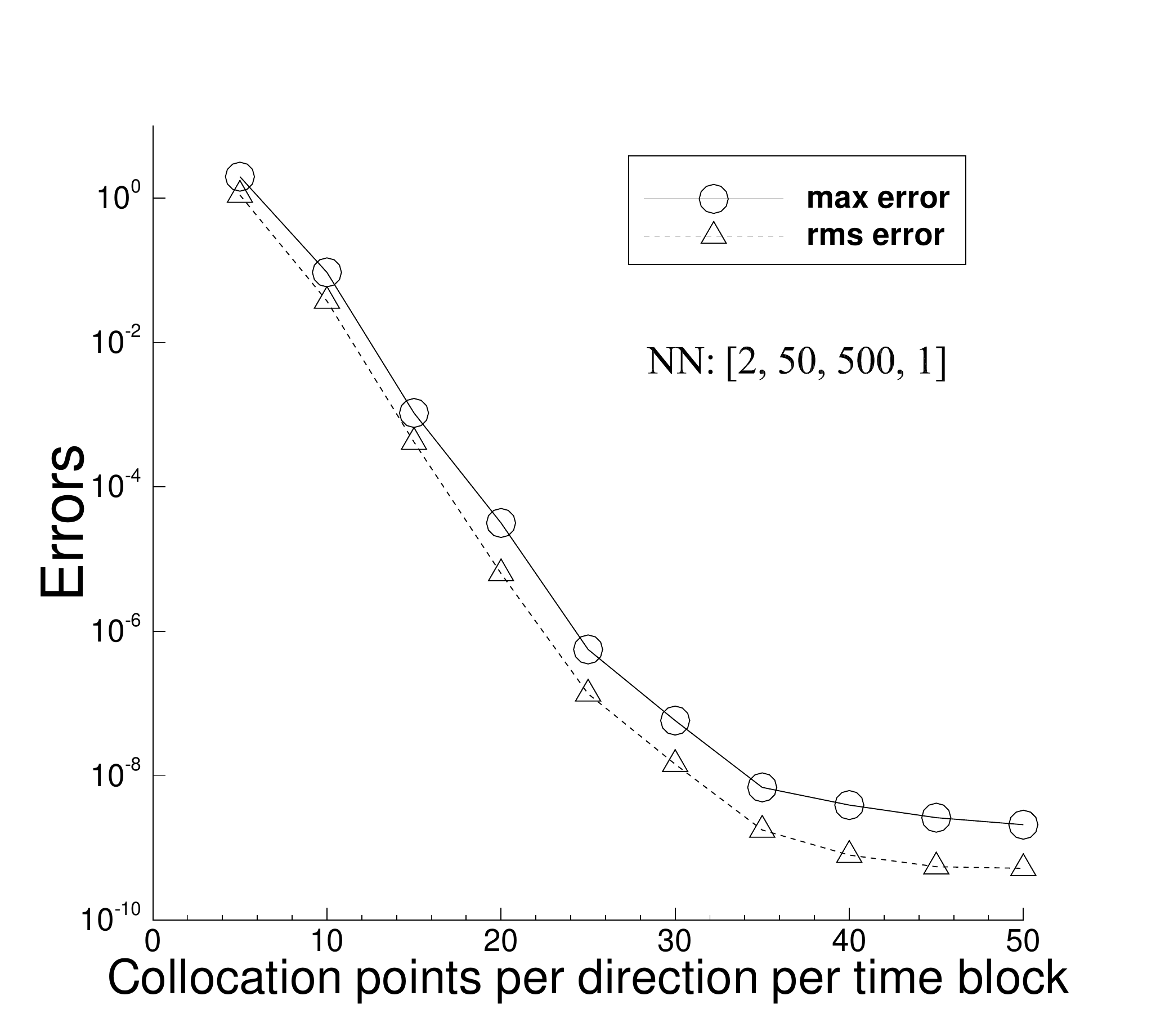}(c)
    \includegraphics[width=2.in]{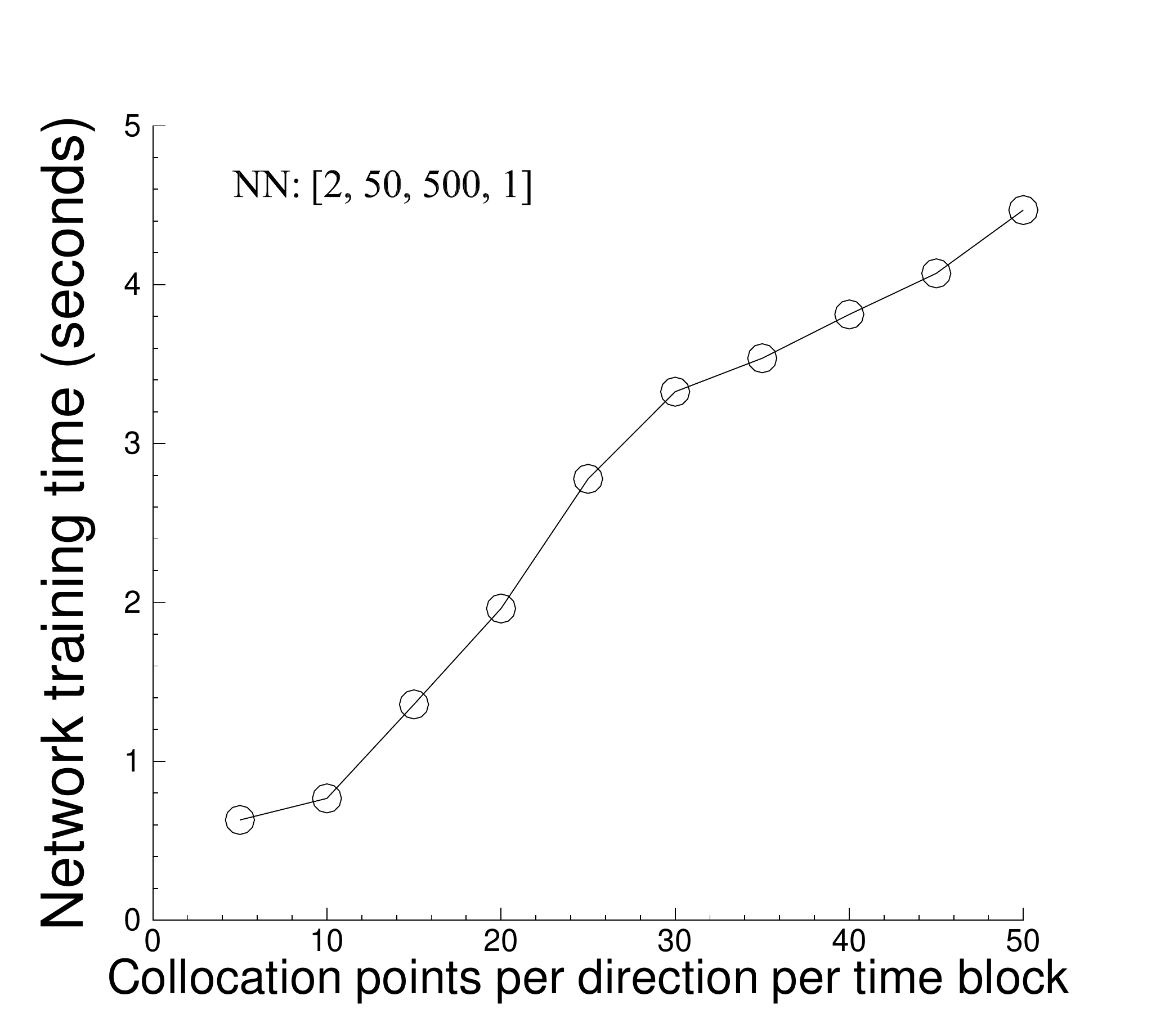}(d)
  }
  \caption{Advection equation:
    The maximum/rms errors in the domain $\Omega$ (a,c)
    and the network training time (b,d) of HLConcELM versus
    the number of collocation points per direction in each time block.
    The results are attained with two network architectures:
    (a,b) $\mbs M_1=[2,500,50,1]$,
    and (c,d) $\mbs M_2=[2,50,500,1]$.
    $\mbs R=(3.0,1.0)$ in (a,b) for the network $\mbs M_1$,
    and $\mbs R=(0.9,0.5)$ in (c,d) for the network $\mbs M_2$.
  }
  \label{fg_8}
\end{figure}

Figure~\ref{fg_8} illustrates the convergence behavior, as well as the growth
in the network training time, of the HLConcELm method
with respect to the number of collocation points.
We have considered two network architectures,
$\mbs M_1=[2, 500, 50, 1]$ and $\mbs M_2=[2, 50, 500, 1]$,
with a narrower last hidden layer in $\mbs M_1$ and
a wider one in $\mbs M_2$.
A uniform set of $Q=Q_1\times Q_1$ collocation points is employed,
with $Q_1$ varied systematically between $Q_1=5$ and $Q_1=50$ in the tests.
The hidden magnitude vector $\mbs R$ computed by the method~\cite{DongY2021}
is used in the simulations, with $\mbs R=(3.0,1.0)$ for the network $\mbs M_1$
and $\mbs R=(0.9,0.5)$ for the network $\mbs M_2$.
Figures~\ref{fg_8}(a) and (b) depict the maximum/rms errors on $\Omega$
and the network training time, respectively, as a function of $Q_1$
obtained with the neural network $\mbs M_1$.
Figures~\ref{fg_8}(c) and (d) show the corresponding results obtained
with the network $\mbs M_2$.
While the convergence behavior is not quite regular, one can observe that
the HLConcELM errors approximately decrease exponentially (before saturation)
with increasing number of collocation points.
The network training time grows approximately linearly with increasing
number of training collocation points.

\begin{figure}
  \centerline{
    \includegraphics[width=2.in]{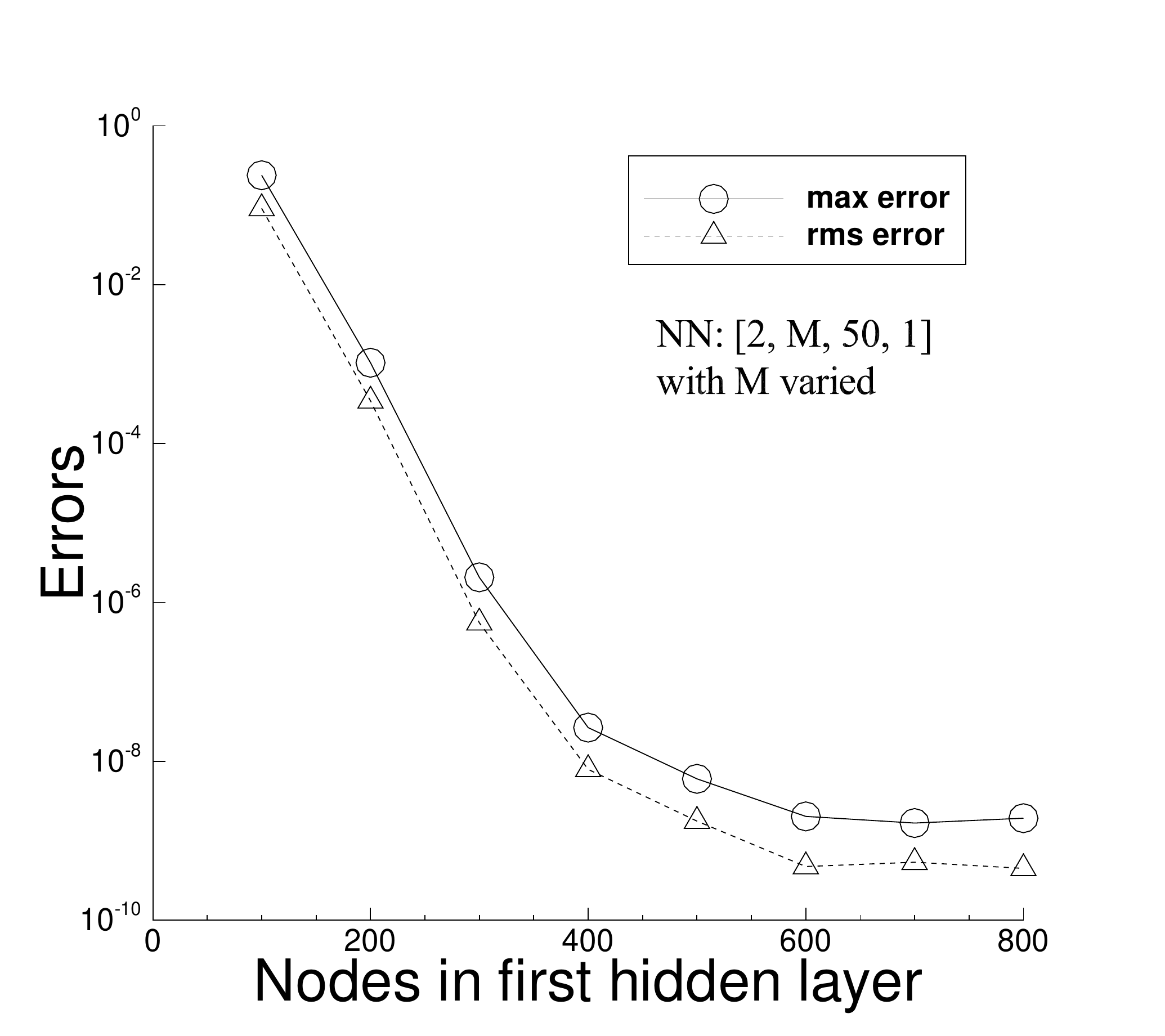}(a)
    \includegraphics[width=2.in]{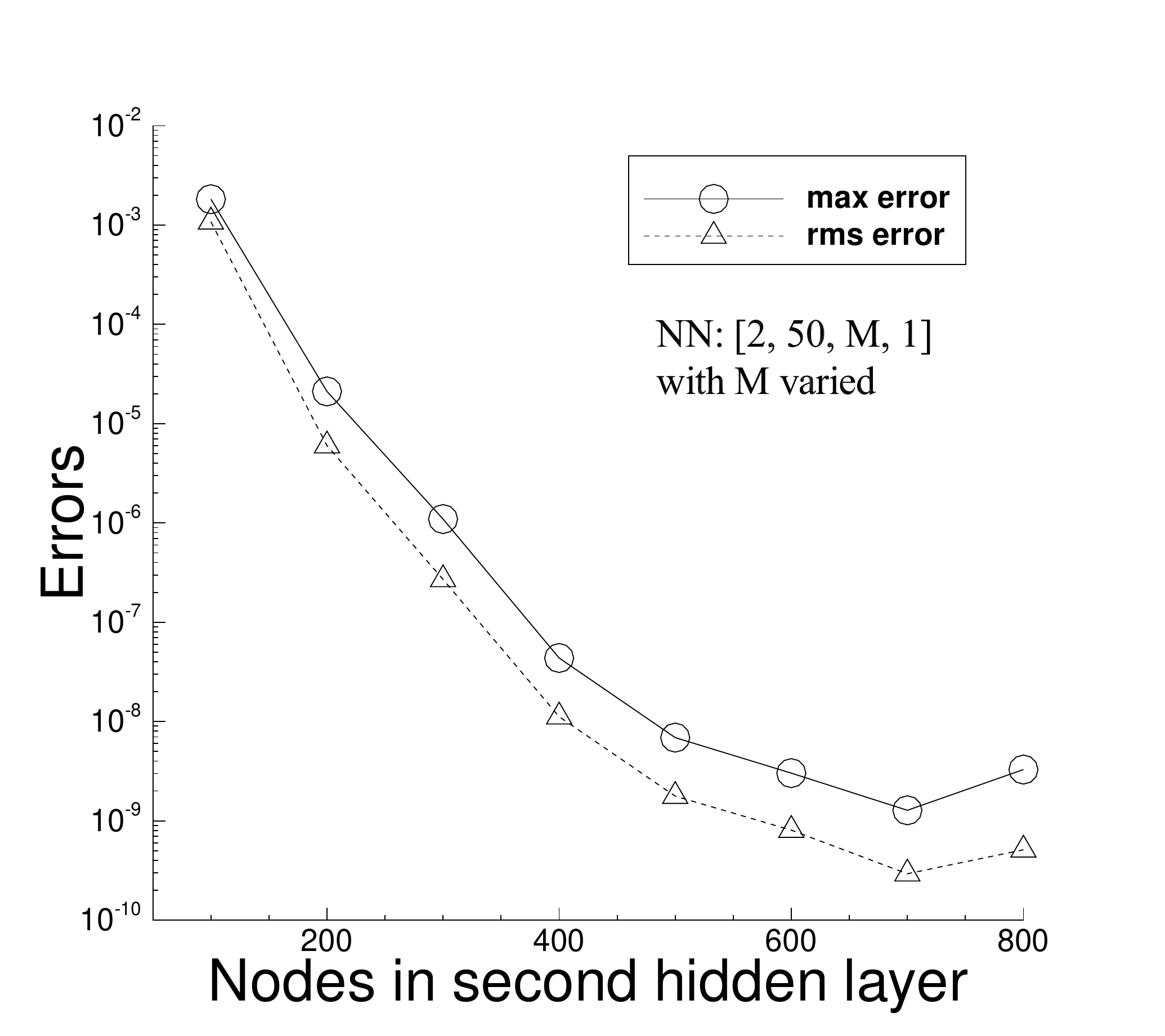}(b)
  }
  \caption{Advection equation:
    (a) The HLConcELM maximum/rms errors on $\Omega$ versus
    the number of nodes in the first hidden layer for the network architecture
    $\mbs M_1=[2, M, 50, 1]$ (varying $M$).
    (b) The HLConcELM maximum/rms errors on $\Omega$  versus the number of
    nodes in the second hidden layer for the network architecture $\mbs M_2=[2, 50, M, 1]$
    (varying $M$). $Q=35\times 35$ in (a,b).
    $\mbs R=(3.0,1.0)$ in (a) for the network $\mbs M_1$,
    and $\mbs R=(0.9,0.5)$ in (b) for the network $\mbs M_2$.
  }
  \label{fg_9}
\end{figure}

Figure~\ref{fg_9} illustrates the convergence behavior of the HLConcELM method
with respect to the number of nodes in the neural network.
Two groups of neural networks are considered in these tests,
with an architecture $\mbs M_1=[2, M, 50, 1]$ for the first group
and $\mbs M_2=[2, 50, M, 1]$ for the second one,
with $M$ varied systematically.
A uniform set of $Q=35\times 35$ collocation points is employed for
training the neural networks.
We use $\mbs R=(3.0,1.0)$ for the architecture $\mbs M_1$ and
$\mbs R=(0.9,0.5)$ for the architecture $\mbs M_2$.
The plots (a) and (b) show the maximum/rms errors of HLConcELM on $\Omega$
as a function of $M$, indicating that the errors decrease approximately 
exponentially (before saturation)
with increasing $M$ in the neural network.


\begin{table}[tb]
  \centering
  \begin{tabular}{l|l|cc|cc}
    \hline
    network & collocation & current & HLConcELM & conventional & ELM  \\ \cline{3-6}
    architecture & points & max error & rms error & max error & rms error \\ \hline
    $[2,500,50,1]$ & $5\times 5$ & $2.48E+0$ & $1.41E+0$ & $1.33E+0$ & $4.57E-1$  \\
    & $10\times 10$ & $2.21E+0$ & $8.26E-1$ & $5.97E-2$ & $1.91E-2$  \\
    & $15\times 15$ & $7.45E-3$ & $2.78E-3$ & $4.31E-2$ & $1.30E-2$  \\
    & $20\times 20$ & $1.22E-3$ & $1.39E-4$ & $3.62E-2$ & $1.08E-2$ \\
    & $25\times 25$ & $1.39E-5$ & $2.00E-6$ & $3.15E-2$ & $9.50E-3$  \\
    & $30\times 30$ & $2.25E-8$ & $6.47E-9$ & $3.24E-2$ & $8.89E-3$  \\ \hline
    $[2,50,500,1]$ & $5\times 5$ & $1.97E+0$ & $1.10E+0$ & $1.86E+0$ & $9.18E-1$  \\
    & $10\times 10$ & $9.33E-2$ & $3.74E-2$ & $4.18E-2$ & $1.65E-2$  \\
    & $15\times 15$ & $1.05E-3$ & $4.16E-4$ & $3.09E-4$ & $8.48E-5$  \\
    & $20\times 20$ & $3.16E-5$ & $6.30E-6$ & $1.97E-4$ & $5.01E-5$ \\
    & $25\times 25$ & $5.60E-7$ & $1.36E-7$ & $7.60E-5$ & $5.85E-6$ \\
    & $30\times 30$ & $5.80E-8$ & $1.45E-8$ & $9.76E-8$ & $3.32E-8$  \\
    \hline
  \end{tabular}
  \caption{Advection equation: Comparison of the maximum/rms errors on $\Omega$
    from the current HLConcELM method and the
    conventional ELM method~\cite{DongL2021}.
    The HLConcELM data in this table correspond to a portion of those in
    Figure~\ref{fg_8}(a) for the network $[2, 500, 50, 1]$ and to
    those in Figure~\ref{fg_8}(c) for the network $[2, 50, 500, 1]$.
    For conventional ELM, the hidden-layer coefficients
    are set to uniform random values generated on $[-R_m,R_m]$
    with $R_m=R_{m0}$. Here $R_{m0}$ is the optimal $R_m$ computed
    by the method of~\cite{DongY2021}, with $R_{m0}=0.065$ for the network
    $[2,500,50,1]$ and $R_{m0}=0.65$ for the network $[2, 50, 500, 1]$.
  }
  \label{tab_2}
\end{table}

Table~\ref{tab_2} provides an accuracy comparison of the current HLConcELM
method and the conventional ELM method~\cite{DongL2021} for solving
the advection equation. Two neural networks are considered here, with
the architectures given by $\mbs M_1=[2, 500, 50, 1]$ and
$\mbs M_2=[2, 50, 500, 1]$, respectively.
The maximum/rms errors of both methods on the domain $\Omega$
corresponding to a sequence of collocation points  are listed in the table.
With the network $\mbs M_1$, whose last hidden layer is narrower,
the conventional ELM exhibits only a fair accuracy with increasing
collocation points, with its maximum errors on the order of $10^{-2}$.
In contrast, the current HLConcELM method produces highly accurate
results with the network $\mbs M_1$, with the maximum error reaching the order of
$10^{-8}$ on the larger set of collocation points.
With the network $\mbs M_2$, whose last hidden layer is wider,
both the conventional ELM and the current HLConcELM
produce highly accurate results with increasing number of collocation points.
These observations are consistent with those in the previous
subsection for the variable-coefficient Poisson equation.

\begin{figure}
  \centerline{
    \includegraphics[width=4.5in]{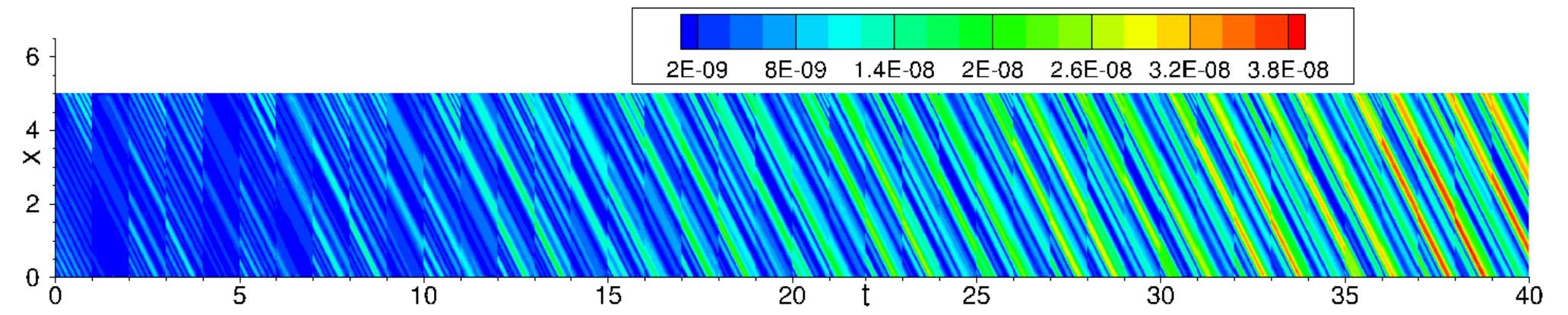}(a)
  }
  \centerline{
    \includegraphics[width=2.in]{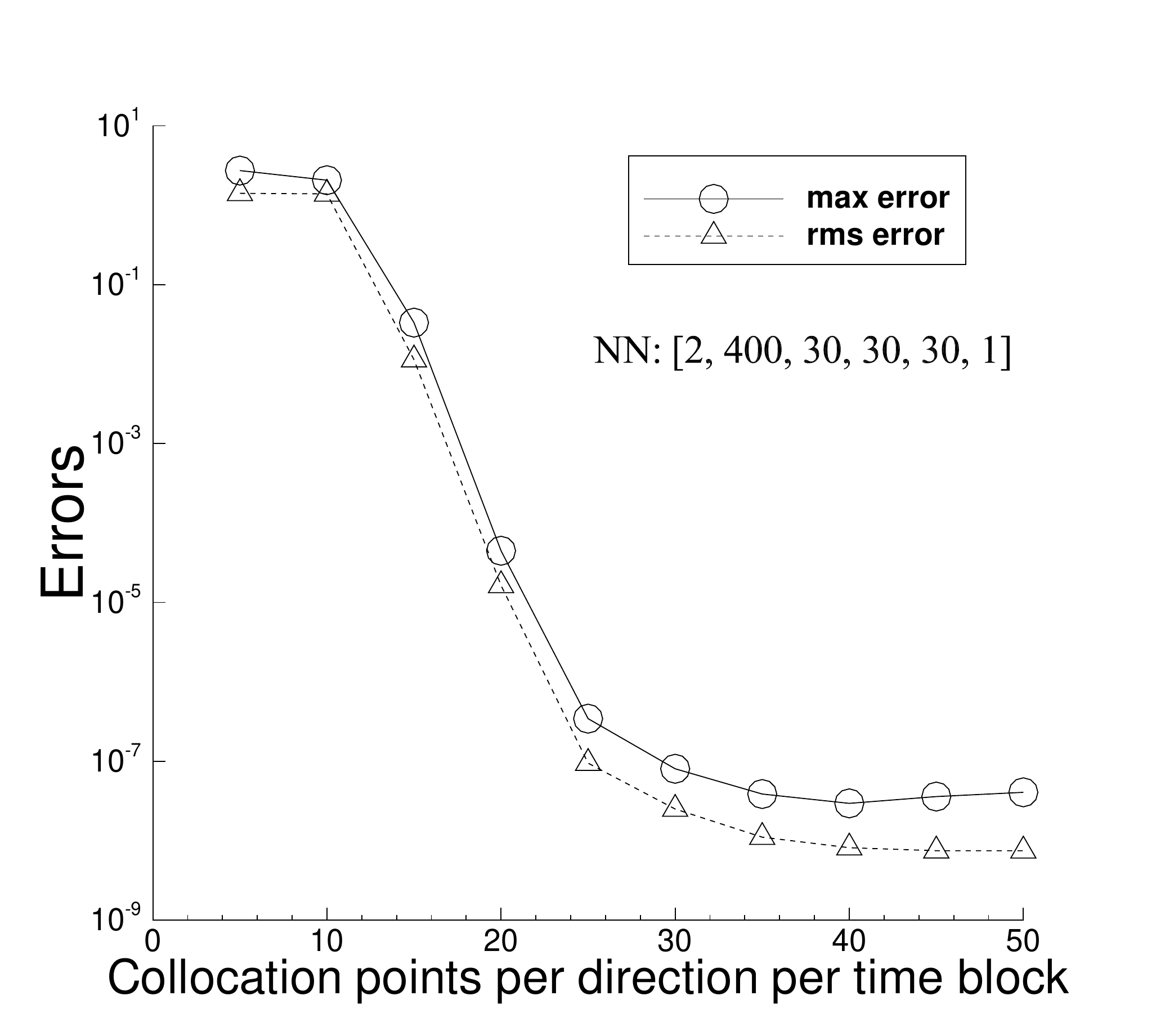}(b)
    \includegraphics[width=2.in]{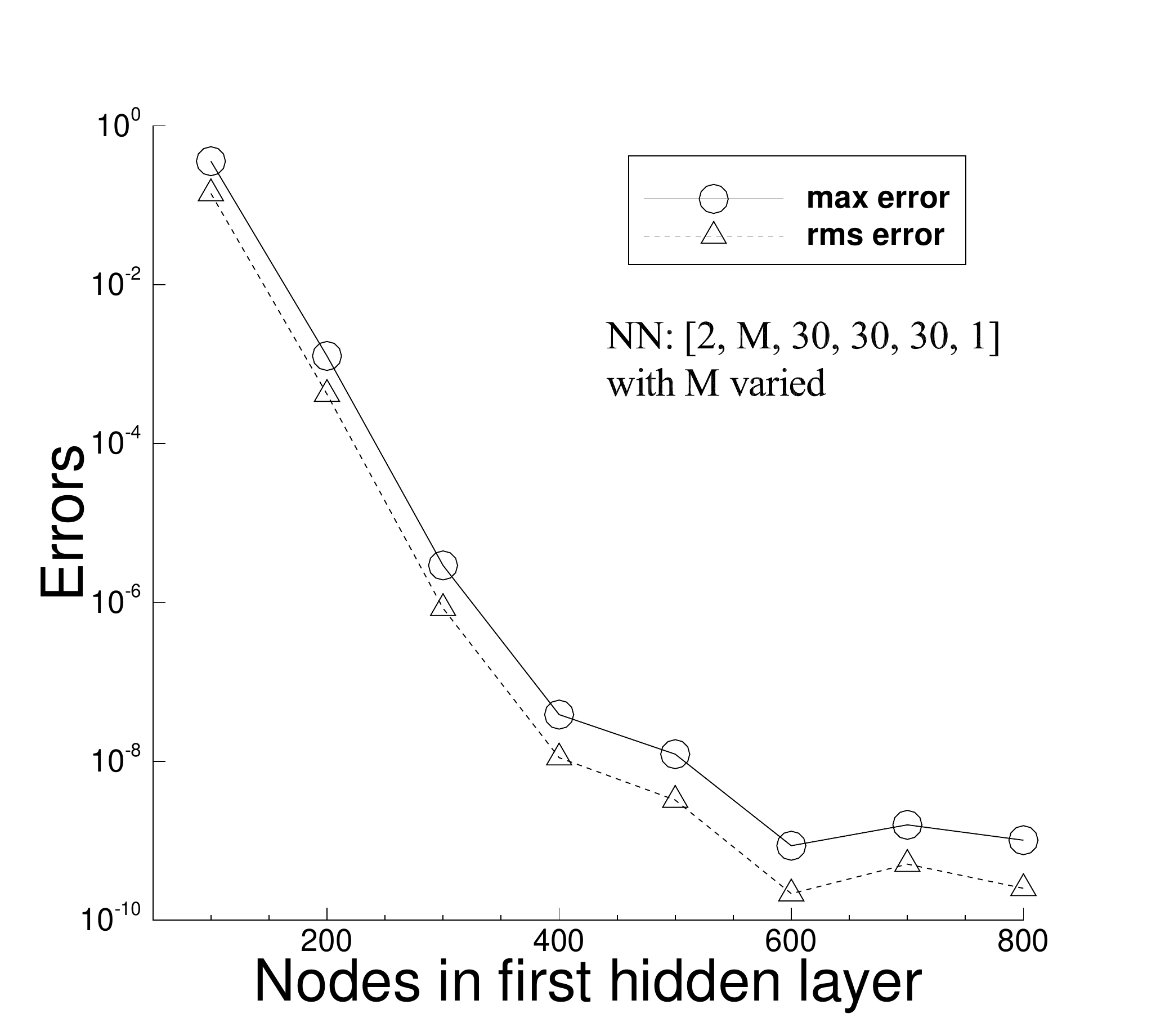}(c)
  }
  \caption{Advection equation (4 hidden layers in neural network):
    (a) Error distribution of the HLConcELM solution on $\Omega$.
    The HLConcELM maximum/rms errors on $\Omega$  versus (b)
    the number of collocation points per direction in each time block,
    and (c) the number of nodes in the first hidden layer ($M$).
    Network architecture $[2, M, 30, 30, 30, 1]$, 
    $40$ time blocks in block time marching.
    $Q=35\times 35$ in (a,c), and is varied in (b).
    $M=400$ in (a,b), and is varied in (c).
    $\mbs R=(3.1,1.0,0.9,0.8)$ in (a,b,c).
  }
  \label{fg_10}
\end{figure}

Figure~\ref{fg_10} illustrates the HLConcELM results obtained on
a deeper neural network containing $4$ hidden layers for solving
the advection equation.
The network architecture is given by $\mbs M=[2, M, 30, 30, 30, 1]$,
where $M$ is either fixed at $M=400$ or varied systematically between
$M=100$ and $M=800$.
A uniform set of $Q=Q_1\times Q_1$ collocation points is used to train
the network, where $Q_1$ is either fixed at $Q_1=35$ or varied systematically
between $Q_1=5$ and $Q_1=50$.
In all simulations we employ a hidden magnitude vector $\mbs R=(3.1,1.0,0.9,0.8)$,
which is computed using the method of~\cite{DongY2021}.
Figure~\ref{fg_10}(a) depicts the distribution of the absolute
error of the HLConcELM solution on $\Omega$, which corresponds to
$M=400$ and $Q_1=35$.
It can be observed that the result is highly accurate, with a maximum error
on the order of $10^{-8}$ in the domain.
Figure~\ref{fg_10}(b) shows the maximum/rms errors of HLConcELM
as a function of $Q_1$, with a fixed $M=400$ in the tests.
Figure~\ref{fg_10}(c) shows the maximum/rms errors of HLConcELM as a function
of $M$ in the neural network, with a fixed $Q_1=35$ for the collocation points.
The exponential convergence of the HLConcELM errors (before saturation)
is unmistakable.


\subsection{Nonlinear Helmholtz Equation}
\label{sec:nonhelm}

\begin{figure}
  \centerline{
    \includegraphics[width=2.in]{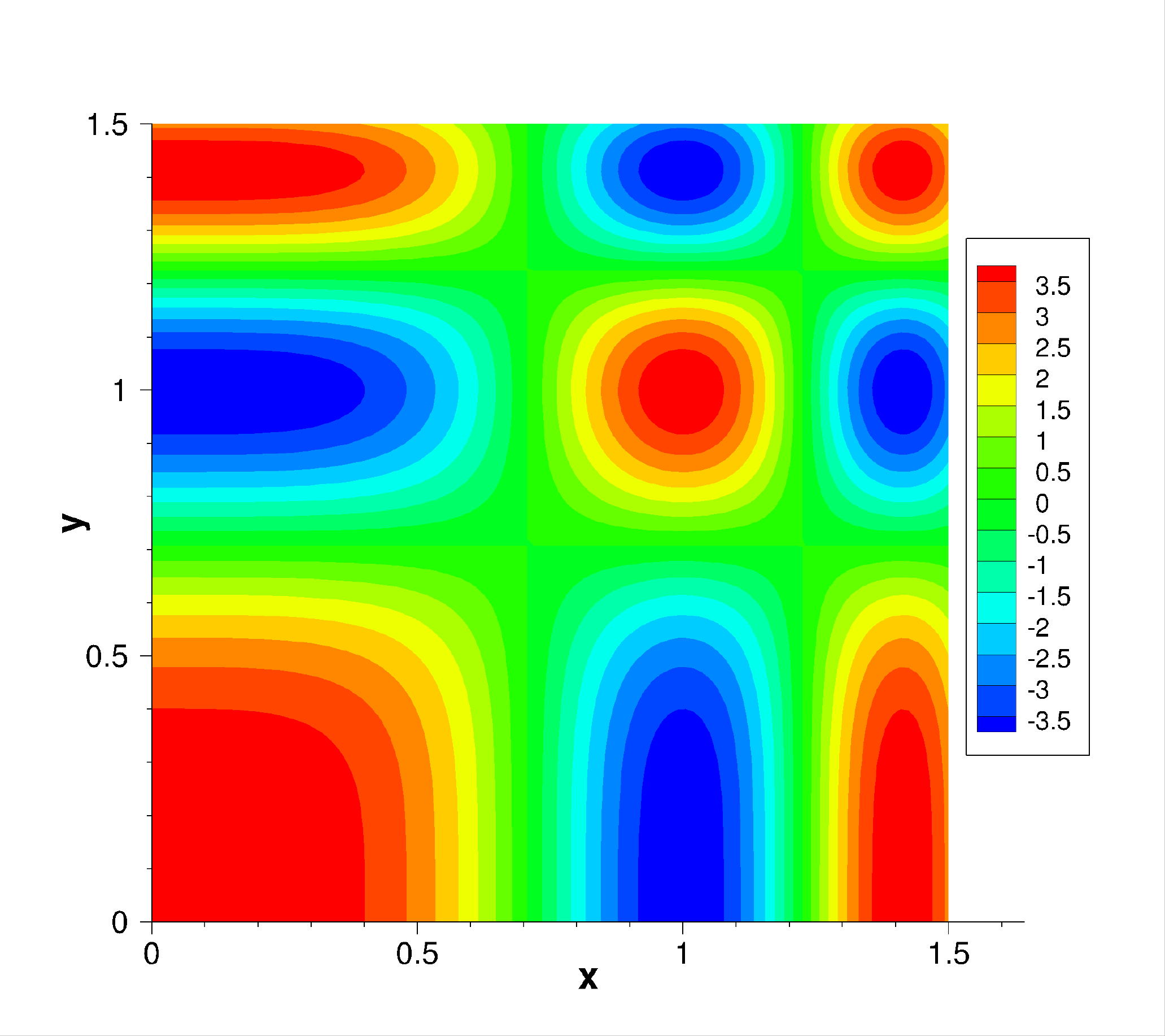}(a)
    \includegraphics[width=2.in]{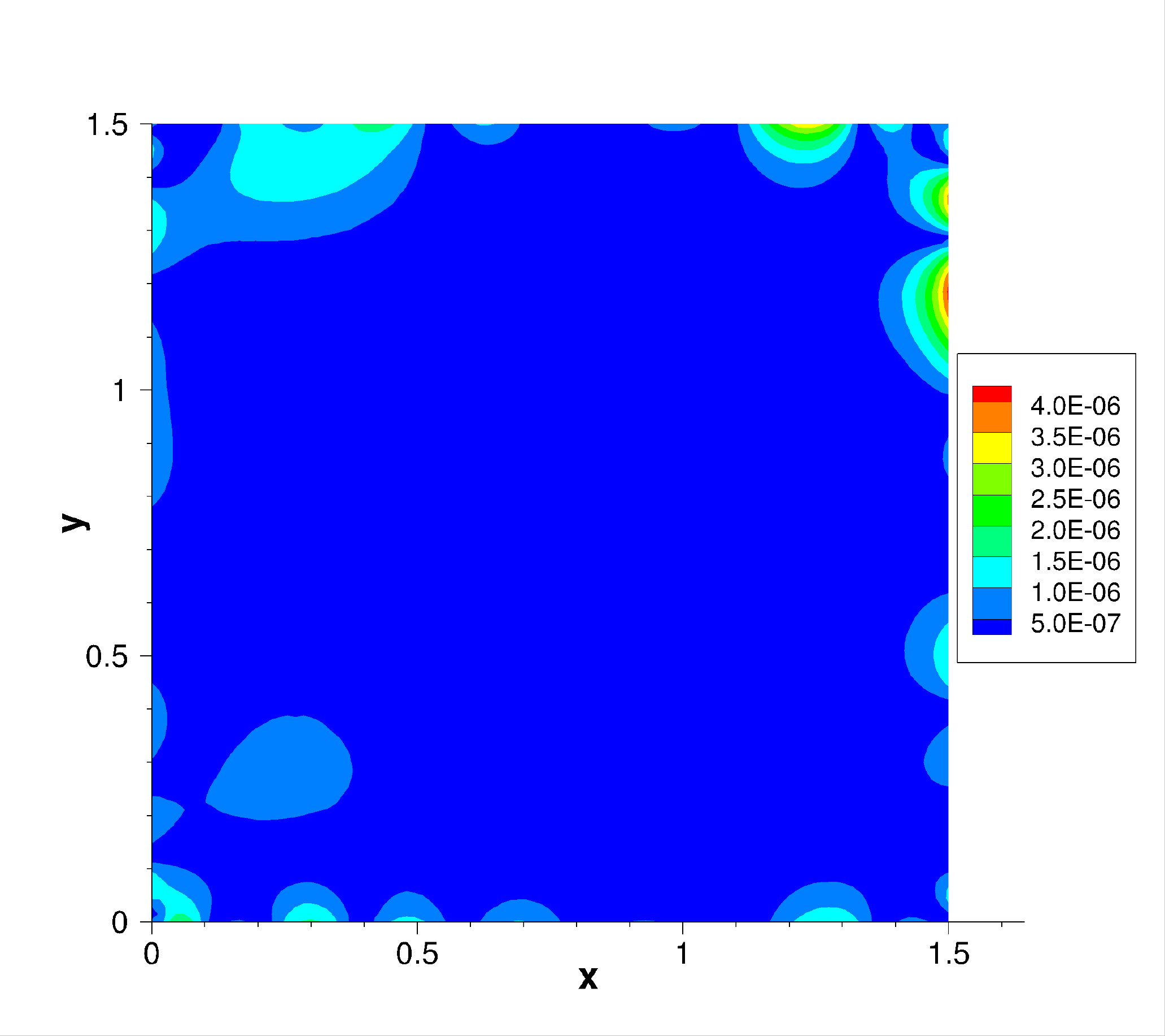}(b)
    \includegraphics[width=2.in]{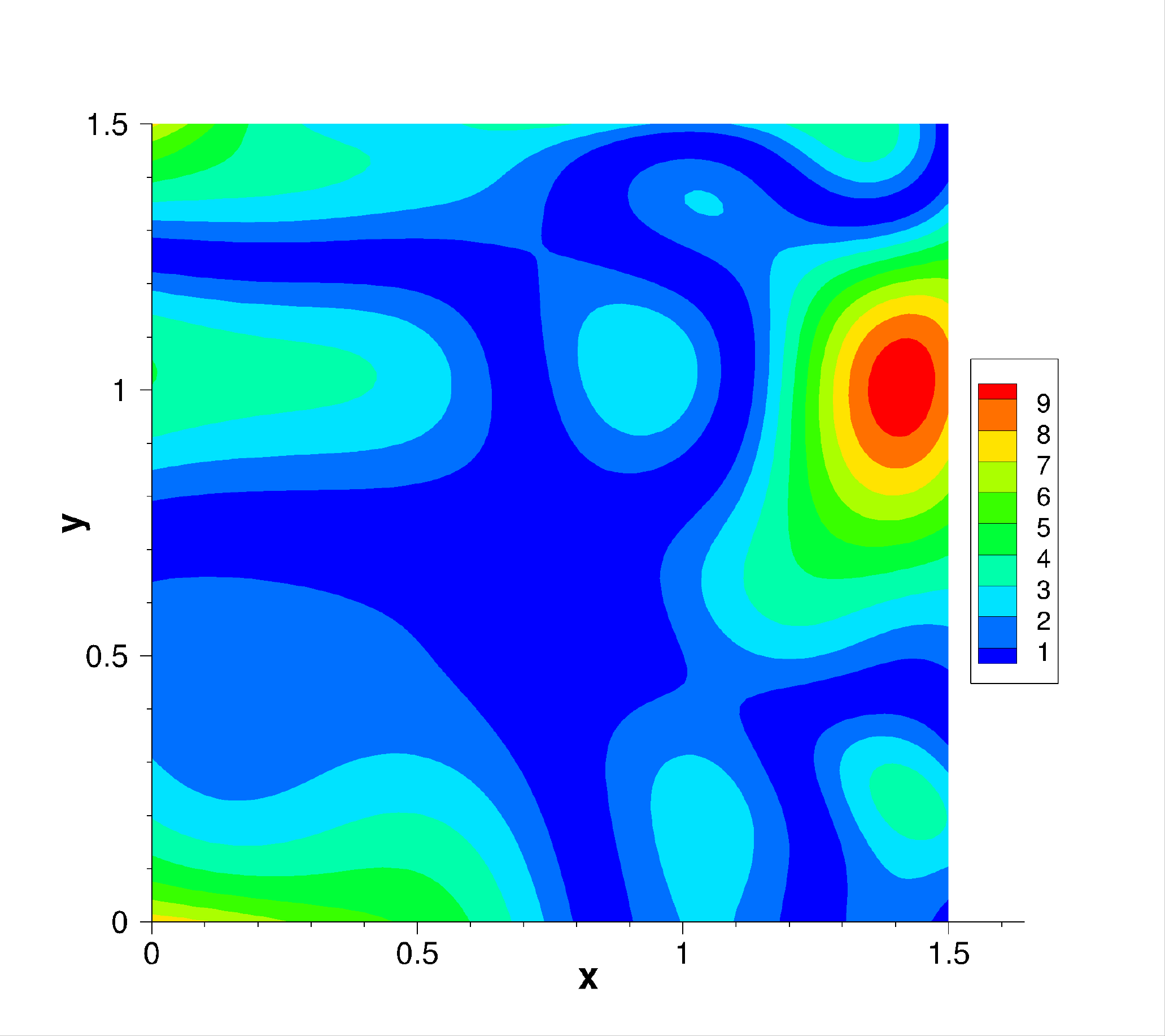}(c)
  }
  \caption{Nonlinear Helmholtz equation:
    Distributions of (a) the exact solution, (b) the absolute error of the HLConcELM
    solution, and (c) the absolute error of the conventional ELM solution.
    In (b,c), network architecture $[2, 500, 30, 1]$, Gaussian activation function,
    $Q=35\times 35$ uniform collocation points.
    $\mbs R=(2.0,3.0)$ in (b) for HLConcELM. $R_m=R_{m0}=0.6$ in (c) for
    conventional ELM.
  }
  \label{fg_11}
\end{figure}

We employ a nonlinear Helmholtz equation to test the HLConcELM
method for the next problem. Consider the 2D domain
$(x,y)\in\Omega = [0,1.5]\times [0, 1.5]$ and the following boundary
value problem on $\Omega$,
\begin{subequations}\label{eq_25}
  \begin{align}
    &
    \frac{\partial^2 u}{\partial x^2} + \frac{\partial^2u}{\partial y^2}
    - 100u + 10\cosh(u) = f(x,y), \\
    &
    u(0,y) = g_1(y),\quad
    u(1.5,y) = g_2(y), \quad
    u(x,0) = h_1(x), \quad
    u(x,1.5) = h_2(x).
  \end{align}
\end{subequations}
In the above equations $u(x,y)$ is the field solution to be sought,
$f(x,y)$ is a prescribed source term, $g_i$ and $h_i$ ($i=1,2$) are
the Dirichlet boundary data.
In this subsection we choose $f$, $g_i$ and $h_i$ ($i=1,2$) such that
the system~\eqref{eq_25} has the following solution,
\begin{equation}\label{eq_26}
  u(x,y) = 4\cos\left(\pi x^2\right)\cos\left(\pi y^2\right).
\end{equation}
The distribution of this solution in the $xy$ plane
is illustrated in Figure~\ref{fg_11}(a).


We employ the HLConcELM method with  neural networks that
contain two input nodes,
representing the $x$ and $y$, and a single output node,
representing the solution $u$. The number of hidden layers and
the number of hidden nodes are varied and will
be specified below. To train the neural network,
we employ a uniform set of $Q=Q_1\times Q_1$ collocation points
on $\Omega$, with $Q_1$ varied in the tests.
The ELM errors reported below are computed on a finer
set of $Q_{eval}=101\times 101$ uniform grid points, as explained before.

Figures~\ref{fg_11}(b) and (c) illustrate the absolute-error
distributions obtained using the HLConcELM method and
the conventional ELM method with the network
architecture $\mbs M=[2, 500, 30, 1]$.
A uniform set of $Q=35\times 35$ collocation points
has been used to train the network with both methods.
The hidden magnitude vector is $\mbs R=(2.0,3.0)$ for HLConcELM,
which is obtained with the method of~\cite{DongY2021}.
For conventional ELM we have employed $R_m=R_{m0}=0.6$, which
is also obtained using the method of~\cite{DongY2021}, for
generating the random hidden-layer coefficients.
The conventional ELM solution
is inaccurate, with the maximum error on order of $10$.
On the other hand, the current HLConcELM method produces
an accurate solution on the same network architecture,
with the maximum error on the order of $10^{-6}$ in the domain.

\begin{figure}
  \centerline{
    \includegraphics[width=2.in]{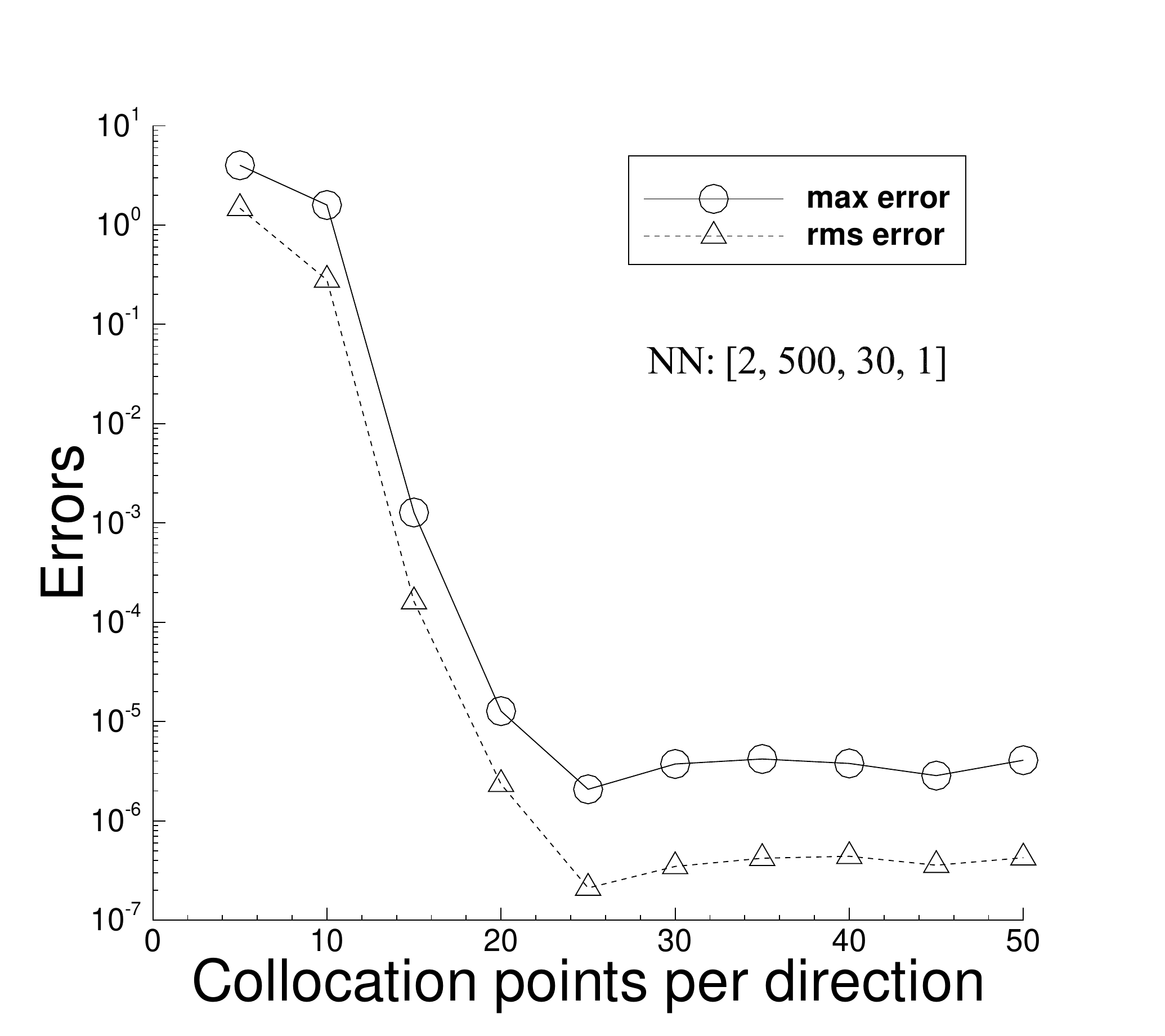}(a)
    \includegraphics[width=2.in]{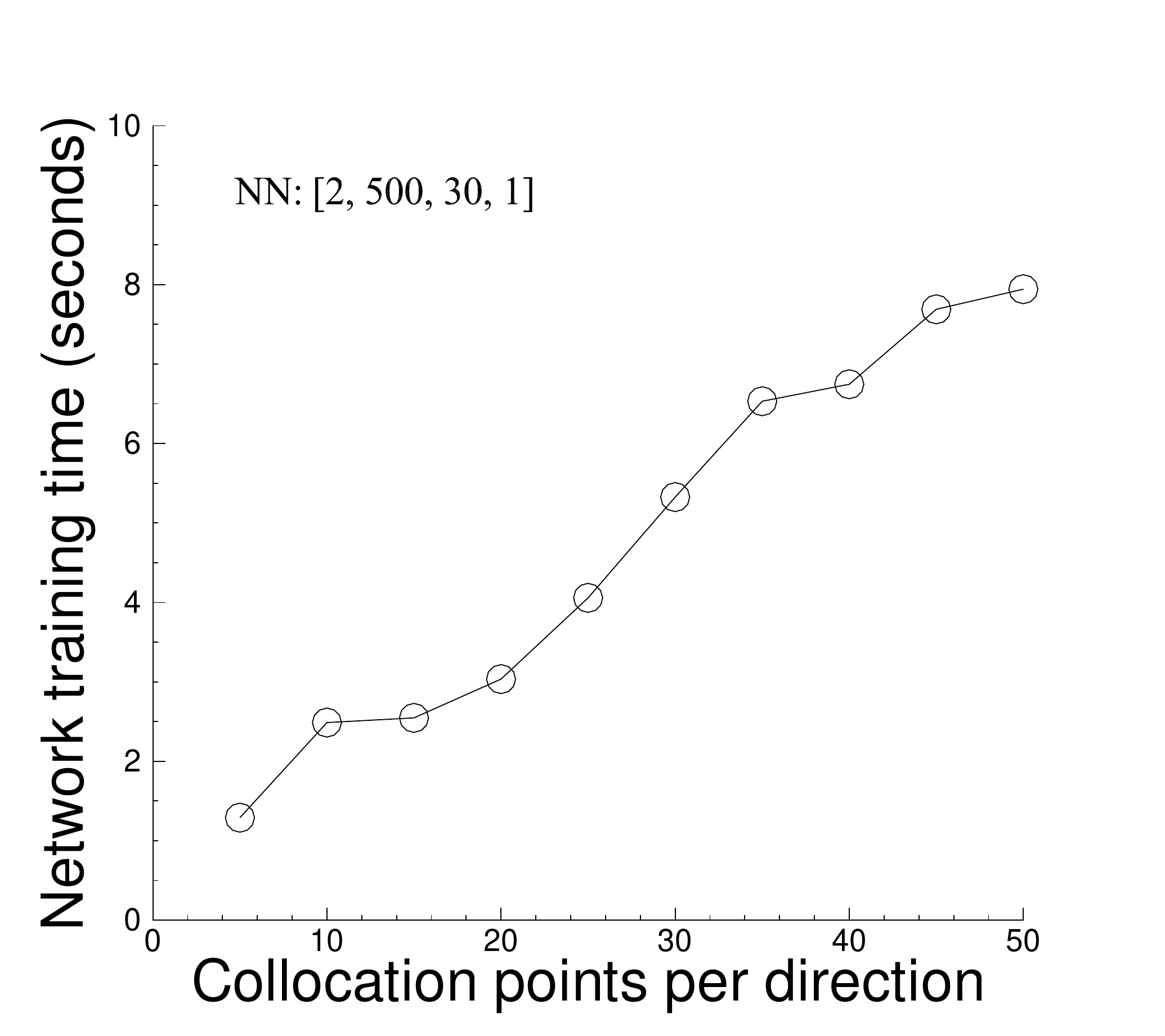}(b)
  }
  \centerline{
    \includegraphics[width=2.in]{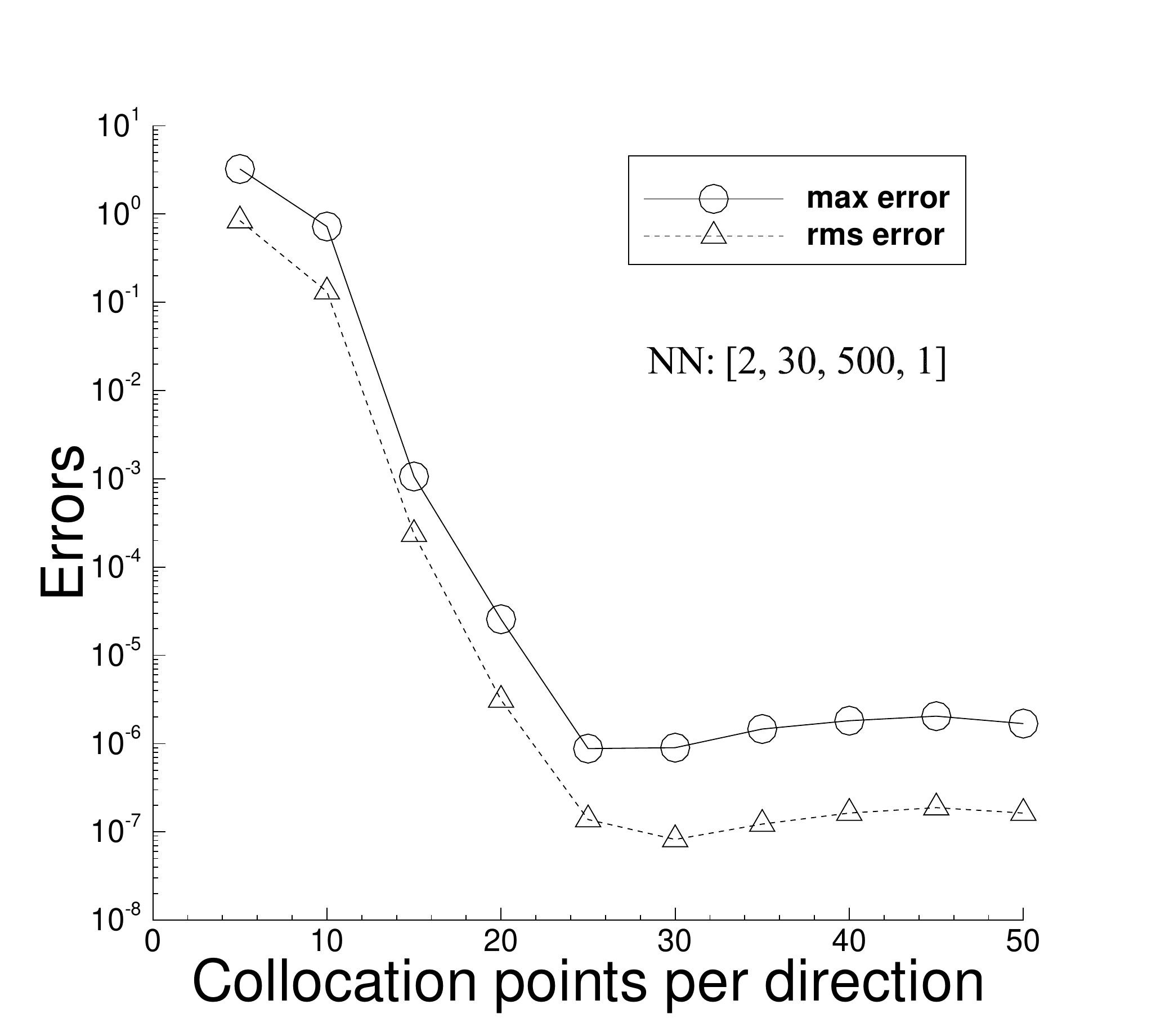}(c)
    \includegraphics[width=2.in]{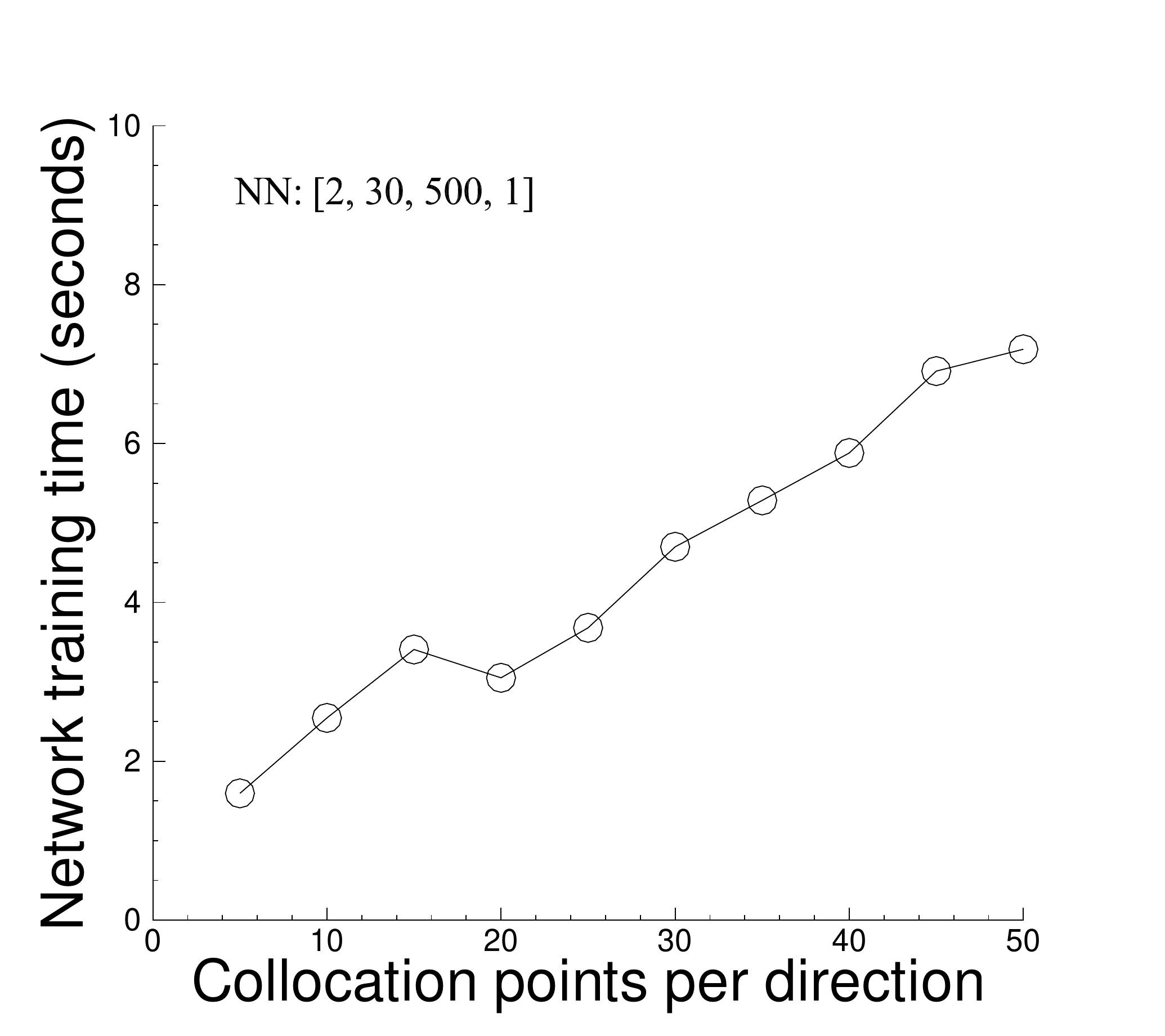}(d)
  }
  \caption{Nonlinear Helmholtz equation:
    The maximum/rms errors (a,c) and the network training
    time (b,d) of the HLConcELM method versus the number of collocation
    points in each direction.
    In (a,b), network architecture $[2, 500, 30, 1]$, $\mbs R=(2.0,3.0)$.
    In (c,d), network architecture $[2, 30, 500, 1]$, $\mbs R=(0.65,0.7)$.
    In (a,b,c,d), uniform collocation points $Q=Q_1\times Q_1$,
    with $Q_1$ varied.
  }
  \label{fg_12}
\end{figure}

Figure~\ref{fg_12} illustrates the convergence behavior and the network training time
with respect to the training collocation points of the HLConcELM method
for solving the nonlinear Helmholtz equation.
Two network architectures are considered here, $\mbs M_1=[2, 500, 30, 1]$
and $\mbs M_2=[2, 30, 500, 1]$.
The number of collocation points in each direction ($Q_1$) is
varied systematically between $Q_1=5$ and $Q_1=50$ in these tests.
We employ $\mbs R=(2.0,3.0)$ for the network $\mbs M_1$ and
$\mbs R=(0.65,0.7)$ for the network $\mbs M_2$, which are obtained
using the method of~\cite{DongY2021}.
Figures~\ref{fg_12}(a) and (b) show the maximum/rms errors and the network
training time of the HLConcELM method as a function of $Q_1$ for the neural
network $\mbs M_1$.
Figures~\ref{fg_12}(c) and (d) show the corresponding results for the
network $\mbs M_2$.
The exponential convergence (before saturation) and the near linear growth
in the network training time observed here for the nonlinear Helmholtz equation
are consistent with
those for the linear problems in previous subsections.

\begin{figure}
  \centerline{
    \includegraphics[width=2.in]{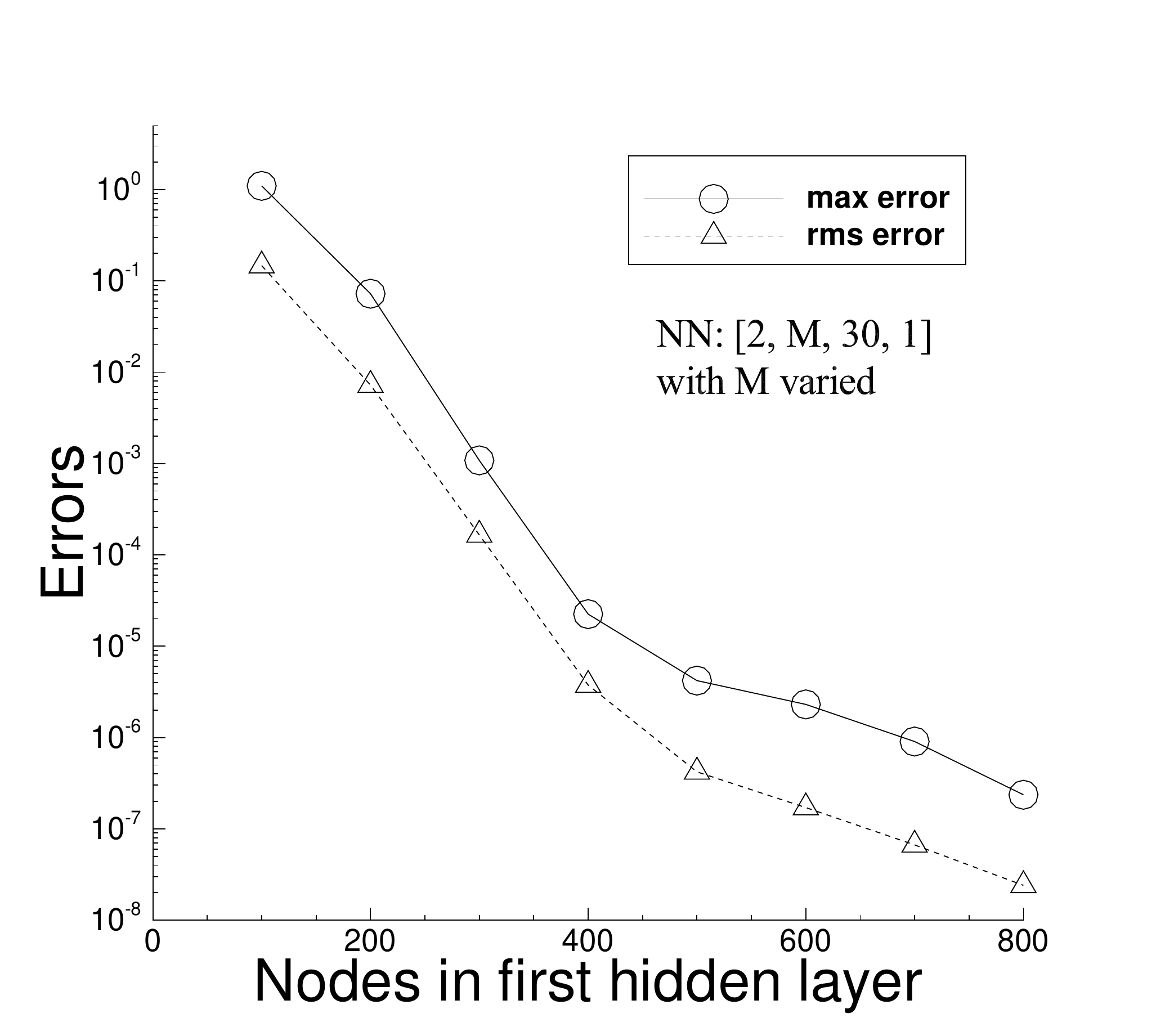}(a)
    \includegraphics[width=2.in]{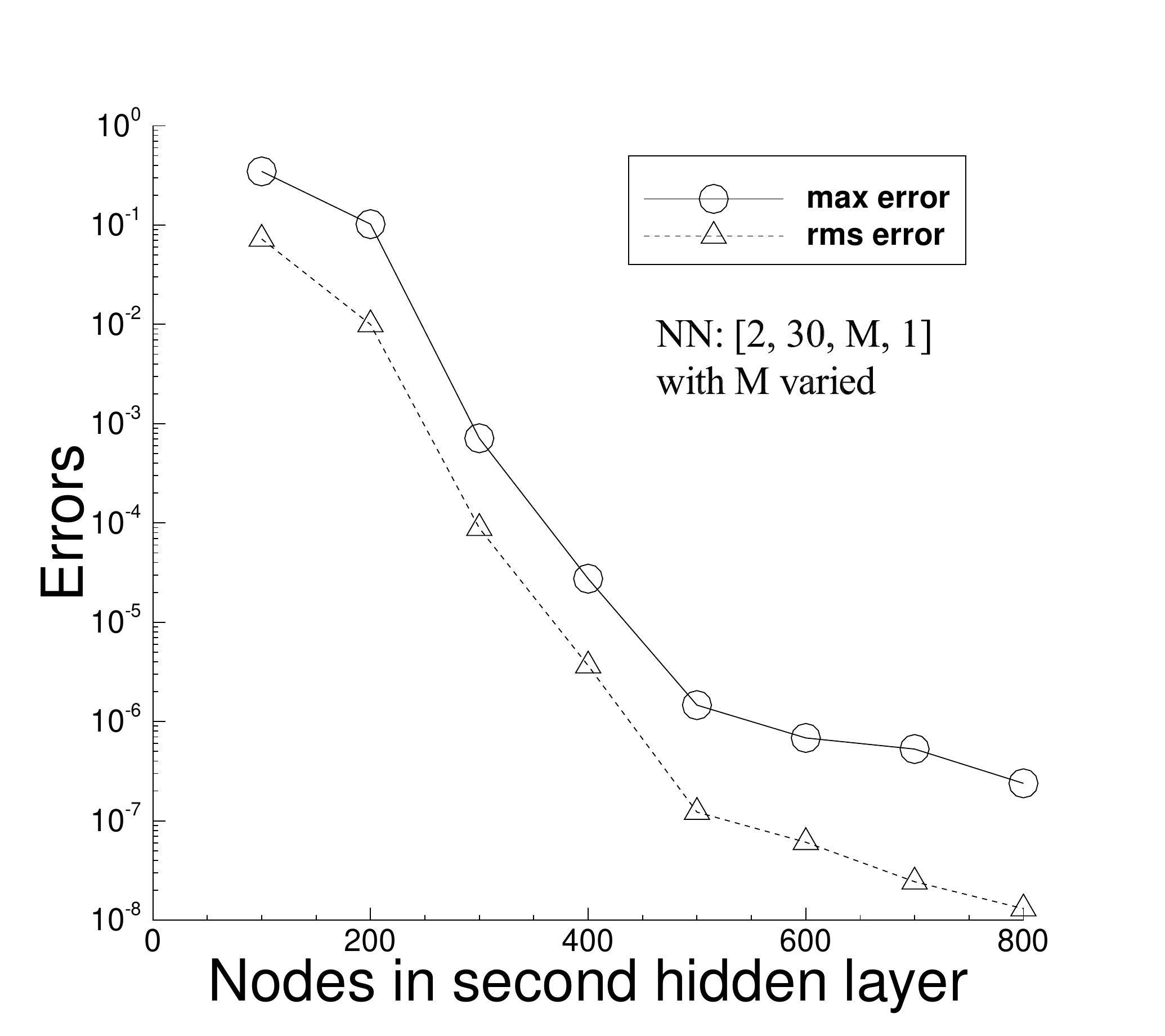}(b)
  }
  \caption{Nonlinear Helmholtz equation: (a) The HLConcELM maximum/rms errors
    versus the number of nodes in
    the first hidden layer with network architecture
    $[2, M, 30, 1]$ ($M$ varied).
    (b) The HLConcELM maximum/rms errors
    versus the number of nodes in
    the second hidden layer with the architecture
    $[2, 30, M, 1]$ ($M$ varied).
    $Q=35\times 35$ collocation points in (a,b).
    $\mbs R=(2.0,3.0)$ in (a), and $\mbs R=(0.65,0.7)$ in (b).
  }
  \label{fg_13}
\end{figure}

Figure~\ref{fg_13} illustrates the error convergence of the HLConcELM
method with respect to the number of nodes in the neural network.
Two groups of neural networks are considered here, with
the architectures $\mbs M_1=[2, M, 30, 1]$ and
$\mbs M_2=[2, 30, M, 1]$, where $M$ is varied  systematically. 
The networks are trained on a uniform set of $Q=35\times 35$ collocation points.
The plots (a) and (b) show the maximum/rms errors of HLConcELM
as a function of $M$ for these two groups of neural networks.
It can be observed that the errors decrease approximately exponentially
with increasing $M$.


\begin{table}[tb]
  \centering
  \begin{tabular}{l|l|cc|cc}
    \hline
    network & collocation & current & HLConcELM & conventional & ELM  \\ \cline{3-6}
    architecture & points & max error & rms error & max error & rms error \\ \hline
    $[2,500,30,1]$ & $5\times 5$ & $4.00E+0$ & $1.48E+0$ & $7.64E+0$ & $2.41E+0$  \\
    & $10\times 10$ & $1.59E+0$ & $2.80E-1$ & $9.69E+0$ & $2.73E+0$  \\
    & $15\times 15$ & $1.27E-3$ & $1.62E-4$ & $9.73E+0$ & $2.71E+0$ \\
    & $20\times 20$ & $1.27E-5$ & $2.34E-6$ & $9.74E+0$ & $2.73E+0$  \\
    & $25\times 25$ & $2.08E-6$ & $2.11E-7$ & $9.74E+0$ & $2.73E+0$ \\
    & $30\times 30$ & $3.74E-6$ & $3.48E-7$ & $9.75E+0$ & $2.74E+0$  \\ \hline
    $[2,30,500,1]$ & $5\times 5$ & $3.23E+0$ & $8.43E-1$ & $3.80E+0$ & $1.15E+0$  \\
    & $10\times 10$ & $7.22E-1$ & $1.32E-1$ & $3.08E+0$ & $7.48E-1$  \\
    & $15\times 15$ & $1.06E-3$ & $2.36E-4$ & $3.86E-4$ & $6.29E-5$ \\
    & $20\times 20$ & $2.56E-5$ & $3.12E-6$ & $3.15E-5$ & $5.67E-6$  \\
    & $25\times 25$ & $8.78E-7$ & $1.38E-7$ & $1.33E-6$ & $2.68E-7$  \\
    & $30\times 30$ & $8.99E-7$ & $8.20E-8$ & $1.76E-6$ & $1.92E-7$ \\
    \hline
  \end{tabular}
  \caption{Nonlinear Helmholtz equation: Comparison of the maximum/rms errors
    from the HLConcELM method and the
    conventional ELM method~\cite{DongL2021}.
    The HLConcELM data  in this table correspond to a portion of those in
    Figure~\ref{fg_12}(a) for the network $[2, 500, 30, 1]$ and to
    those in Figure~\ref{fg_12}(c) for the network $[2, 30, 500, 1]$.
    For conventional ELM, the hidden-layer coefficients
    are set to uniform random values generated on $[-R_m,R_m]$
    with $R_m=R_{m0}$. Here $R_{m0}$ is the optimal $R_m$ obtained
    using the method of~\cite{DongY2021}, with $R_{m0}=0.6$ for the network
    $[2,500,30,1]$ and $R_{m0}=0.65$ for the network $[2, 30, 500, 1]$.
  }
  \label{tab_3}
\end{table}

Table~\ref{tab_3} compares the numerical errors of the current HLConcELM  method
and the conventional ELM method for solving the nonlinear Helmholtz equation
on two network architectures,
$\mbs M_1=[2, 500, 30, 1]$ and $\mbs M_2=[2, 30, 500, 1]$,
trained on a sequence of uniform sets of collocation points.
The  HLConcELM method produces highly accurate results on
both neural networks. On the other hand, while the conventional ELM
produces accurate results on the network $\mbs M_2$, its solution
on the network $\mbs M_1$ is utterly inaccurate.


\begin{figure}
  \centerline{
    \includegraphics[width=2in]{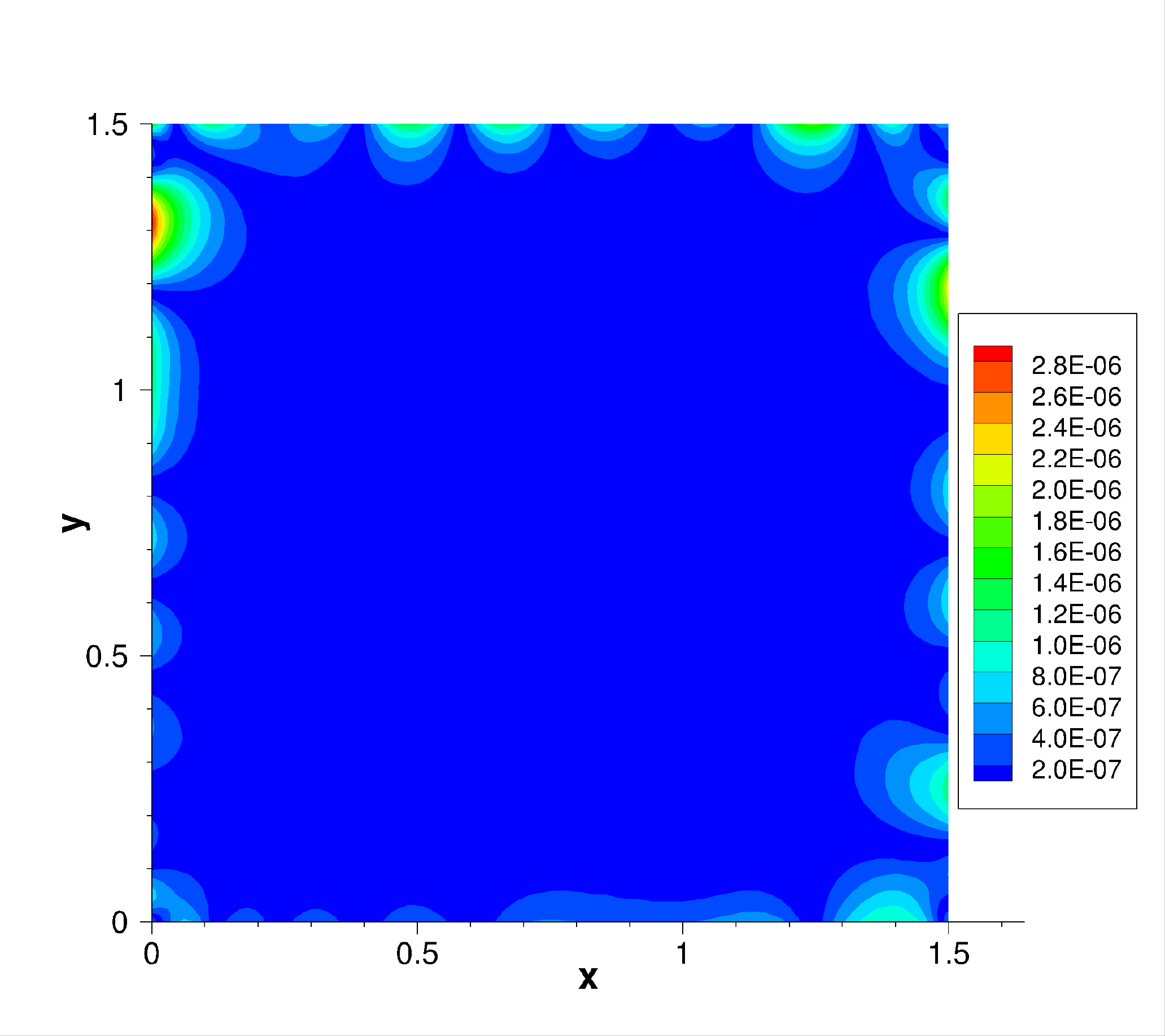}(a)
    \includegraphics[width=2in]{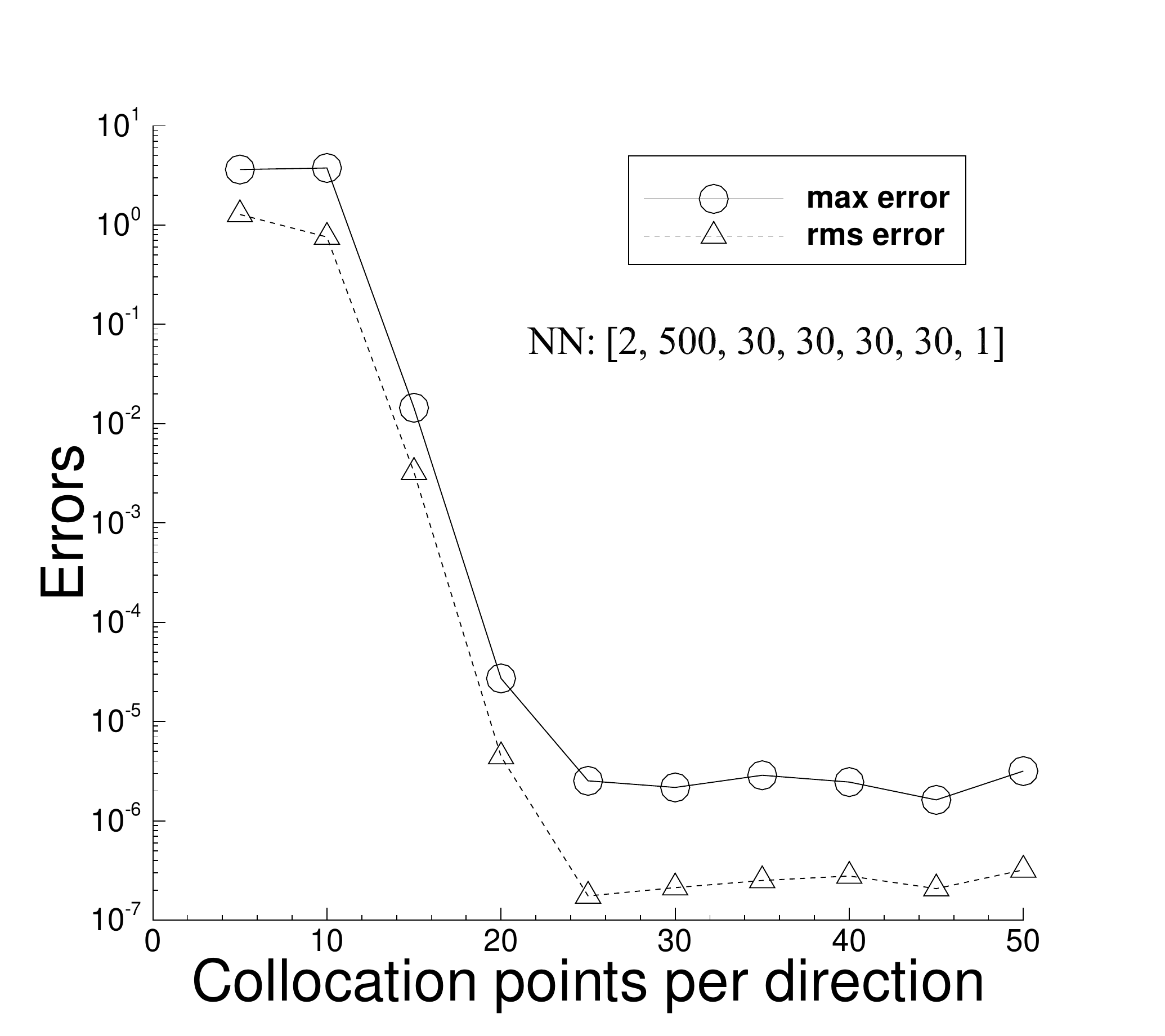}(b)
    \includegraphics[width=2in]{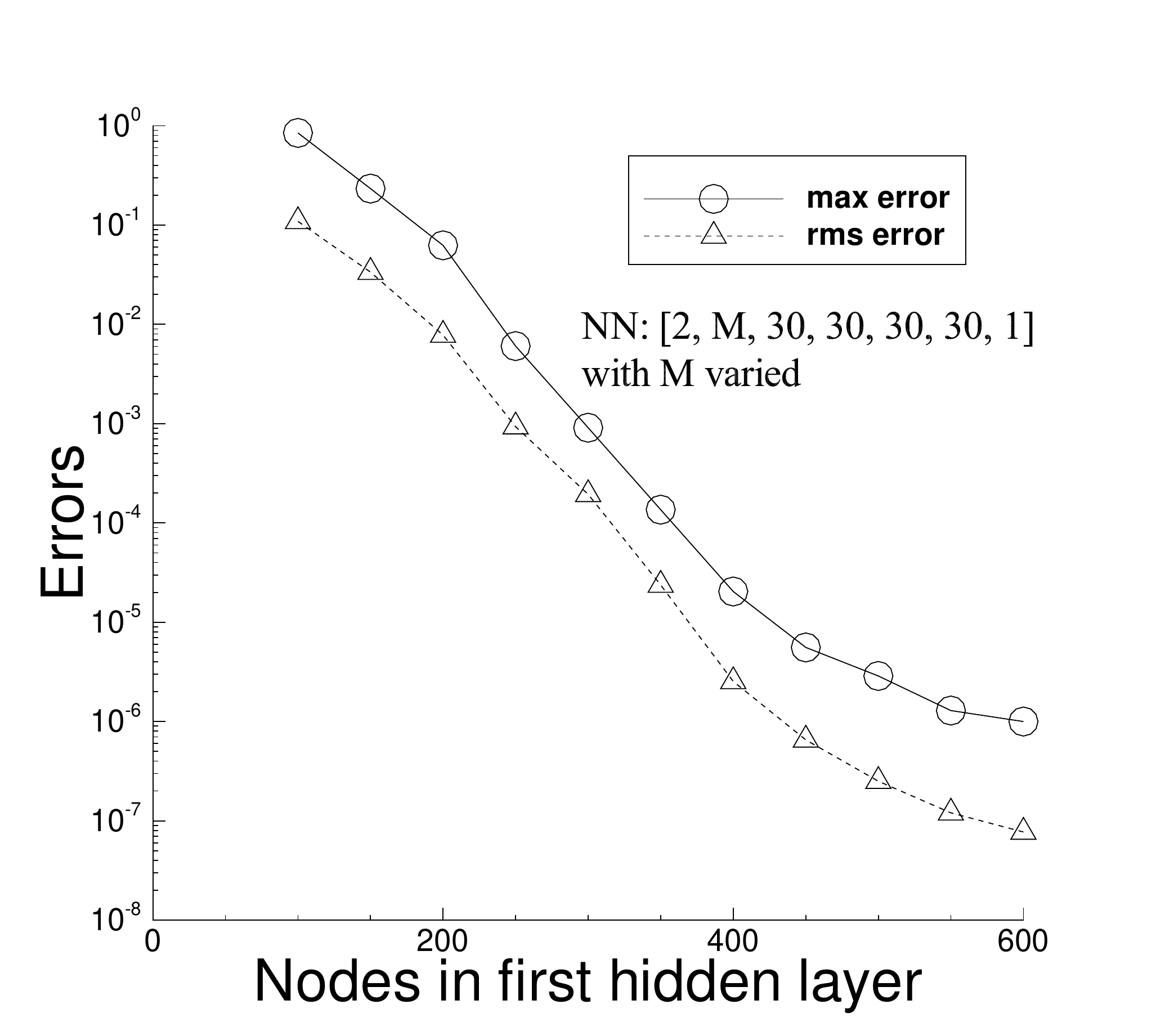}(c)
  }
  \caption{Nonlinear Helmholtz equation (5 hidden layers in neural network):
    (a) Error distribution of
    the HLConcELM solution.
    The maximum/rms errors of the HLConcELM solution versus (b) the
    number of collocation points in each direction and (c) the number of
    nodes in the first hidden layer ($M$).
    Neural network
    architecture $[2, M, 30, 30, 30, 30, 1]$, $Q=Q_1\times Q_1$ uniform collocation points.
    $M=500$ in (a,b), and is varied in (c).
    $Q_1=35$ in (a,c), and is varied in (b).
    $\mbs R=(2.1, 0.1, 2.0, 2.5, 0.5)$ in (a,b,c).
  }
  \label{fg_14}
\end{figure}

Figure~\ref{fg_14} illustrates the HLConcELM results computed on
a deeper neural network with $5$ hidden layers.
The network architecture is given by $\mbs M=[2, M, 30, 30, 30, 30, 1]$,
where $M$ is either fixed at $M=500$ or varied systematically in the tests.
The network is trained on a uniform set of $Q=Q_1\times Q_1$ collocation points,
where $Q_1$ is either fixed at $Q_1=35$ or varied systematically.
Figure~\ref{fg_14}(a) depicts the absolute-error distribution of the HLConcELM
solution obtained with $M=500$ and $Q_1=35$, indicating a quite high accuracy
with the maximum error on the order of $10^{-6}$ in the domain.
Figures~\ref{fg_14}(b) and (c) shows the HLConcELM errors as a function
of $Q_1$ and $M$, respectively. The exponential convergence (prior to saturation)
of these errors is evident. 


\subsection{Burgers' Equation}
\label{sec:burger}

\begin{figure}
  \centerline{
    \includegraphics[width=2.5in]{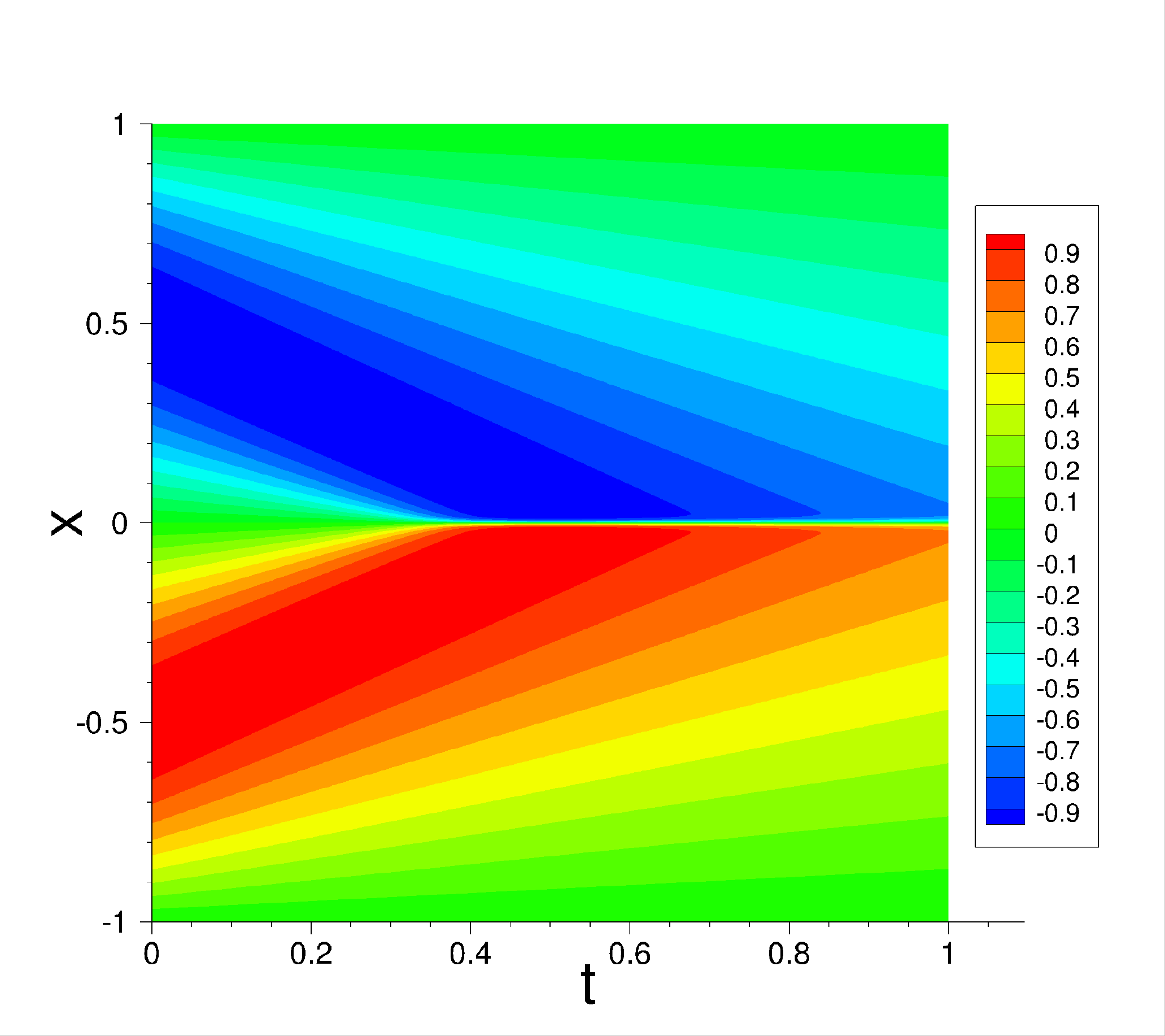}
  }
  \caption{Burgers' equation: distribution of the exact solution.}
  \label{fg_15}
\end{figure}

In the next benchmark example we use the viscous Burgers' equation to test
the performance of the HLConcELM method. Consider the spatial-temporal
domain, $(x,t)\in \Omega=[-1,1]\times [0,1]$, and the following initial/boundary
value problem on $\Omega$,
\begin{subequations}\label{eq_27}
  \begin{align}
    &
    \frac{\partial u}{\partial t} + u\frac{\partial u}{\partial x}
    = \nu\frac{\partial^2 u}{\partial x^2}, \\
    &
    u(-1,t) = u(1,t) = 0, \\
    &
    u(x,0) = -\sin(\pi x),
  \end{align}
\end{subequations}
where $\nu=\frac{1}{100\pi}$, and $u(x,t)$ denotes the field function
to be solved for.
This problem has the following exact solution~\cite{Basdevant1986},
\begin{equation}\label{eq_28}
  u(x,t) = -\frac{\int_{-\infty}^{\infty}\sin\pi(x-\eta)f(x-\eta)e^{-\frac{\eta^2}{4\nu t}}d\eta}
  { \int_{-\infty}^{\infty}f(x-\eta)e^{-\frac{\eta^2}{4\nu t}}d\eta  },
\end{equation}
where $f(y) = e^{-\frac{\cos(\pi y)}{2\pi\nu}}$.
Figure~\ref{fg_15} illustrates the distribution of this solution
on the spatial-temporal domain, which indicates that a sharp gradient
develops in the domain over time.

\begin{figure}
  \centerline{
    \includegraphics[height=2.1in]{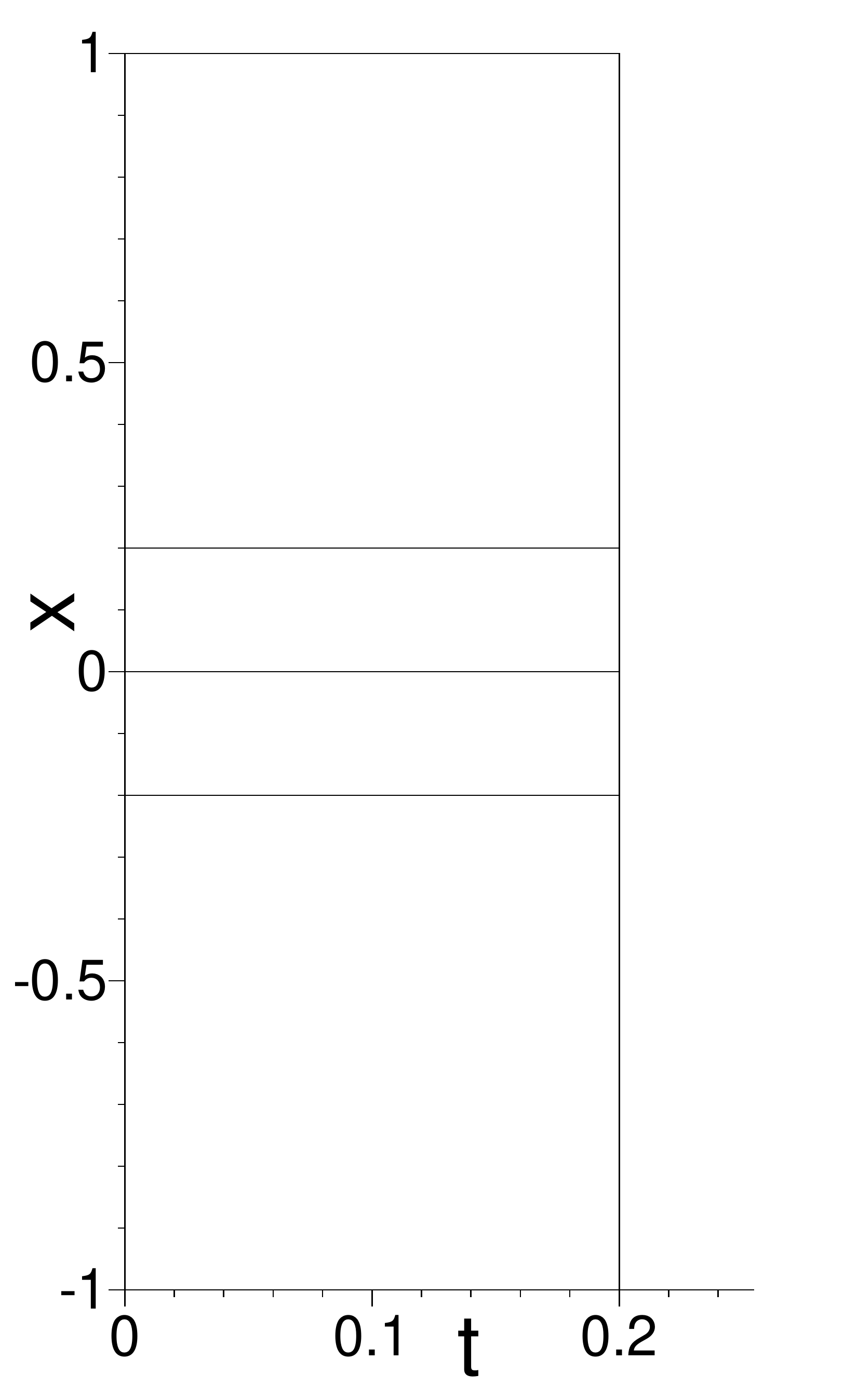}(a)
    \includegraphics[height=2.1in]{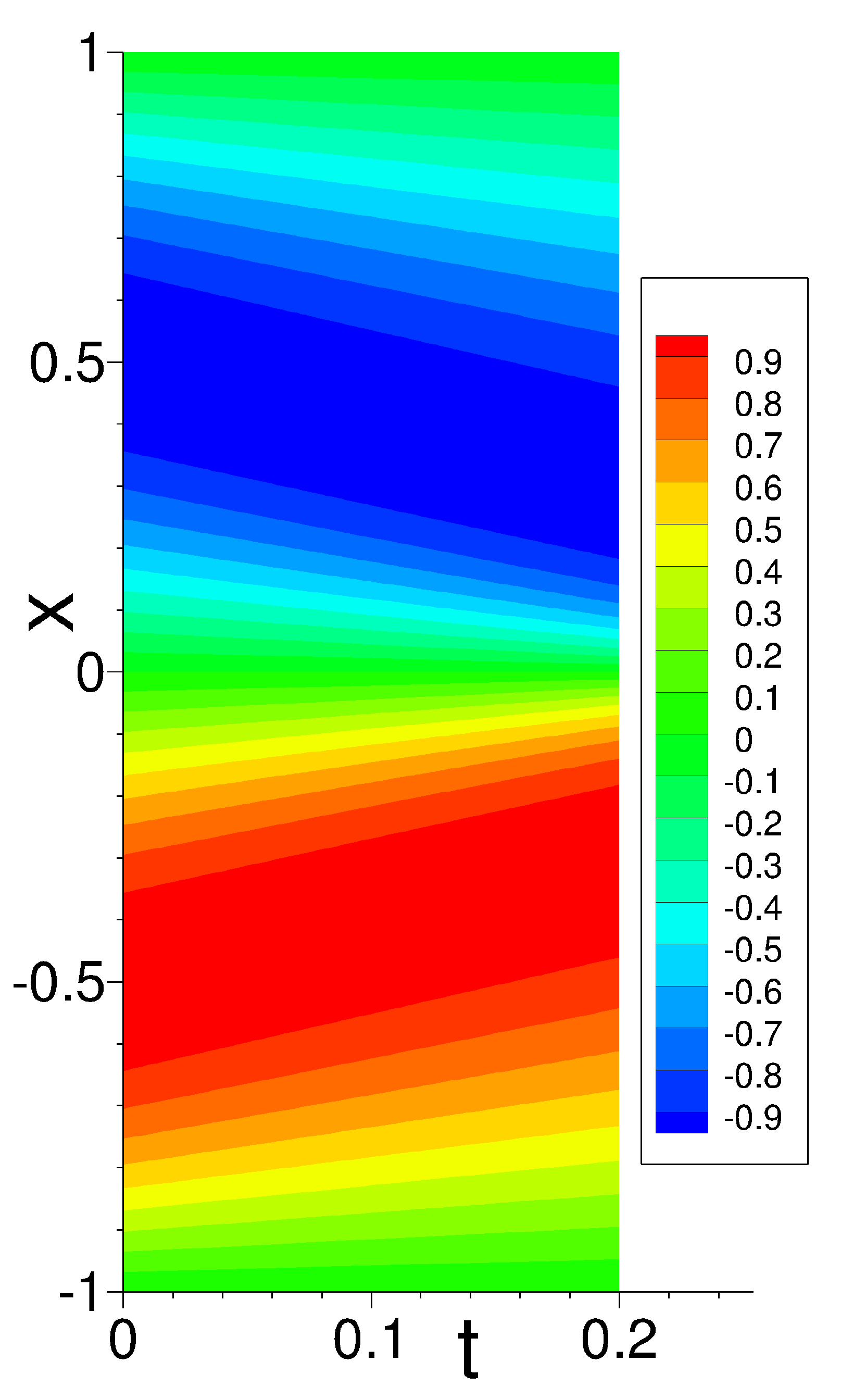}(b)
    \includegraphics[height=2.1in]{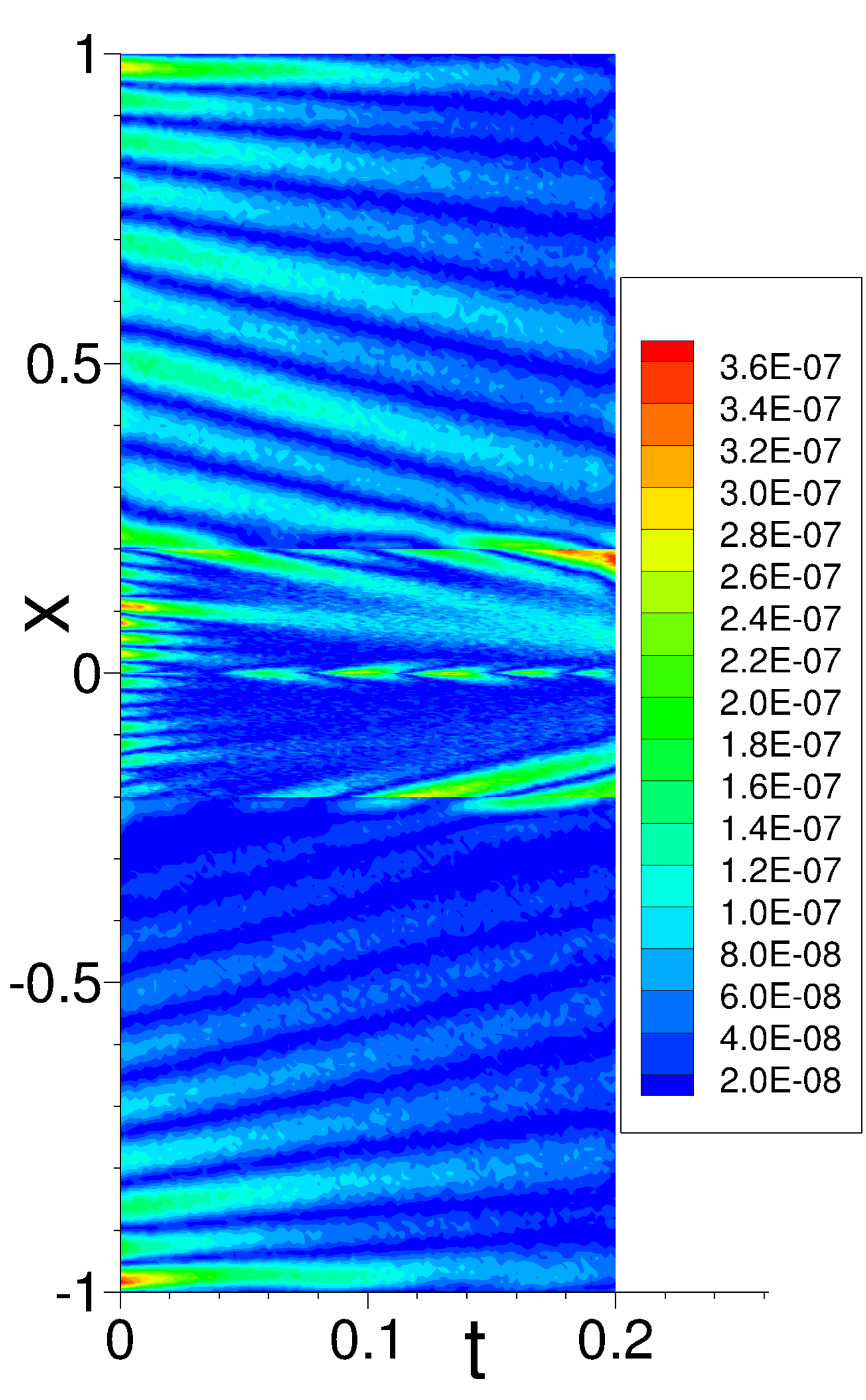}(c)
  }
  \caption{Burgers's equation on the smaller domain $\Omega_1$ ($t\in[0,0.2]$):
    (a) Configuration of the $4$ sub-domains in the locHLConcELM simulation.
    Distributions of the locHLConcELM solution (b)
    and its absolute error (c) on $\Omega_1$.
    Local neural-network architecture: [2, 200, 30, 1],
    $Q=21\times 21$ uniform collocation points per sub-domain,
    $\mbs R=(0.9,0.05)$.
  }
  \label{fg_16}
\end{figure}


We will first solve the problem~\eqref{eq_27} on a smaller domain
(with a smaller temporal dimension)
$\Omega_{1}=[-1,1]\times[0,0.2]$, before the sharp gradient develops,
in order to investigate the convergence behavior of the HLConcELM method.
Then we will compute this problem on the larger domain $\Omega$ using HLConcELM.

On the smaller domain $\Omega_1$ we solve the system~\eqref{eq_27} by
the locHLConcELM method (local version of HLConcELM, see Remark~\ref{rem_4}).
We partition $\Omega_1$ along the $x$ direction into $4$ sub-domains;
see Figure~\ref{fg_16}(a). These sub-domains are non-uniform,
and the $x$ coordinates of the sub-domain boundaries are given
by the vector $\bm{\mathcal{X}}=[-1, -0.2, 0, 0.2, 1]$.
We impose $C^1$ continuity conditions in $x$ across the interior
sub-domain boundaries. We employ a HLConcELM for the local neural network
on each sub-domain, which contains two input nodes (representing
the $x$ and $t$ of the sub-domain) and a single output node (representing
the solution $u$ on the sub-domain). The specific architectures of
the neural network will be provided below.
On each sub-domain we employ a uniform set of $Q=Q_1\times Q_1$ collocation
points ($Q_1$ points in both $x$ and $t$ directions)
for the network training, with $Q_1$ varied in the tests.
We train the overall neural network, which consists of the local neural networks
coupled together by the $C^1$ continuity conditions, by
the nonlinear least squares method;
see Section~\ref{sec:pde} and also~\cite{DongL2021}.

Figures~\ref{fg_16}(b) and (c) illustrate the distributions of
the HLConcELM solution and its absolute error on the domain $\Omega_1$, respectively.
These results are obtained by locHLConcELMs with an
architecture $[2, 200, 30, 1]$ and a uniform set of $Q=21\times 21$
collocation points on each sub-domain.
The hidden magnitude vector is $\mbs R=(0.9, 0.05)$, which is obtained
using the method of~\cite{DongY2021}.
The locHLConcELM method produces an accurate solution, with the maximum
error on the order of $10^{-7}$ on $\Omega_1$.

\begin{figure}
  \centerline{
    \includegraphics[height=2.in]{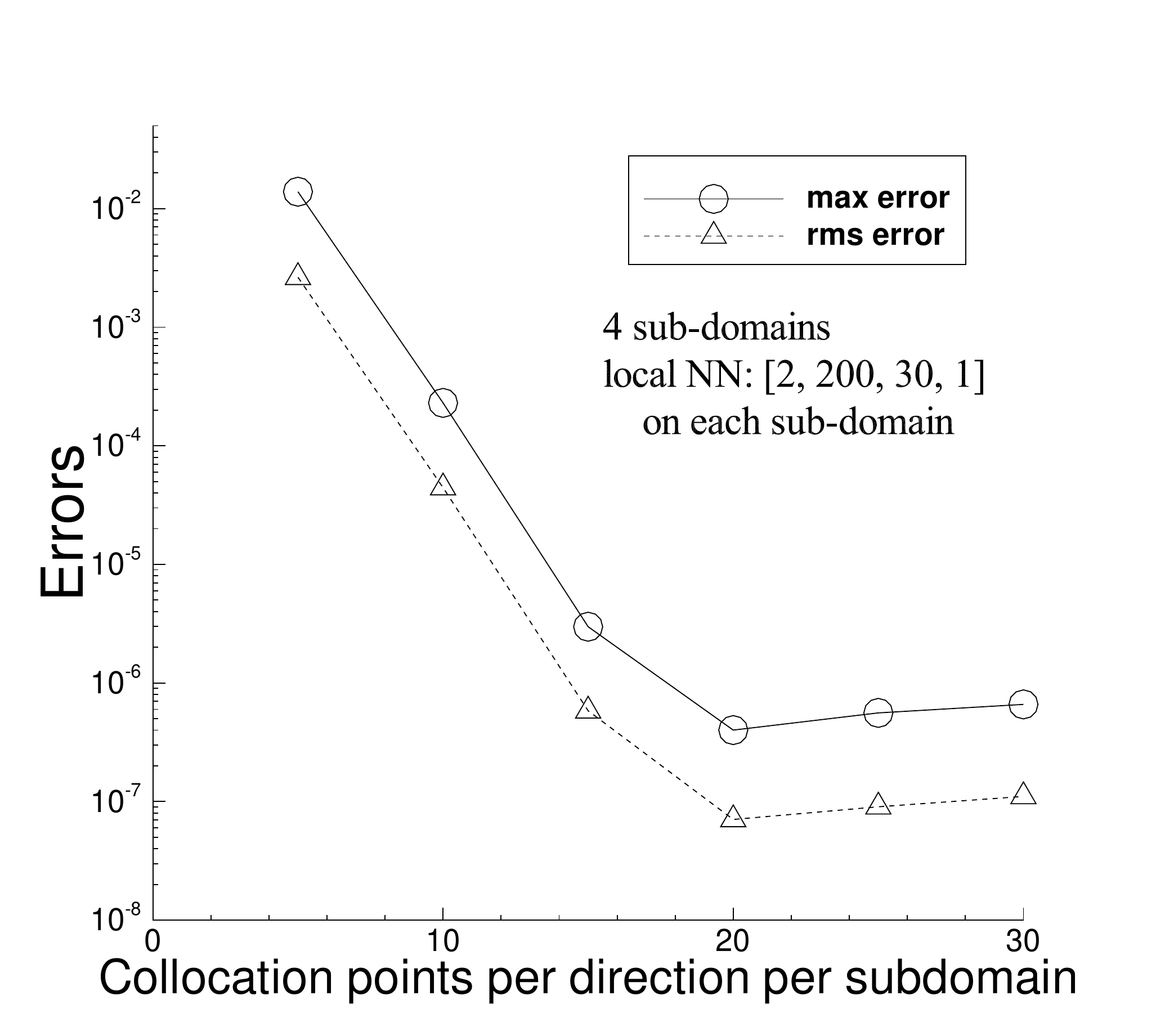}(a)
    \includegraphics[height=2.in]{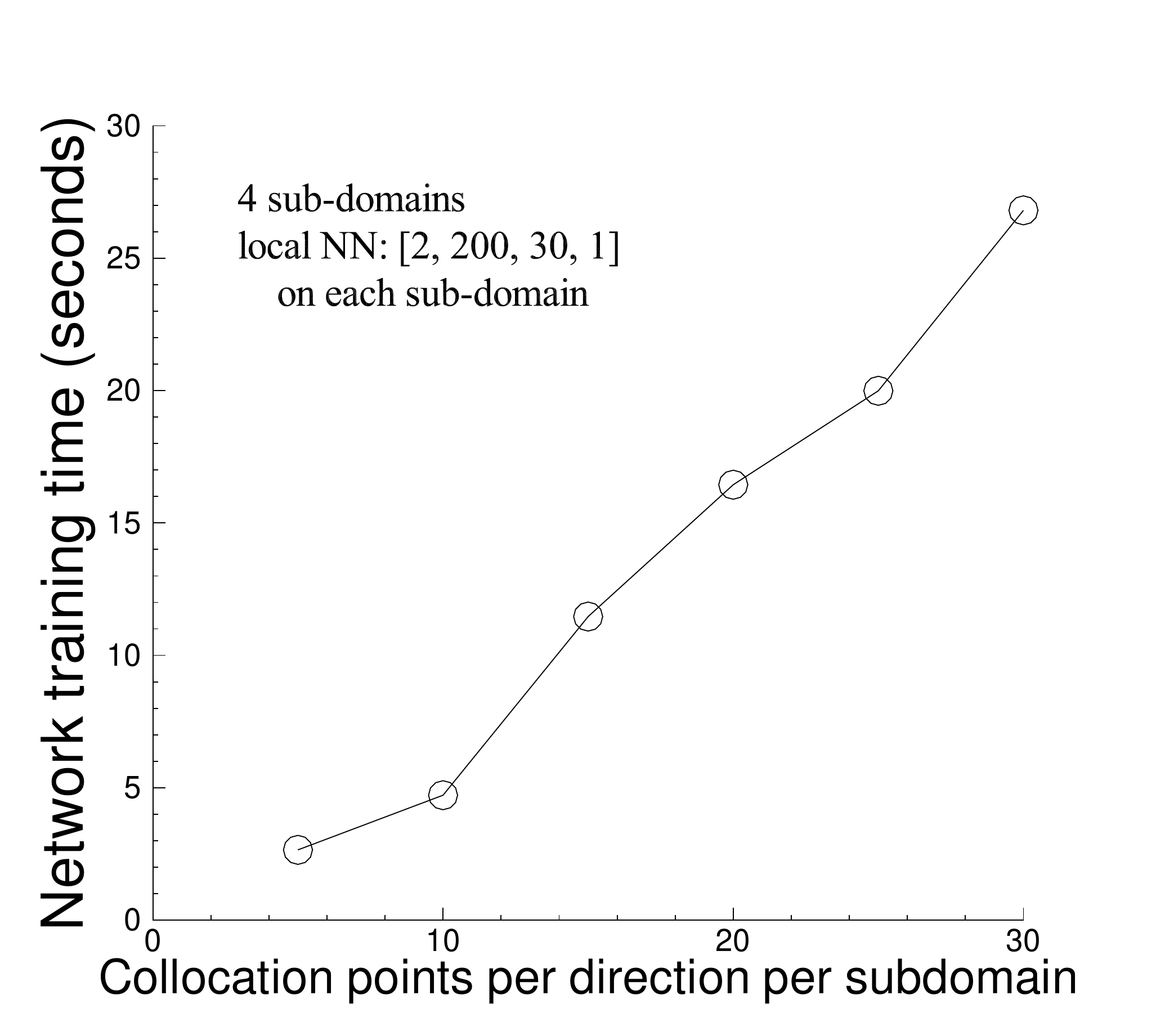}(b)
  }
  \caption{Burgers' equation on the smaller domain $\Omega_1$:
    (a) The locHLConcELM maximum/rms errors and (b) the network
    training time versus the number of collocation points per direction
    in each sub-domain.
    Local network architecture: [2, 200, 30, 1], $Q=Q_1\times Q_1$ uniform
    collocation points ($Q_1$ varied), $\mbs R=(0.9,0.05)$.
  }
  \label{fg_17}
\end{figure}

\begin{figure}
  \centerline{
    \includegraphics[height=2.in]{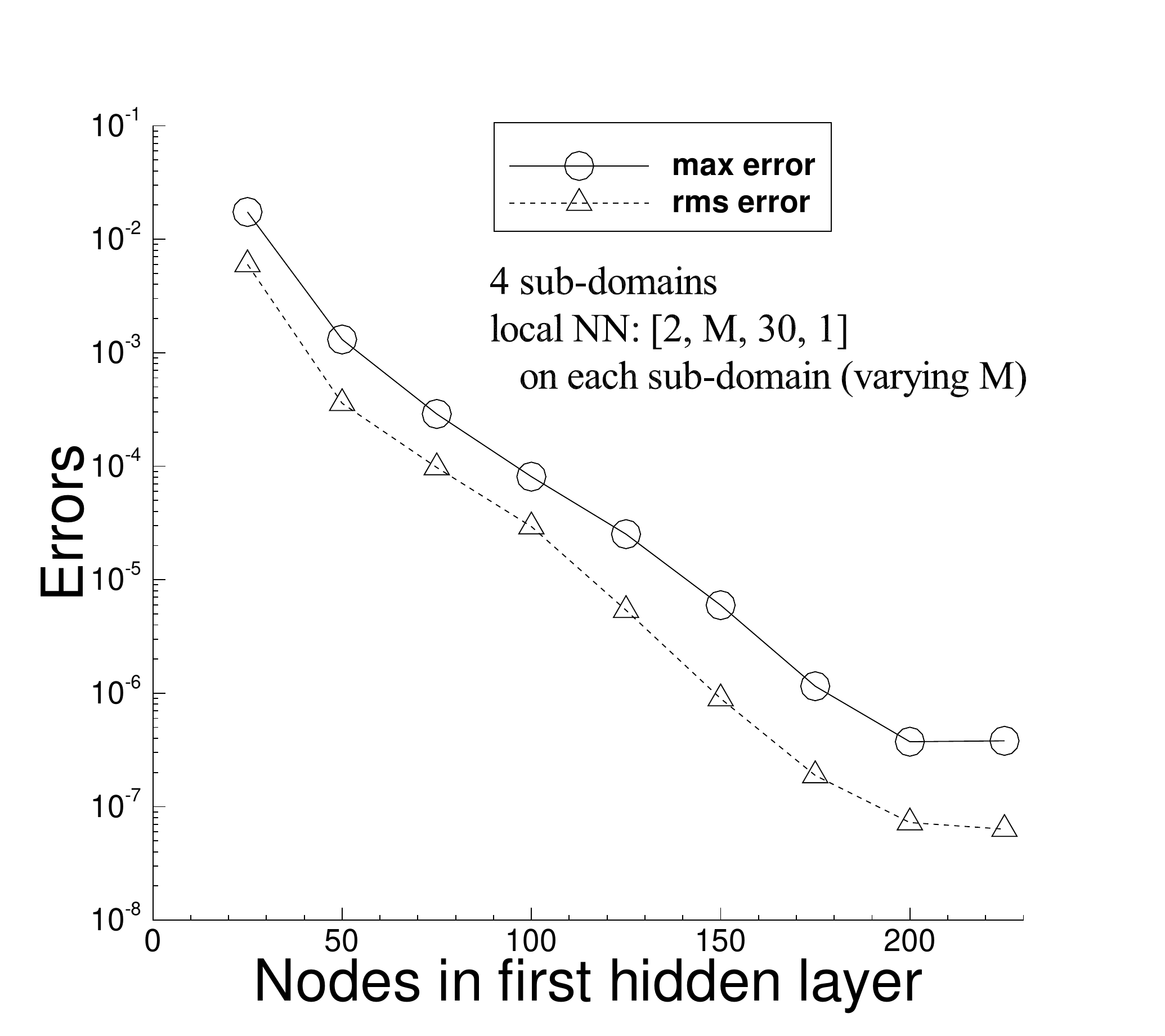}(a)
    \includegraphics[height=2.in]{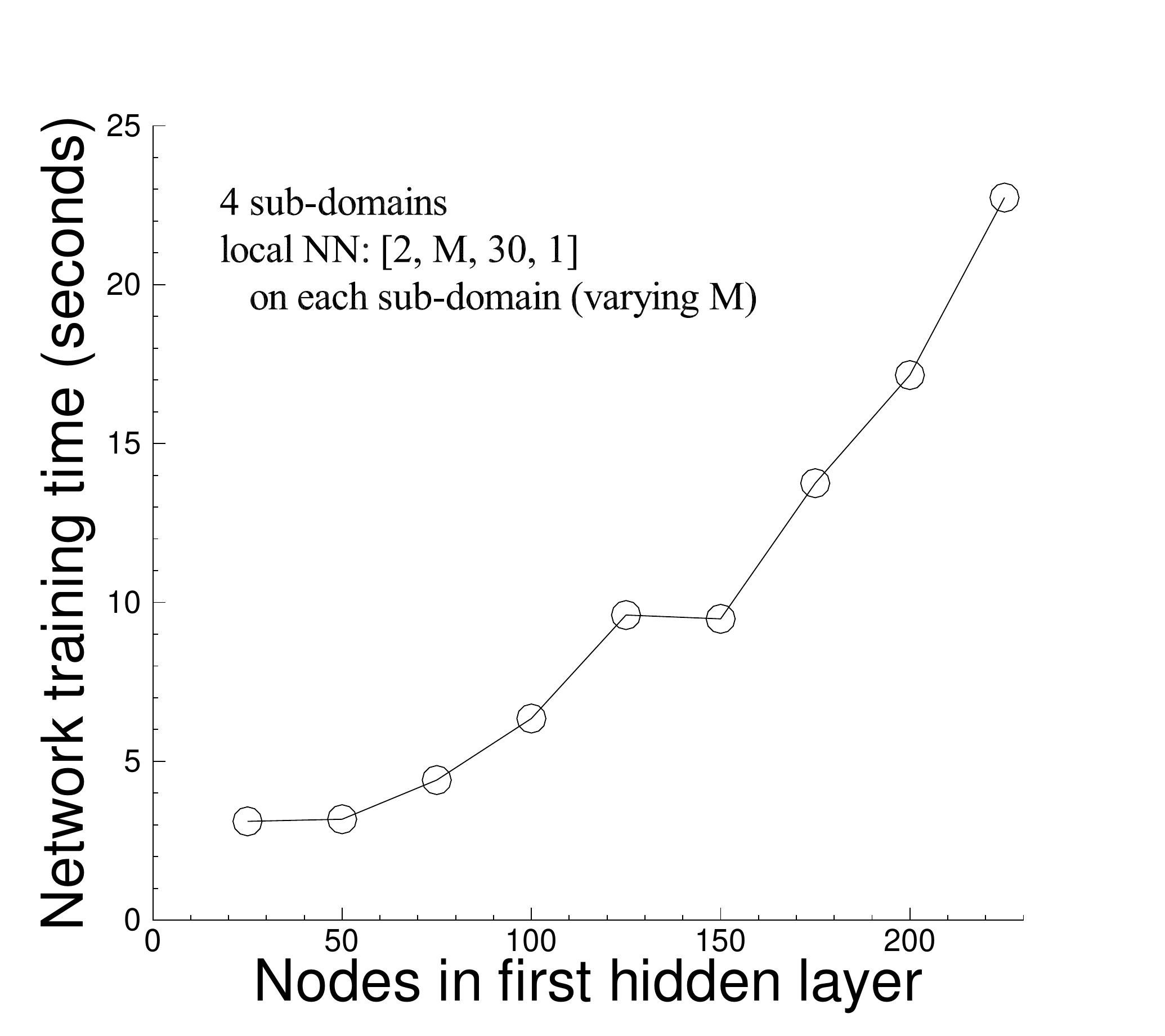}(b)
  }
  \caption{Burgers' equation on the smaller domain $\Omega_1$:
    (a) The locHLConcELM maximum/rms errors and (b) the network
    training time versus the number of nodes in the first hidden layer ($M$).
    Local network architecture [2, M, 30, 1] ($M$ varied), 
    $Q=21\times 21$ uniform collocation points per sub-domain,
    $\mbs R=(0.9,0.05)$.
  }
  \label{fg_18}
\end{figure}

Figure~\ref{fg_17} illustrates the convergence behavior and the network training
time of the locHLConcELM method with respect to the increase of the collocation points
for the smaller domain $\Omega_1$.
The local network architecture is given by $[2, 300, 30, 1]$, and
the collocation points are varied systematically between $Q=5\times 5$
and $Q=30\times 30$ in the tests.
The plots (a) and (b) show the locHLConcELM errors and the network training
time as a function of the number of collocation points in each direction, respectively.
We observe that the locHLConcELM errors decrease exponentially (before saturation)
and the network training time grows approximately linearly with increasing
collocation points.

Figure~\ref{fg_18} is an illustration of the convergence behavior and the network
training time of the locHLConcELM method with respect to the size
of the neural network.
The local network architecture is given by $[2, M, 30, 1]$, with
$M$ varied systematically.
We employ a fixed uniform set of $Q=21\times 21$ collocation points,
and a hidden magnitude vector $\mbs R=(0.9,0.05)$ obtained using
the method of~\cite{DongY2021}.
One observes that the errors decrease exponentially
and that the network training time grows superlinearly with increasing $M$.

\begin{table}[tb]
  \centering
  \begin{tabular}{l|l|cc|cc}
    \hline
    local network & collocation & current & locHLConcELM & conventional & locELM  \\ \cline{3-6}
    architecture & points & max error & rms error & max error & rms error \\ \hline
    $[2,200,30,1]$ & $5\times 5$ & $1.39E-2$ & $2.64E-3$ & $4.12E-2$ & $1.21E-2$ \\
    & $10\times 10$ & $2.30E-4$ & $4.44E-5$ & $6.62E-2$ & $2.84E-2$ \\
    & $15\times 15$ & $2.98E-6$ & $5.83E-7$ & $7.42E-2$ & $3.16E-2$ \\
    & $20\times 20$ & $4.01E-7$ & $7.06E-8$ & $7.97E-2$ & $3.39E-2$ \\
    & $25\times 25$ & $5.59E-7$ & $9.04E-8$ & $8.44E-2$ & $3.58E-2$ \\
    & $30\times 30$ & $6.60E-7$ & $1.11E-7$ & $8.86E-2$ & $3.76E-2$ \\
    \hline
  \end{tabular}
  \caption{Burgers' equation on the smaller domain $\Omega_1$:
    Comparison of the maximum/rms errors
    from the locHLConcELM method and the
    conventional locELM method~\cite{DongL2021}.
    The locHLConcELM data in this table correspond to those in
    Figure~\ref{fg_17}(a).
    For conventional locELM, the random hidden-layer coefficients
    are set to uniform random values generated on $[-R_m,R_m]$
    with $R_m=R_{m0}=0.175$, where $R_{m0}$ is the optimal $R_m$ obtained
    using the method of~\cite{DongY2021}.
  }
  \label{tab_4}
\end{table}

Table~\ref{tab_4} provides an accuracy comparison of the HLConcELM
method and the conventional locELM method~\cite{DongL2021} for solving
the Burgers' equation on the smaller domain $\Omega_1$.
With both methods, we employ $4$ sub-domains as shown in Figure~\ref{fg_16}(a),
a local neural network architecture $[2, 200, 30, 1]$, and
a sequence of uniform collocation points ranging from $Q=5\times 5$
and $Q=30\times 30$.
The current locHLConcELM method is significantly more accurate than
the conventional locELM method, with their maximum errors on the order of $10^{-7}$
and $10^{-2}$ respectively.

\begin{figure}
  \centerline{
    \includegraphics[height=2.0in]{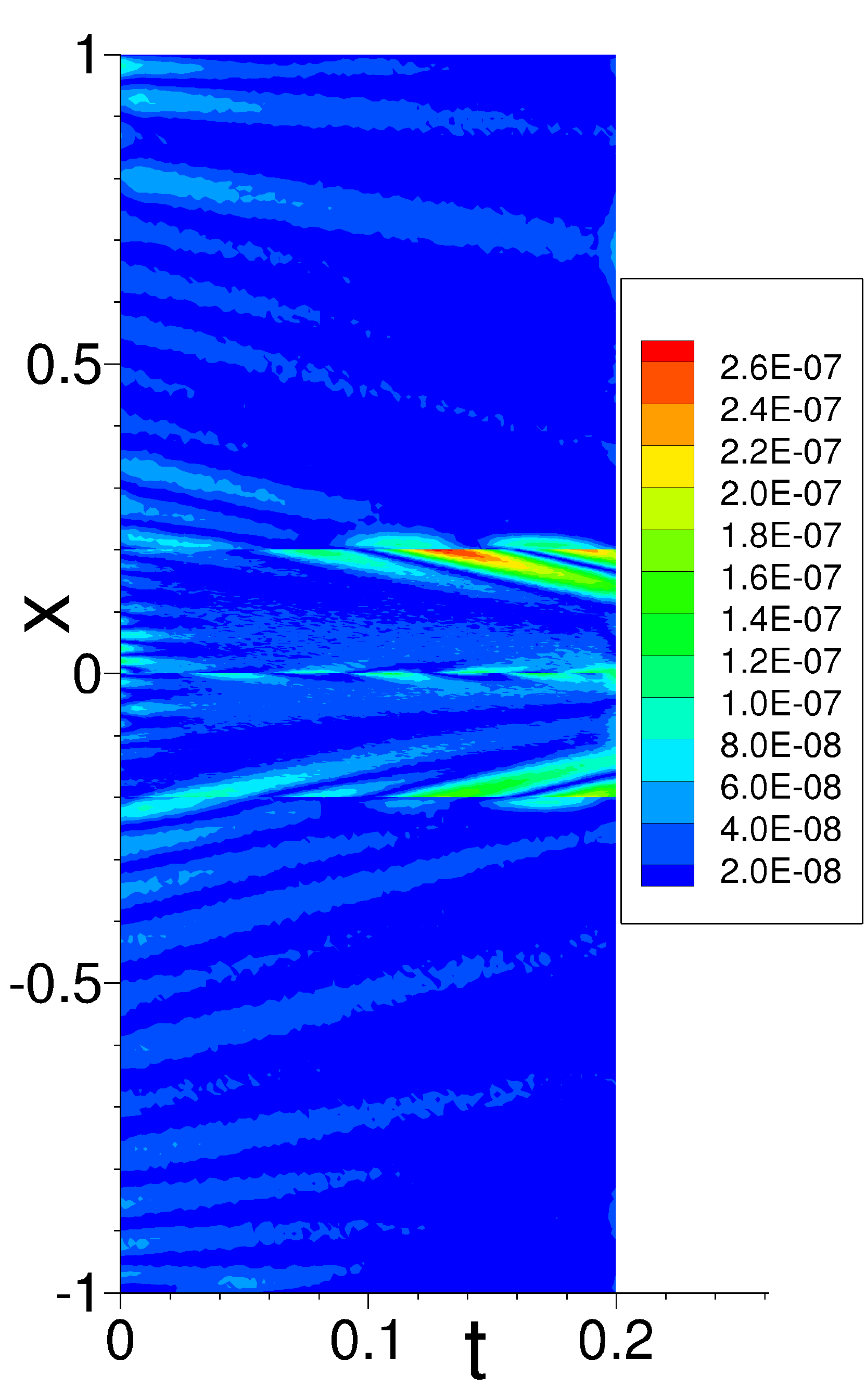}(a)
    \includegraphics[height=2.0in]{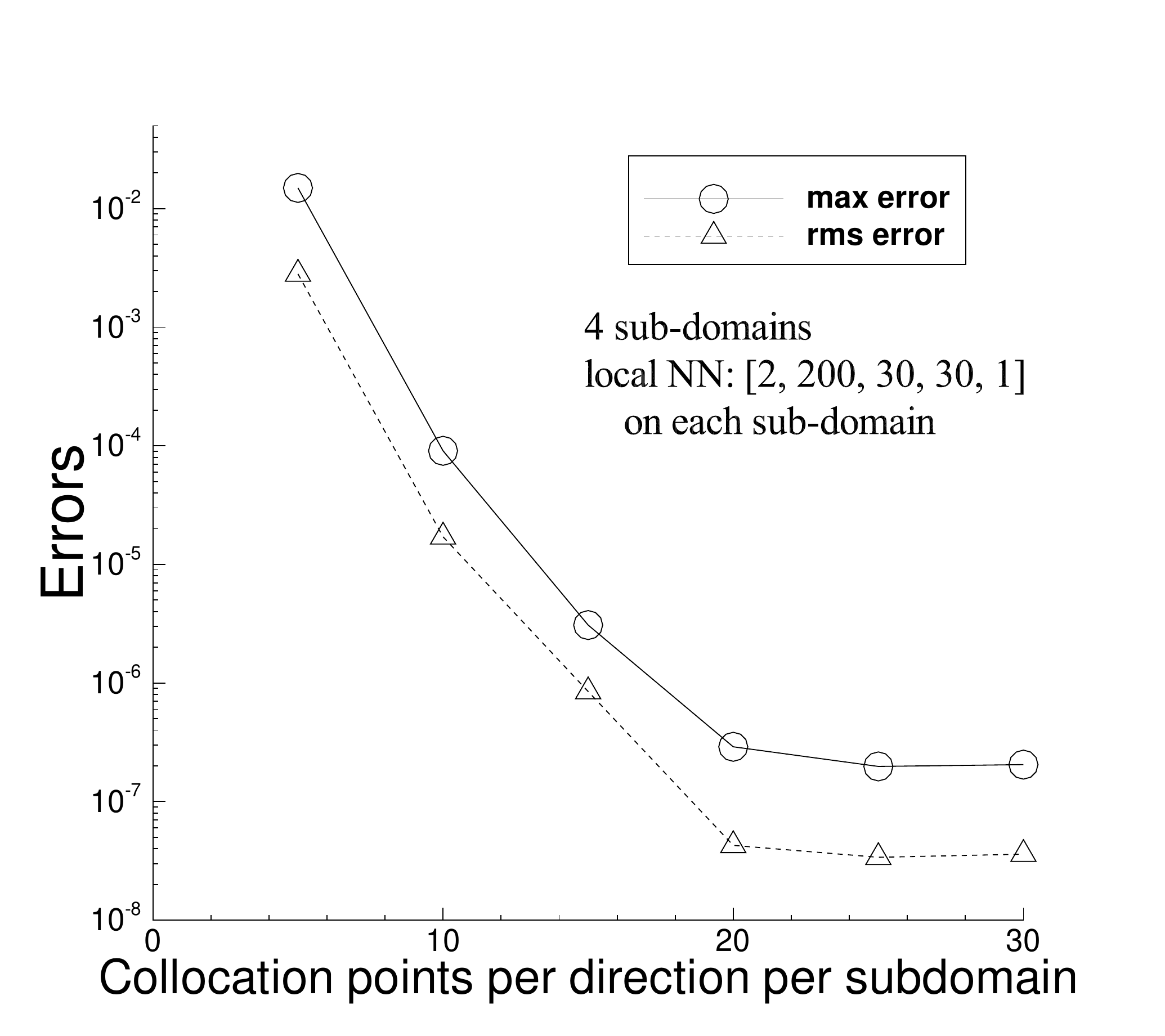}(b)
    \includegraphics[height=2.0in]{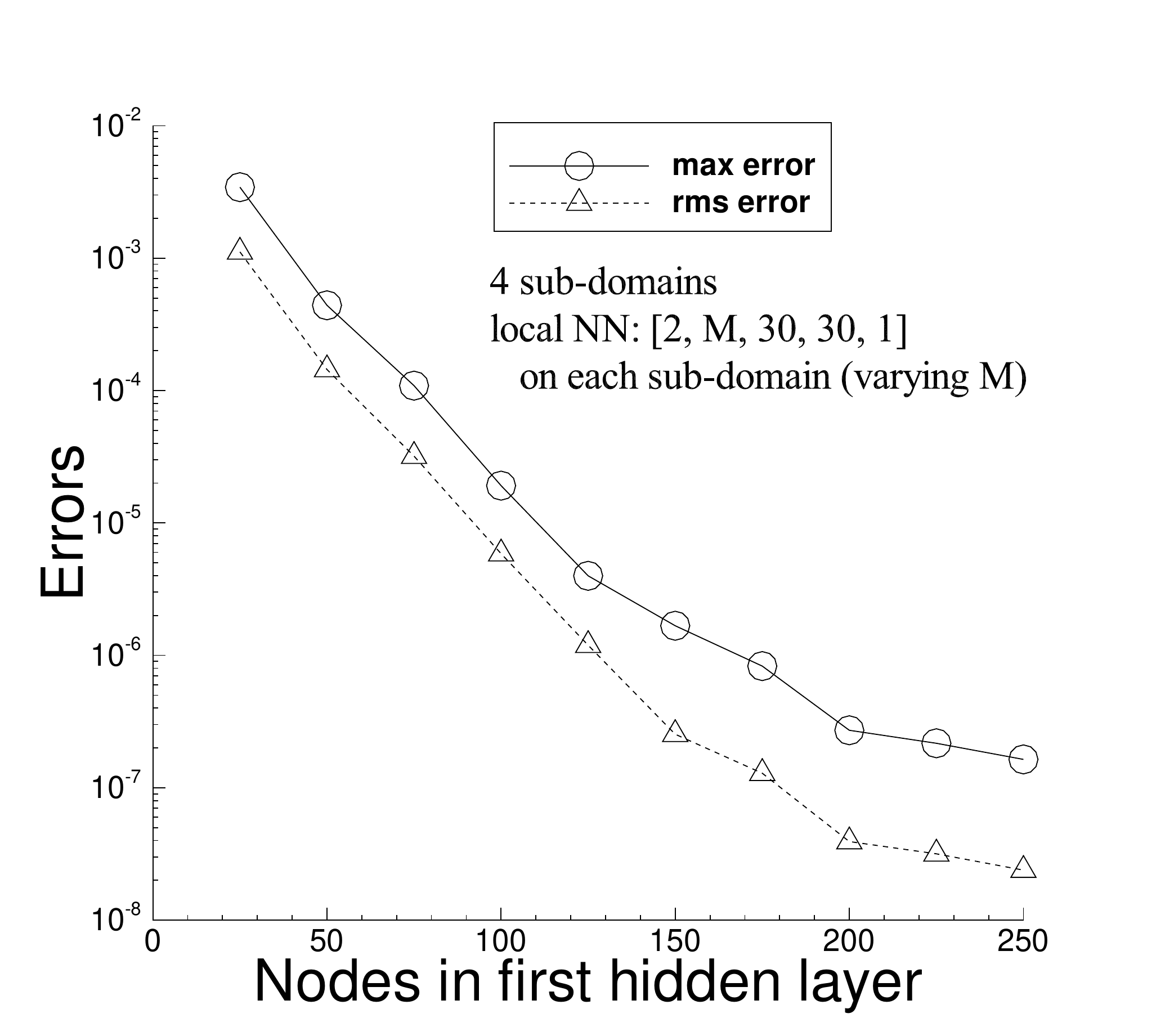}(c)
  }
  \caption{Burgers' equation on the smaller domain $\Omega_1$ ($3$ hidden layers in network):
    (a) Error distribution of the locHLConcELM solution.
    The locHLConcELM maximum/rms errors versus (b) the number of
    collocation points per direction in each sub-domain,
    and (c) the number of nodes in
    the first hidden layer ($M$).
    Local network architecture: $[2, M, 30, 30, 1]$.
    $M=200$ in (a,b), and is varied in (c).
    $Q=21\times 21$ in (a,c), and is varied in (b).
    $\mbs R=(1.0, 0.035, 0.03)$ in (a,b,c).
  }
  \label{fg_19}
\end{figure}

Figure~\ref{fg_19} illustrates the characteristics of the locHLConcELM solution
obtained with $3$ hidden layers in the local neural network on
the smaller domain $\Omega_1$.
Here we employ a local network architecture $[2, M, 30, 30, 1]$, with $M$
either fixed at $M=200$ or varied systematically, and
a uniform set of $Q=Q_1\times Q_1$ collocation points, with
$Q_1$ either fixed at $Q_1=21$ or varied systematically.
The hidden magnitude vector is $\mbs R=(1.0, 0.035, 0.03)$,
obtained using the method of~\cite{DongY2021}.
Figure~\ref{fg_19}(a) is an illustration of the absolute-error
distribution on $\Omega_1$ corresponding to $M=200$ and $Q_1=21$,
demonstrating a high accuracy with the maximum error on the order of $10^{-7}$.
Figures~\ref{fg_19}(b) and (c) show the exponential convergence behavior
of the HLConcELM errors with respect to $Q_1$ and $M$.

\begin{figure}
  \centerline{
    \includegraphics[width=2.in]{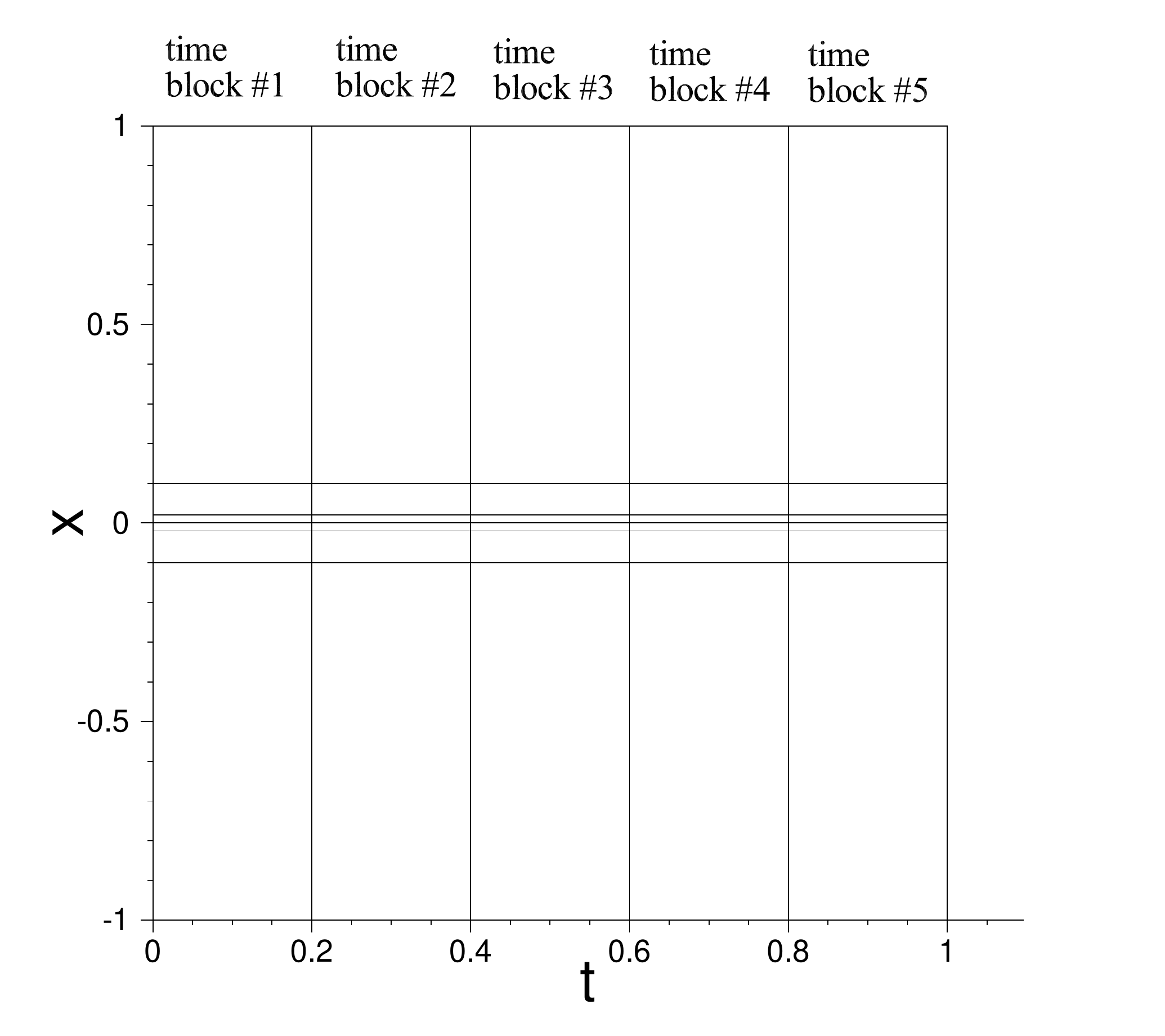}(a)
    \includegraphics[width=2.in]{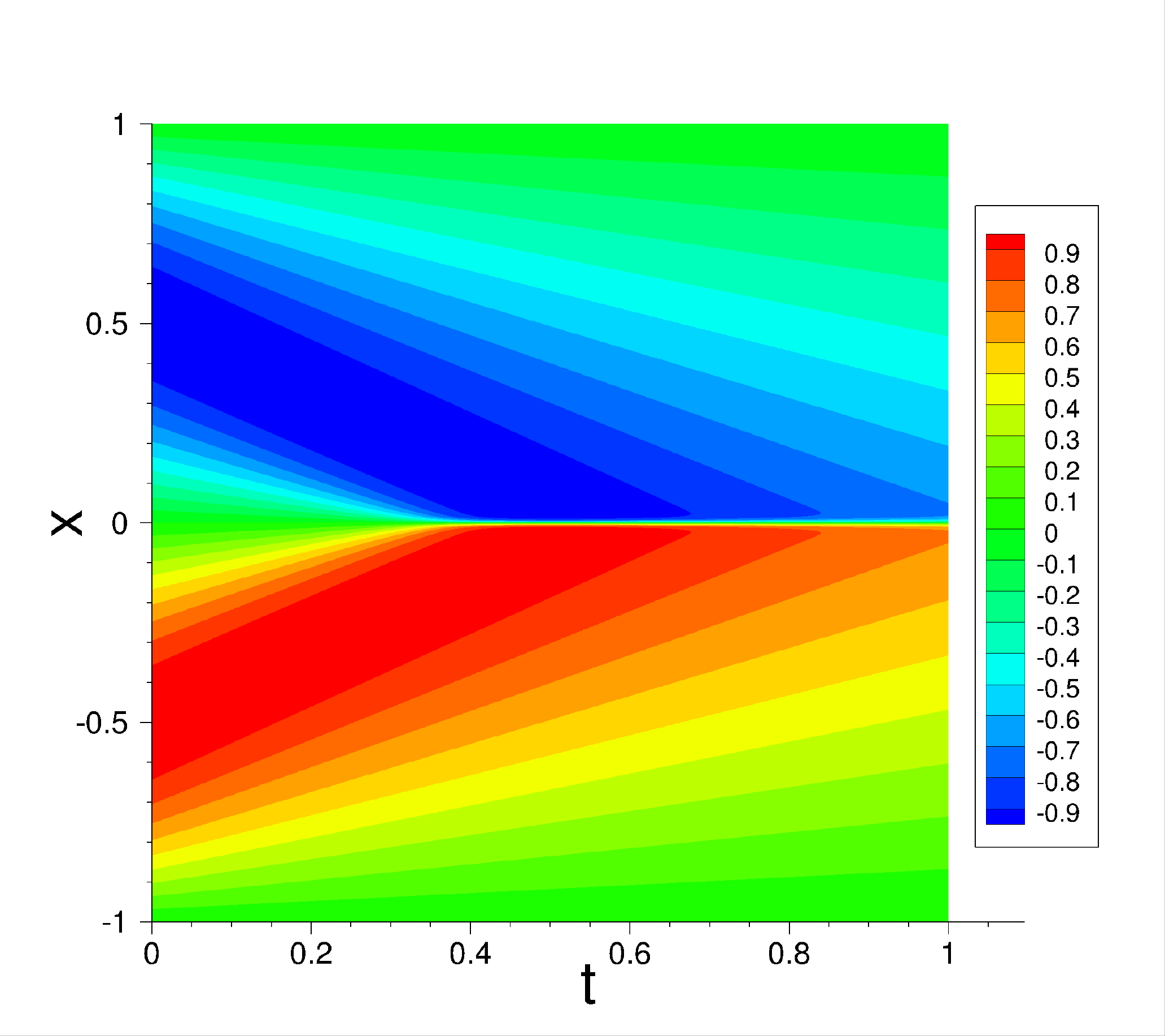}(b)
    \includegraphics[width=2.in]{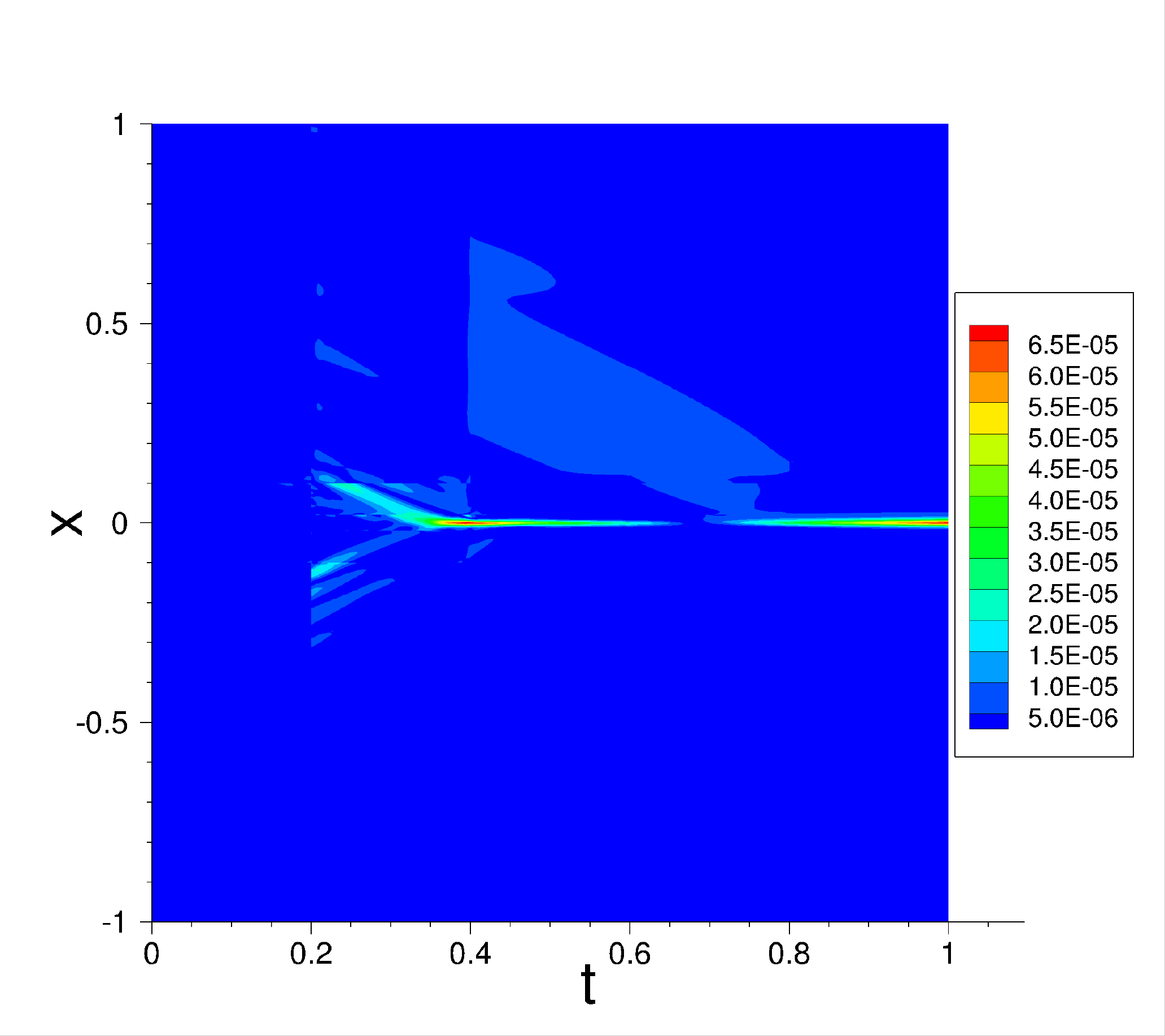}(c)
  }
  \caption{Burgers' equation on the larger domain $\Omega$ ($t\in[0,1]$):
    (a) Configuration of the $5$ uniform time blocks and the $6$ non-uniform sub-domains
    on each time block.
    Distributions of (b) the locHLConcELM solution and (c) its absolute error.
    Local network architecture $[2, 300, 1]$ on each sub-domain,
    $Q=21\times 21$ per sub-domain,
    $\mbs R=2.0$.
  }
  \label{fg_20}
\end{figure}

\begin{figure}
  \centerline{
    \includegraphics[width=1.9in]{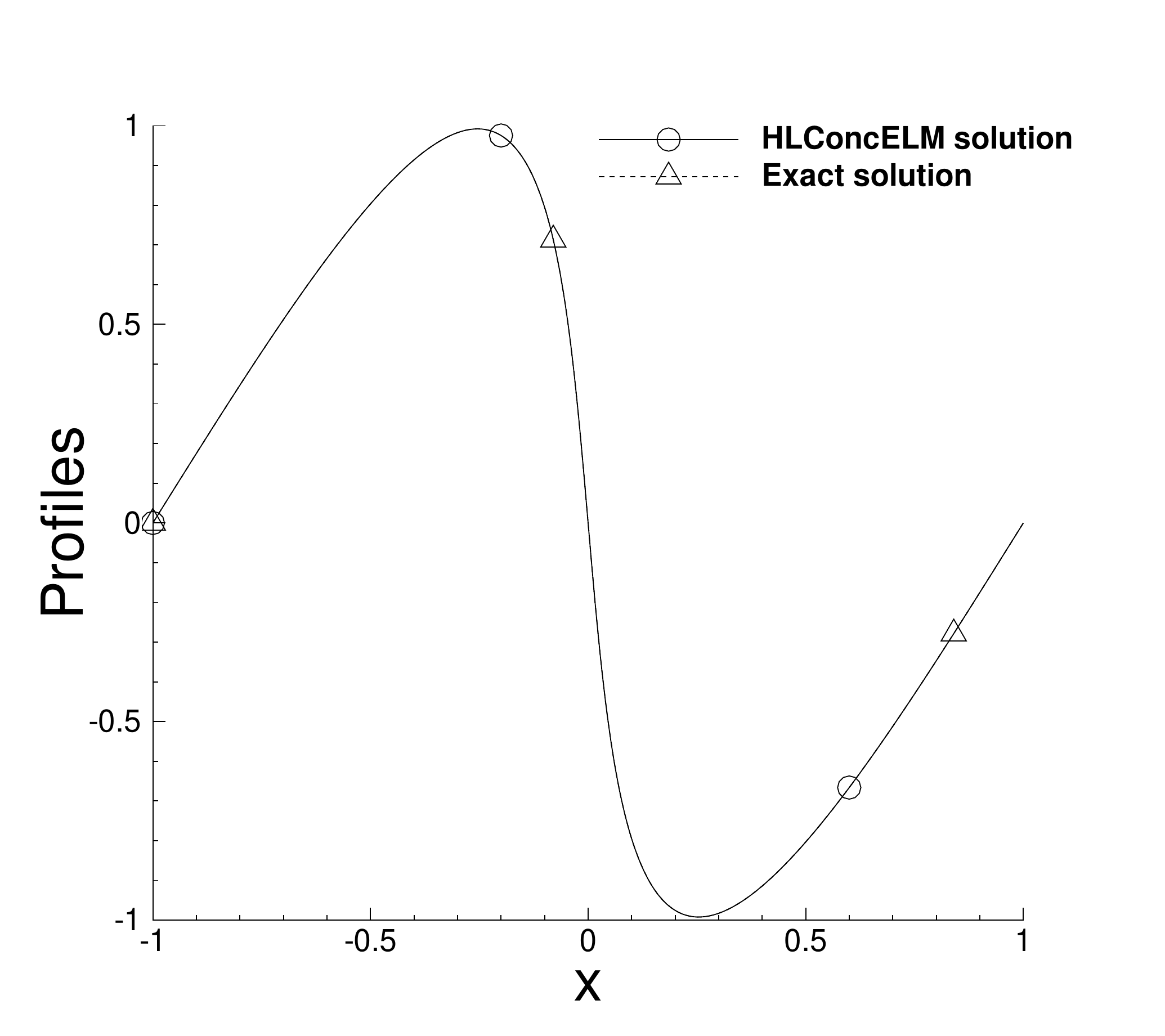}(a)
    \includegraphics[width=1.9in]{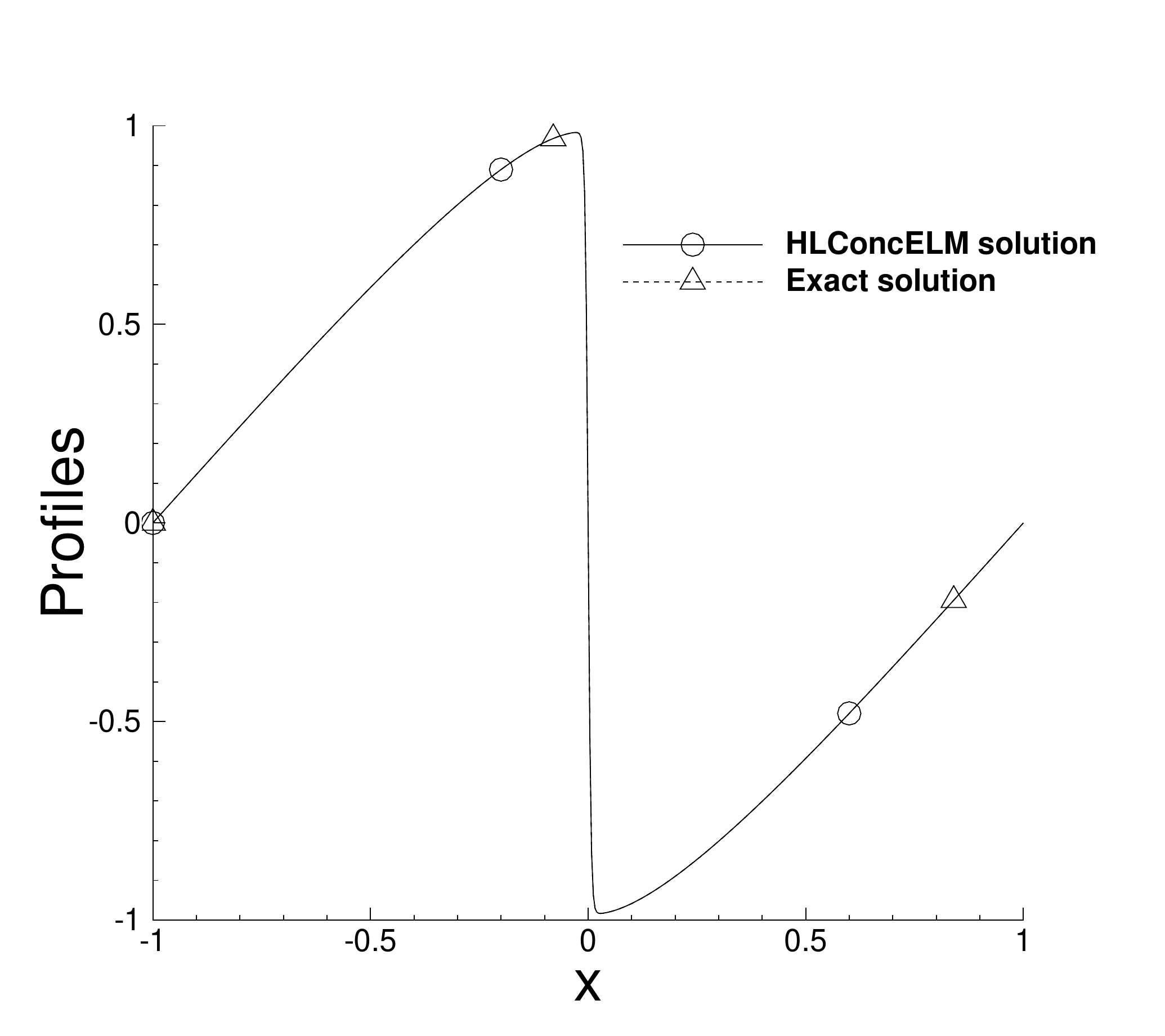}(b)
    \includegraphics[width=1.9in]{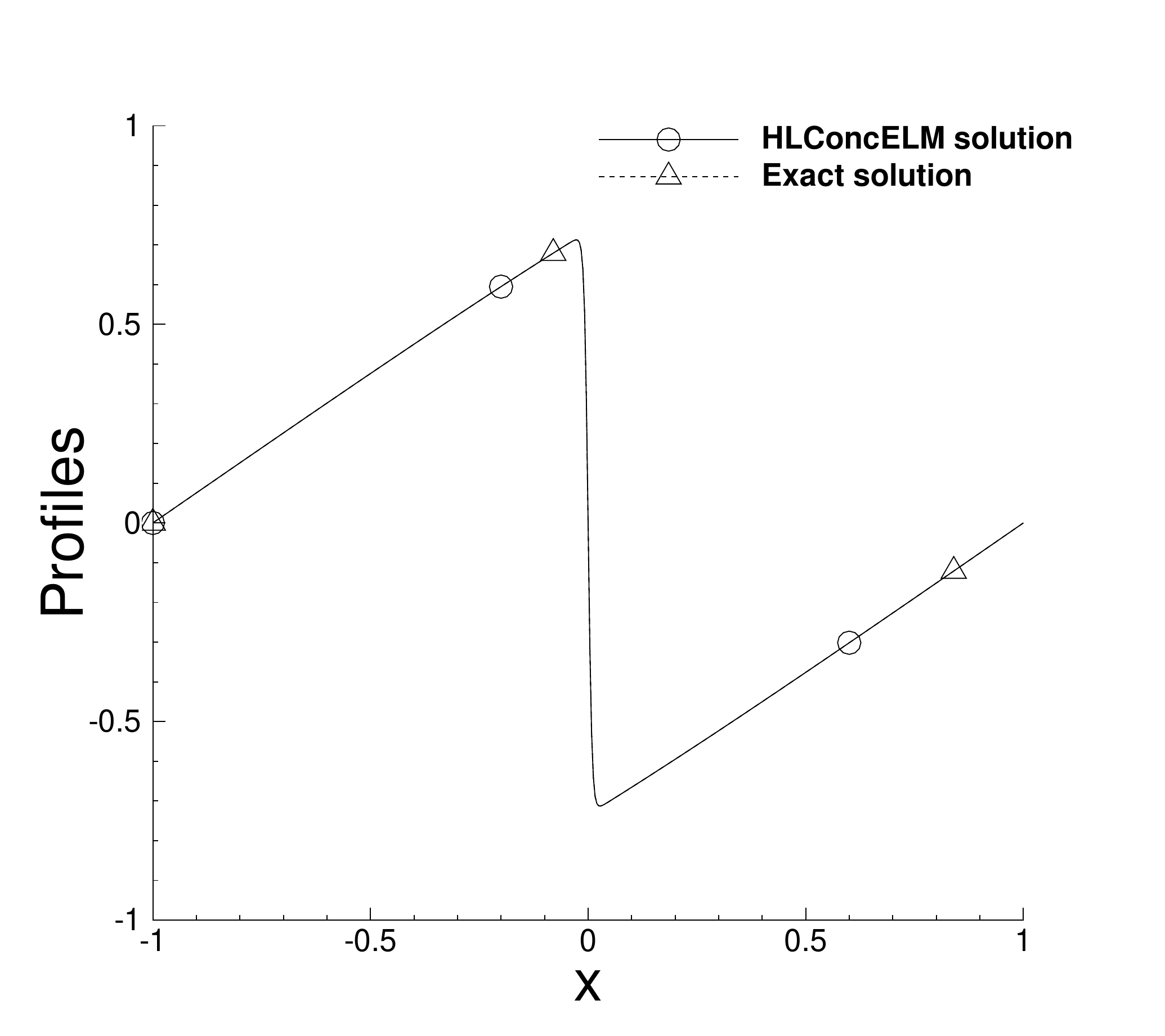}(c)
  }
  \centerline{
    \includegraphics[width=1.9in]{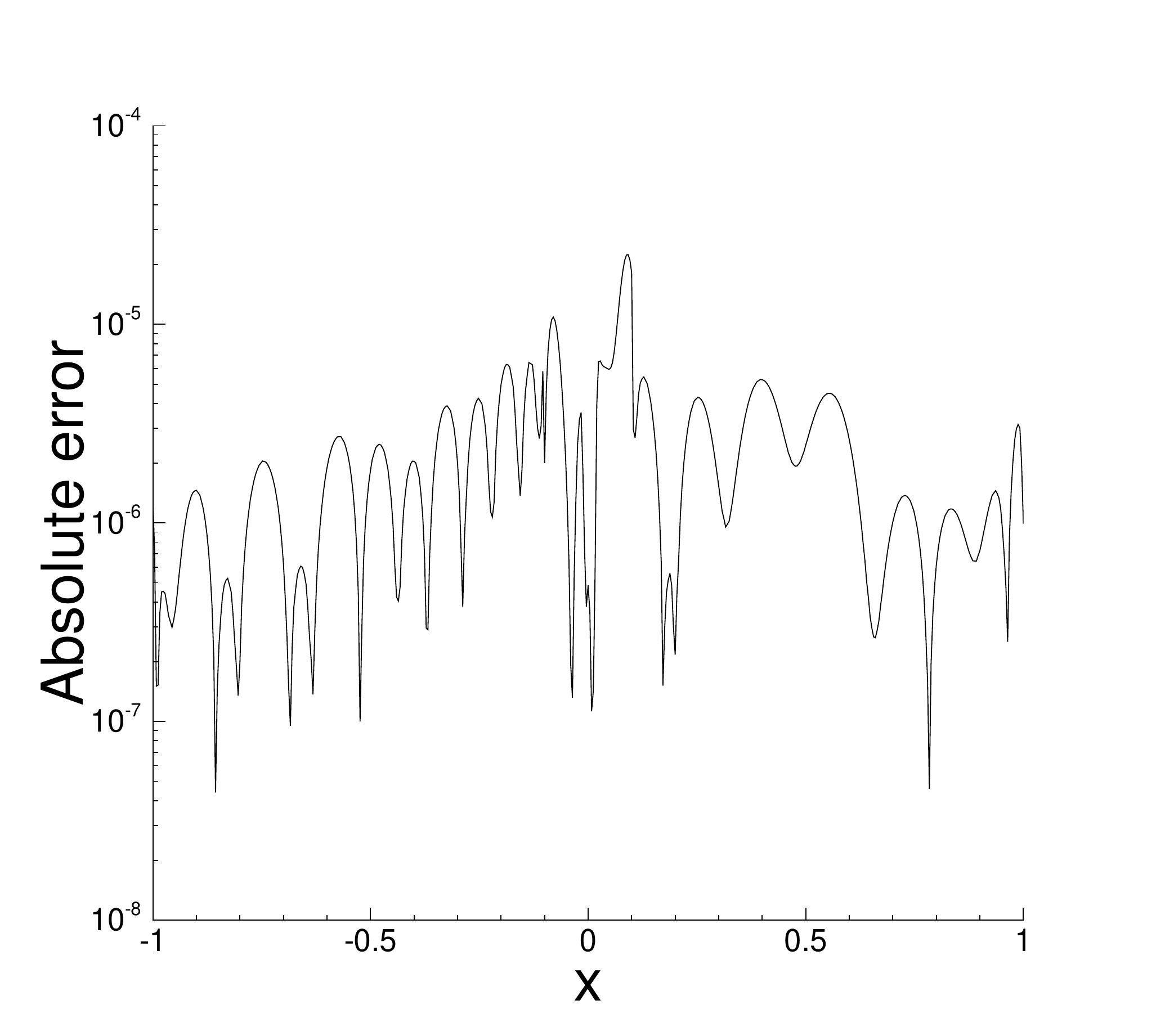}(d)
    \includegraphics[width=1.9in]{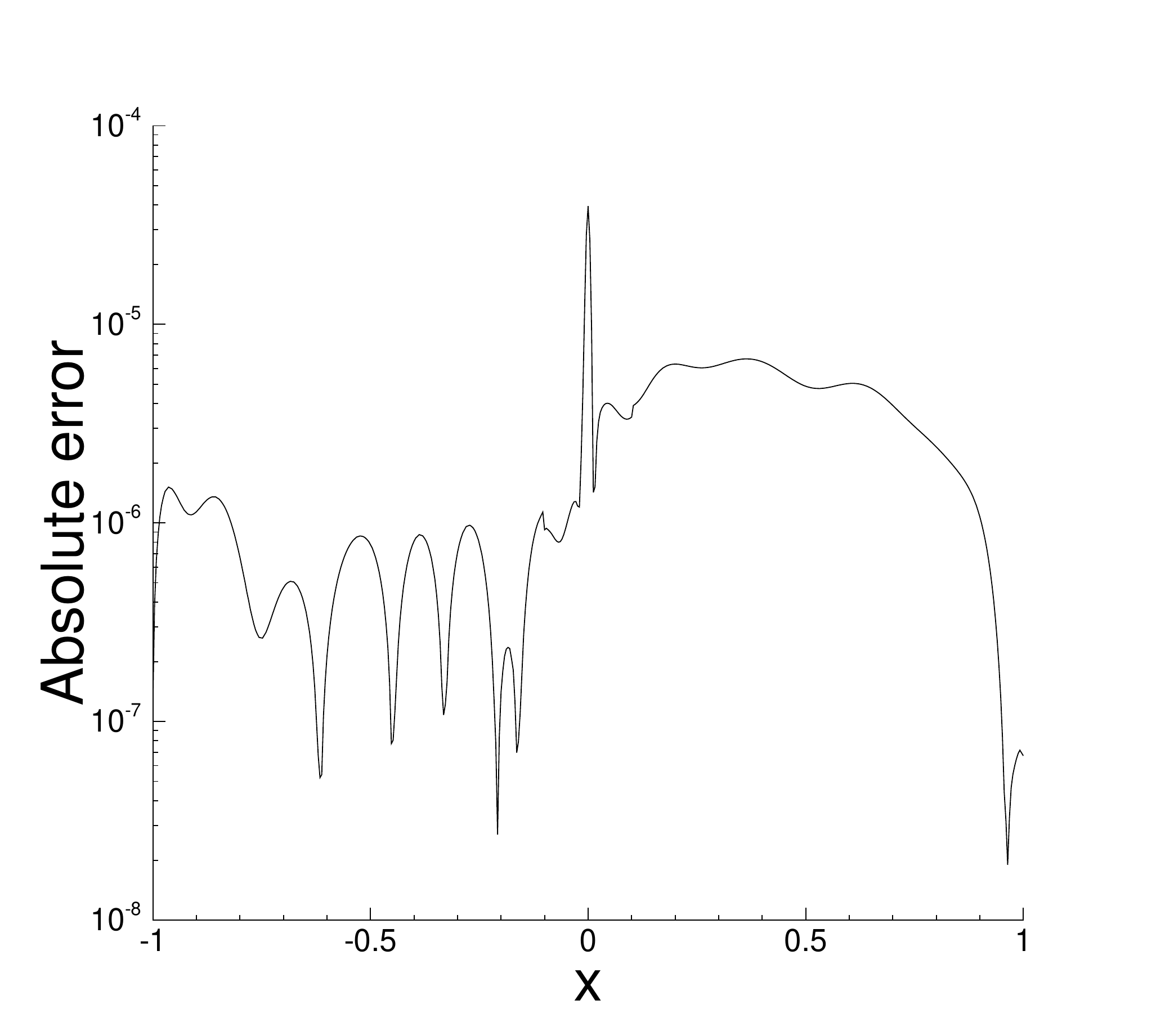}(e)
    \includegraphics[width=1.9in]{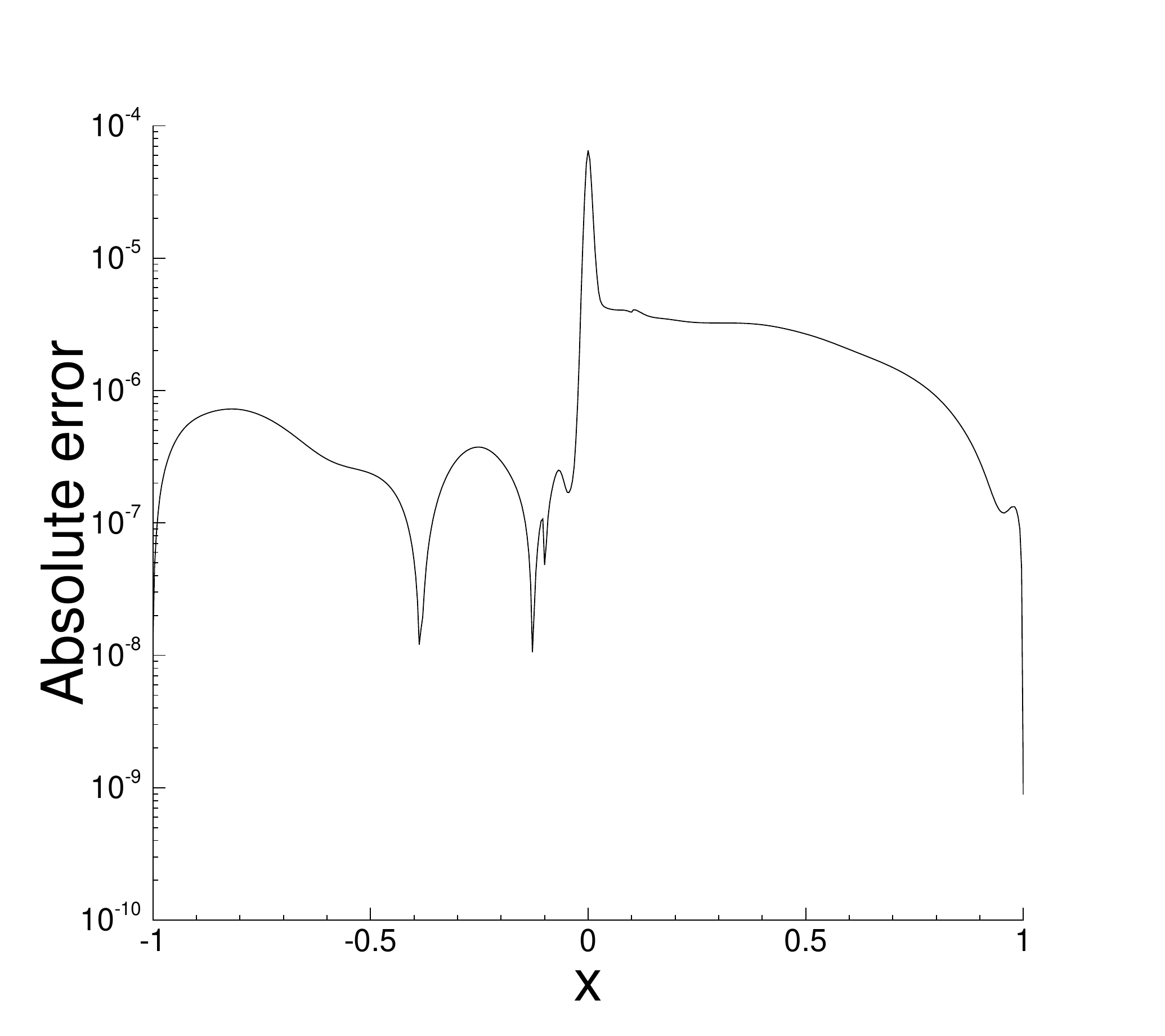}(f)
  }
  \caption{Burgers' equation on the larger domain $\Omega$:
    Profiles of the locHLConcELM solution (top row) and its absolute error (bottom row)
    at the time instants: $t=0.25$ (a,d), $t=0.5$ (b,e), and $t=1.0$ (c,f).
    The profiles of the exact solution are also shown in (a,b,c),
    which overlap with those of the locHLConcELM solution.
    The locHLConcELM simulation parameters and configurations follow those
    of Figure~\ref{fg_20}.
  }
  \label{fg_21}
\end{figure}

Let us next consider the larger domain $\Omega$ ($t\in[0,1]$) and solve
the system~\eqref{eq_27} using the current method.
We employ the locHLConcELM method together with the block time marching scheme
(see Remarks~\ref{rem_3} and~\ref{rem_4}) in the simulation.
Specifically, we divide the temporal dimension into $5$ uniform time blocks,
and partition each time block into $6$ non-uniform sub-domains along the
$x$ direction. Figure~\ref{fg_20}(a) illustrates the configuration of
the time blocks and the sub-domains on each time block,
where the $x$ coordinates of the sub-domain boundaries are
given by $\bm{\mathcal{X}}=[-1, -0.1, -0.02, 0, 0.02, 0.1, 1]$.
We employ a local neural network architecture $\mbs M=[2, 300, 1]$ and
a uniform set of $Q=21\times 21$ collocation points on each sub-domain.
The hidden magnitude vector is $\mbs R=2.0$, which is obtained using
the method of~\cite{DongY2021}.
Figures~\ref{fg_20}(b) and (c) show the distributions of the locHLConcELM
solution and its absolute error on $\Omega$. The data indicate
that the current method achieves a quite high accuracy with
the sharp gradient present in the domain, with the maximum error on
the order of $10^{-5}$.

Figure~\ref{fg_21} compares profiles of the locHLConcELM solution
and the exact solution~\eqref{eq_28} for the Burgers' equation
at three time instants $t=0.25$, $0.5$ and $1.0$.
The error profiles of the locHLConcELM solution have also been included
in this figure. The simulation configuration and the parameters used here
correspond to those of Figure~\ref{fg_20}.
It is evident that the current locHLConcELM method has achieved a quite
high accuracy for this problem.


\subsection{KdV Equation}
\label{sec:kdv}

\begin{figure}
  \centerline{
    \includegraphics[width=2.in]{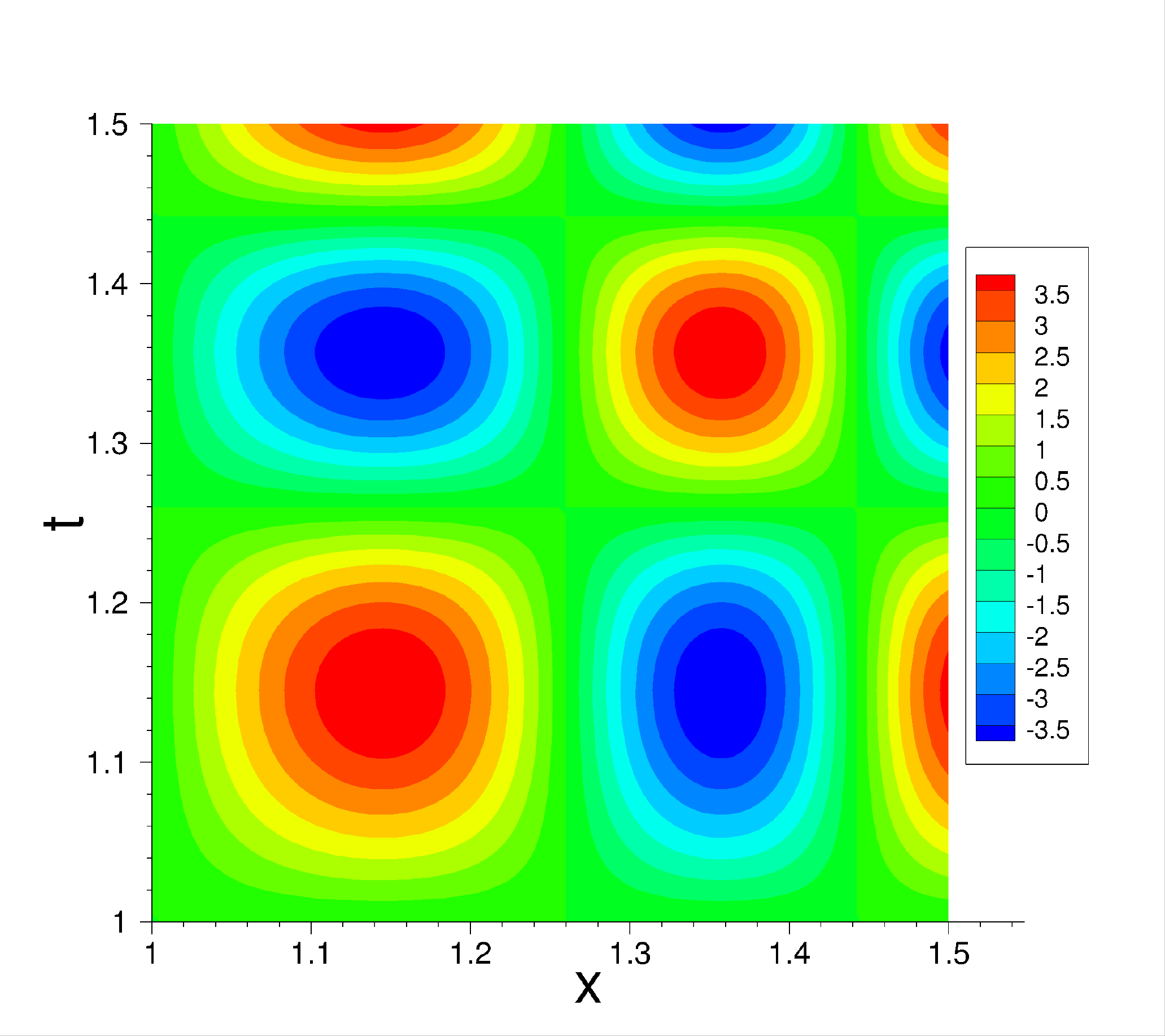}(a)
    \includegraphics[width=2.in]{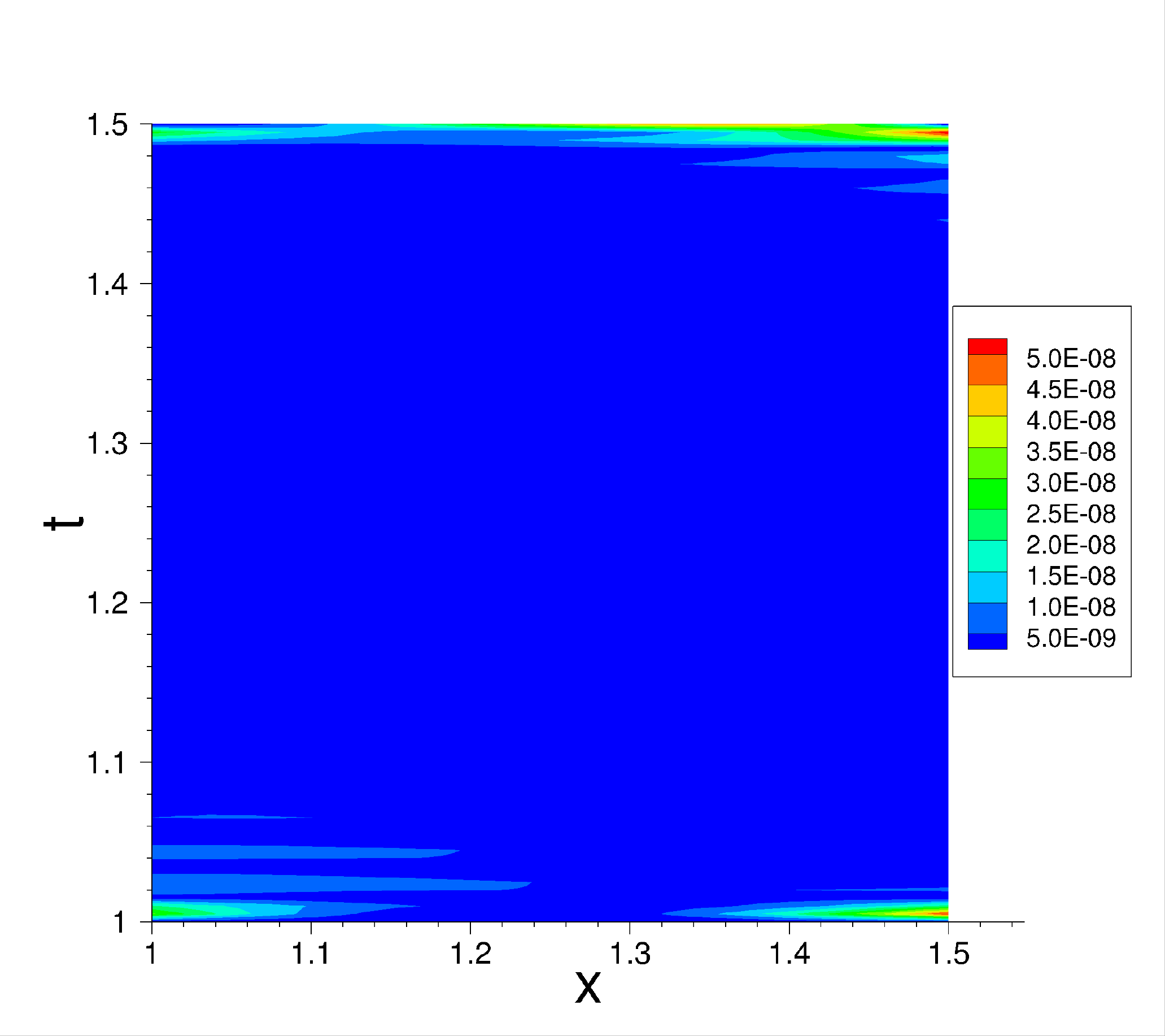}(b)
    \includegraphics[width=2.in]{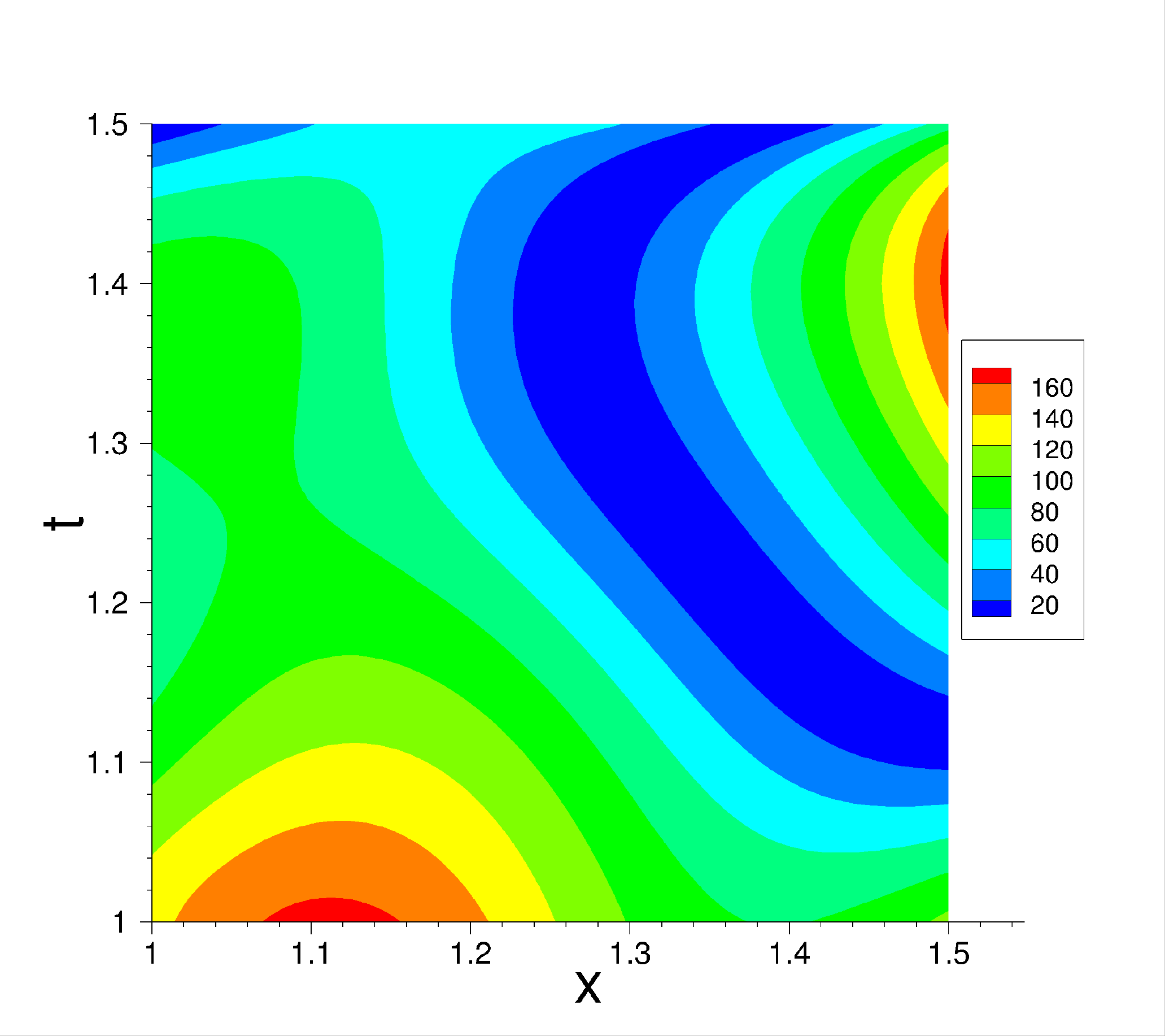}(c)
  }
  \caption{KdV equation: Distributions of (a) the exact solution, 
    (b) the absolute error of the HLConcELM solution, and (c)
    the absolute error of the conventional ELM solution.
    In (b,c), neural network architecture $\mbs M=[2,800,50,1]$,
    $Q=35\times 35$ uniform collocation points.
    $\mbs R=(3.2,0.01)$ in (b) for HLConcELM.
    $R_m=R_{m0}=0.27$ in (c) for conventional ELM.
  }
  \label{fg_22}
\end{figure}

\begin{figure}
  \centerline{
    \includegraphics[width=2.in]{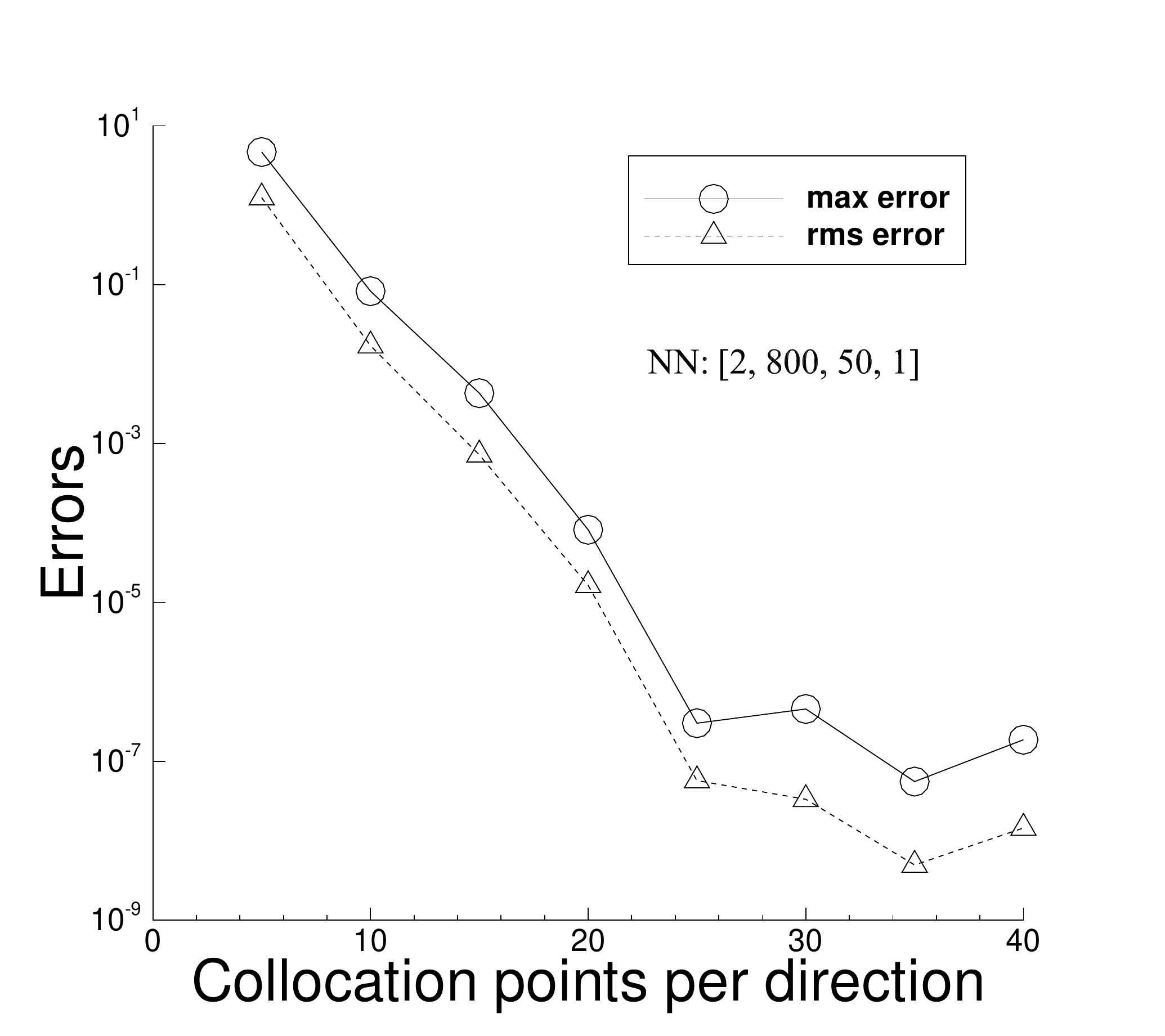}(a)
    \includegraphics[width=2.in]{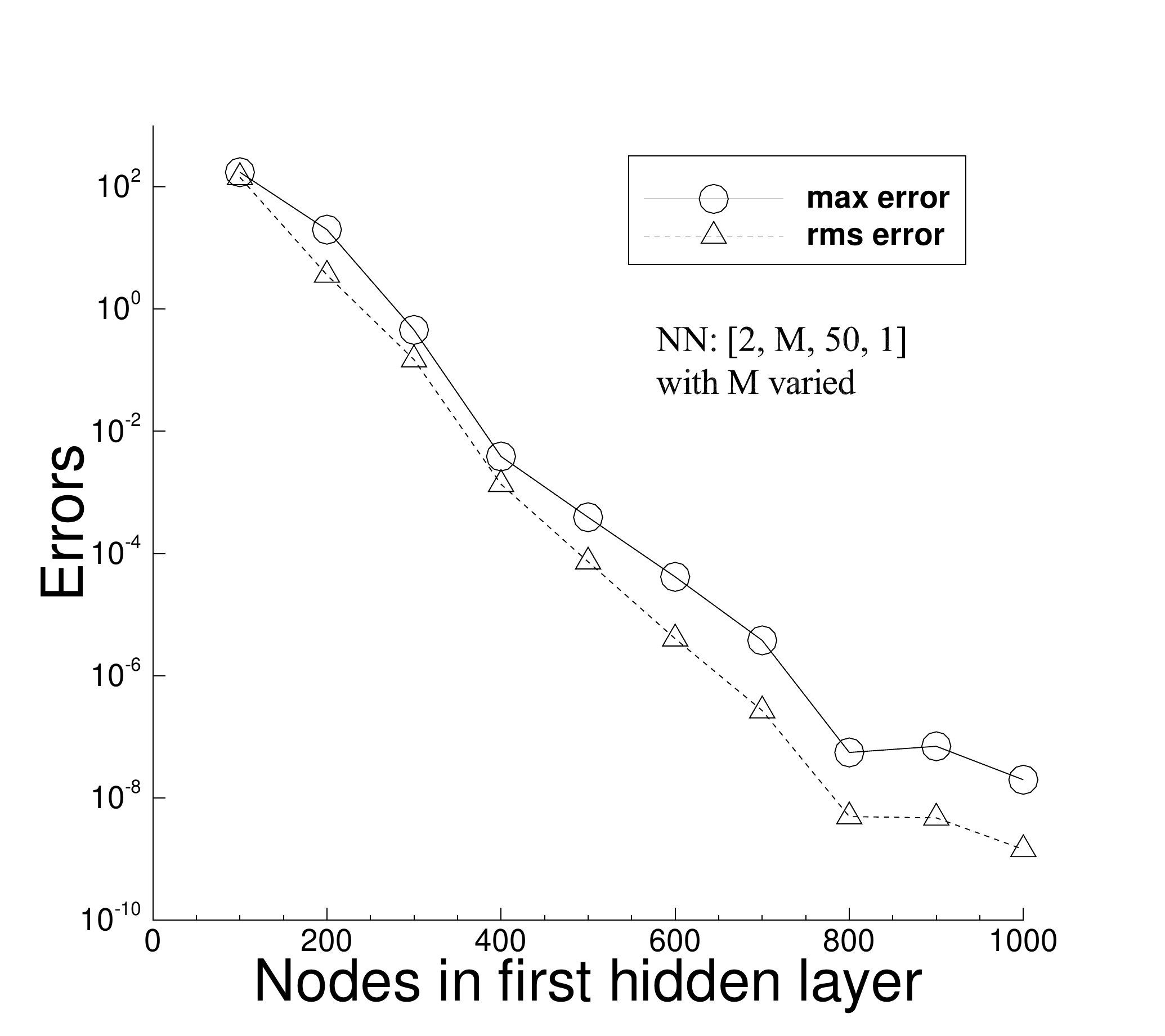}(b)
  }
  \caption{KdV equation: The HLConcELM maximum/rms errors 
    versus (a) the number of collocation points per direction,
    and (b) the number of nodes in the first hidden layer ($M$).
    Network architecture $[2, M, 50, 1]$.
    $M=800$ in (a), varied in (b). $Q=35\times 35$ in (b), varied in (a).
    $\mbs R=(3.2,0.01)$ in (a,b).
  }
  \label{fg_23}
\end{figure}

In the next benchmark problem we employ the Korteweg-de Vries (KdV)
equation to test the HLConcELM method.
Consider the spatial-temporal domain $(x,t)\in\Omega=[1.0,1.5]\times[1.0,1.5]$
and the following initial/boundary value problem,
\begin{subequations}\label{eq_29}
  \begin{align}
    &
    \frac{\partial u}{\partial t} - u\frac{\partial u}{\partial x}
    + \frac{\partial^3u}{\partial x^3} = f(x,t), \\
    &
    u(1,t) = g_1(t), \quad
    u(1.5,t) = g_2(t), \quad
    \frac{\partial u}{\partial x}(1,t) = g_3(t), \\
    &
    u(x,1) = h(x).
  \end{align}
\end{subequations}
In the above equations $u(x,t)$ is the field solution to be sought,
$f(x,t)$ is a prescribed source term, $g_i$ ($i=1,2,3$) and $h$
are the data for the boundary and initial conditions.
We choose $f$, $g_i$ ($i=1,2,3$) and $h$ such that
the system~\eqref{eq_29} has the following analytic solution,
\begin{equation}\label{eq_30}
  u(x,t) = 4\sin\left(\pi x^3\right) \sin\left(\pi t^3\right).
\end{equation}
Figure~\ref{fg_23}(a) shows the distribution of this
solution in the spatial-temporal $xt$ plane.


To solve the problem~\eqref{eq_29} using the HLConcELM method,
we employ neural networks with two input nodes (representing
the $x$ and $t$) and a single output node (representing $u$),
with the Gaussian activation function for all the hidden nodes.
A uniform set of $Q=Q_1\times Q_1$ collocation points on the domain
$\Omega$ is used to train the neural network, where
$Q_1$ is varied systematically in the tests.

Figures~\ref{fg_22}(b) and (c) illustrate the absolute-error
distributions obtained using the HLConcELM method
and the conventional ELM method.
Here the network architecture is given by $\mbs M=[2, 800, 50, 1]$,
and a uniform set of $Q=35\times 35$ collocation points is used
for both methods. The hidden magnitude vector is
$\mbs R=(3.2, 0.01)$ with HLConcELM, which is
obtained using the method of~\cite{DongY2021}.
For conventional ELM, we set the hidden-layer coefficients
to random values generated on $[-R_m,R_m]$ with $R_m=R_{m0}$,
where $R_{m0}=0.27$ is the optimal $R_m$ obtained using
the method of~\cite{DongY2021}.
The conventional ELM solution is observed to be inaccurate
(maximum error on the order of $10^2$), because of the
narrow last hidden layer in the network.
In contrast, the HLConcELM solution
is highly accurate, with the maximum error
on the order of $10^{-8}$ in the domain.

Figure~\ref{fg_23} illustrates the convergence
behavior of the HLConcELM errors
with respect to the collocation points and
the number of nodes in the network.
Here the network architecture is $\mbs M=[2, M, 50, 1]$,
with $M$ either fixed at $M=800$ or varied between $M=100$
and $M=1000$. A uniform set of $Q=Q_1\times Q_1$ collocation
points is used, with $Q_1$ either fixed at $Q_1=35$ or varied
between $Q_1=5$ and $Q_1=40$.
The hidden magnitude vector is $\mbs R=(3.2,0.01)$,
obtained from the method of~\cite{DongY2021}.
The two plots (a) and (b) depict the maximum/rms errors
of the HLConcELM solution as a function of $Q_1$ and $M$,
respectively.
One can observe the familiar exponential decrease
in the errors with increasing $Q_1$ or $M$.

\begin{table}[tb]
  \centering
  \begin{tabular}{l|l|ll|ll}
    \hline
    neural & collocation & \quad current & ELM & conventional & ELM  \\ \cline{3-6}
    network & points & max error & rms error & max error & rms error \\ \hline
    $[2,800,50,1]$ & $5\times 5$ & $4.67E+0$ & $1.24E+0$ & $1.66E+2$ & $8.18E+1$  \\
    & $10\times 10$ & $8.28E-2$ & $1.70E-2$ & $1.78E+2$ & $1.02E+2$  \\
    & $15\times 15$ & $4.30E-3$ & $7.21E-4$ & $1.88E+2$ & $9.36E+1$  \\
    & $20\times 20$ & $8.21E-5$ & $1.64E-5$ & $1.78E+2$ & $9.80E+1$ \\
    & $25\times 25$ & $3.02E-7$ & $5.71E-8$ & $1.75E+2$ & $8.89E+1$ \\
    & $30\times 30$ & $4.55E-7$ & $3.32E-8$ & $1.66E+2$ & $7.76E+1$ \\
    & $35\times 35$ & $5.55E-8$ & $4.95E-9$ & $1.67E+2$ & $7.77E+1$ \\
    & $40\times 40$ & $1.87E-7$ & $1.45E-8$ & $1.67E+2$ & $7.74E+1$ \\
    \hline
  \end{tabular}
  \caption{KdV equation:
    Comparison of the maximum/rms errors
    from the HLConcELM method and the
    conventional ELM method.
    Network architecture: $[2, 800, 50, 1]$.
    The HLConcELM data in this table correspond to those in
    Figure~\ref{fg_23}(a).
    For conventional ELM, the hidden-layer coefficients
    are set to uniform random values generated on $[-R_m,R_m]$
    with $R_m=R_{m0}$. Here $R_{m0}=0.27$ is the optimal $R_m$ obtained
    using the method of~\cite{DongY2021}.
  }
  \label{tab_5}
\end{table}

Table~\ref{tab_5} provides an accuracy comparison of
the HLConcELM method and the conventional ELM
method~\cite{DongL2021} for solving the KdV equation
on a network architecture $\mbs M=[2, 800, 50, 1]$
corresponding to a sequence of collocation points.
The HLConcELM solution is
highly accurate, while the conventional ELM solution
exhibits no accuracy at all on such a neural network.

\begin{figure}
  \centerline{
    \includegraphics[width=2in]{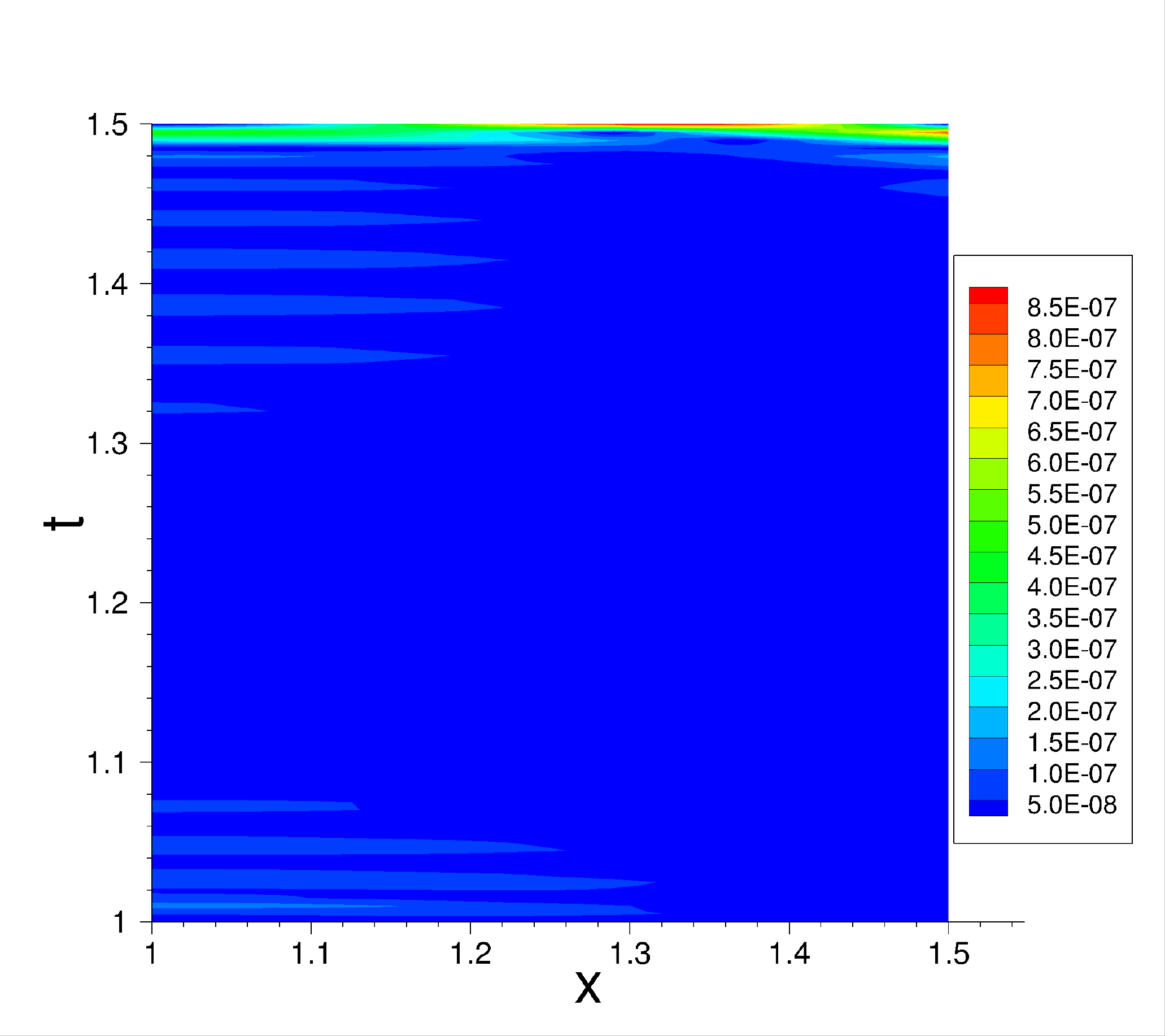}(a)
    \includegraphics[width=2in]{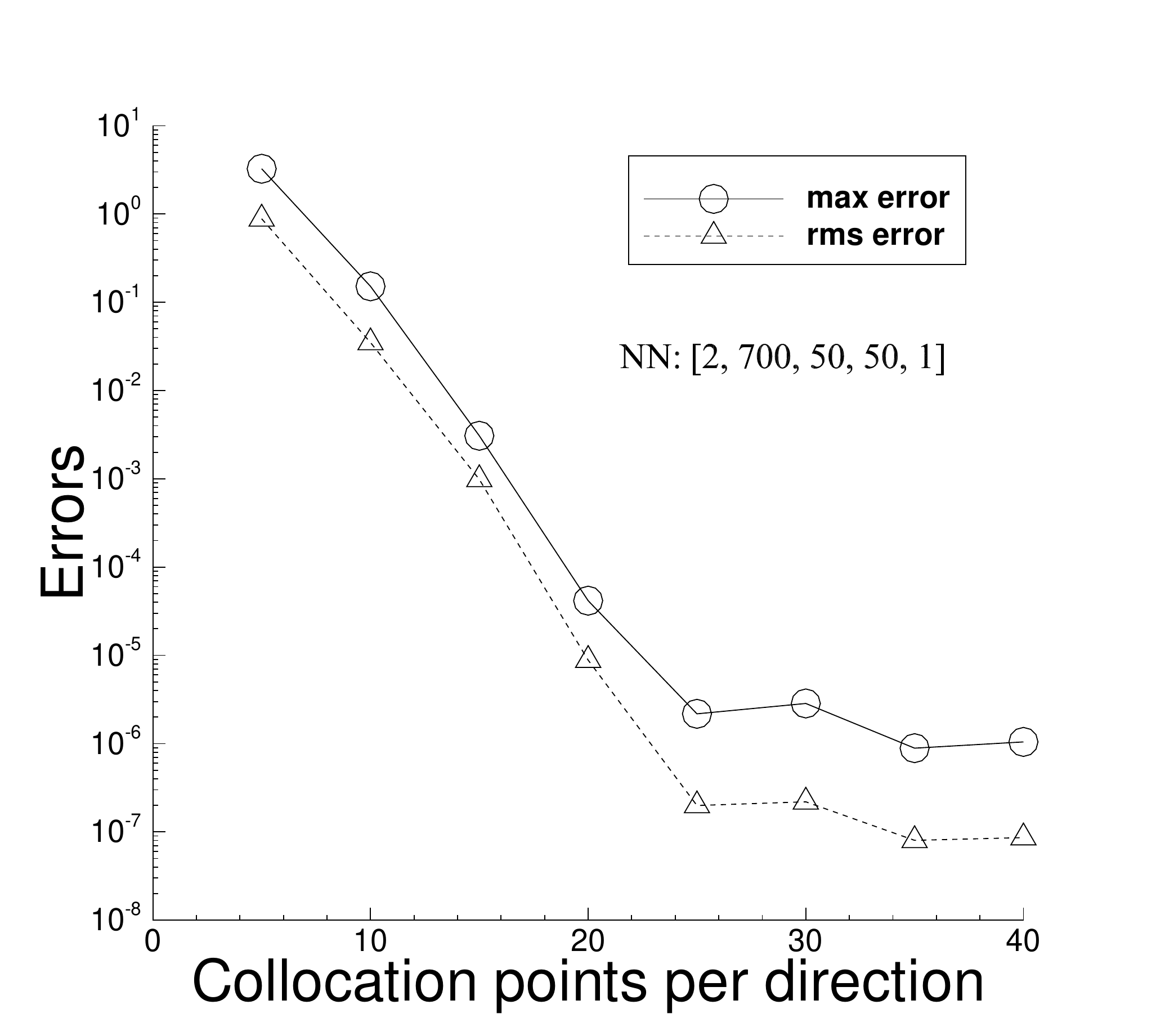}(b)
    \includegraphics[width=2in]{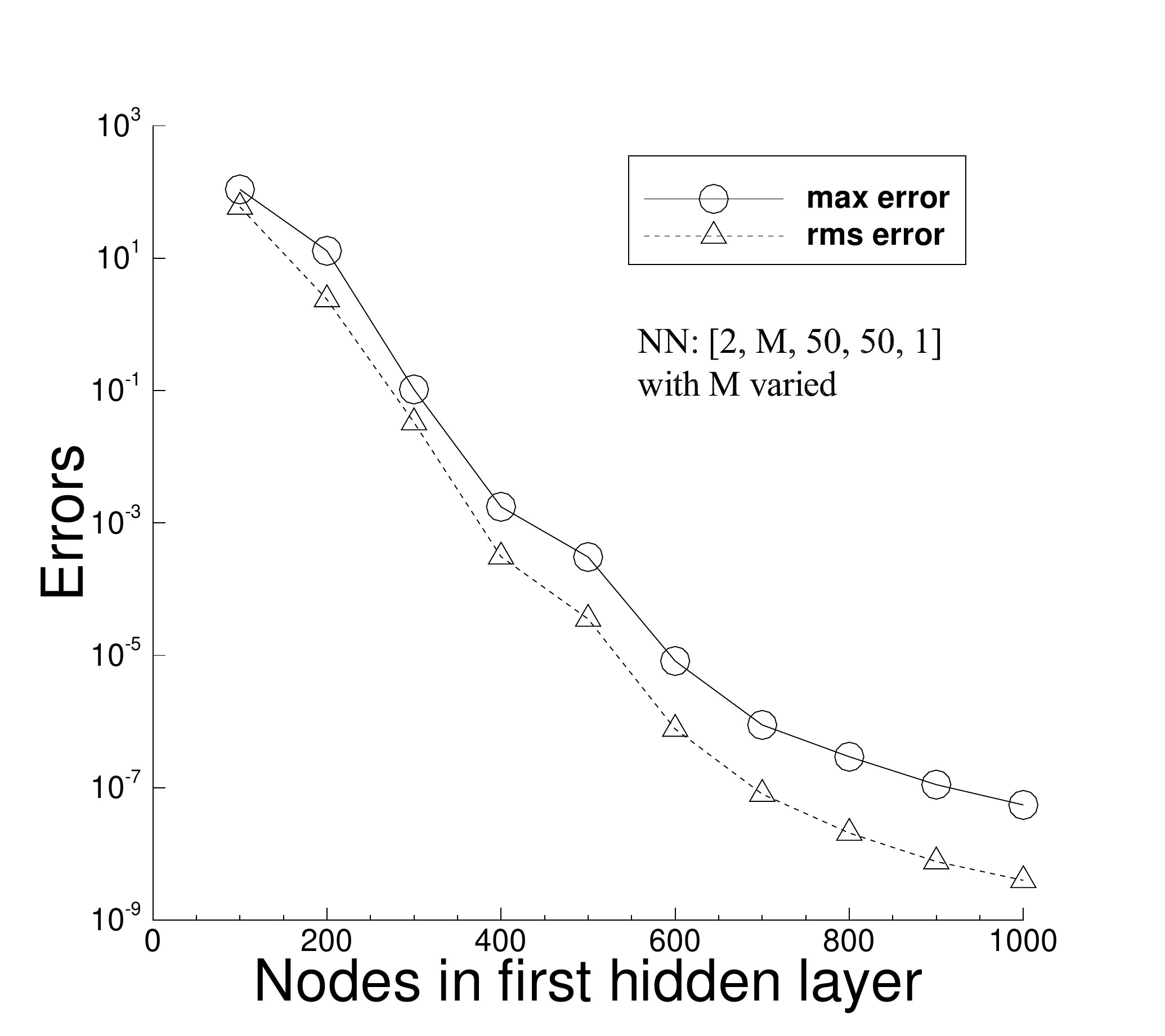}(c)
  }
  \caption{KdV equation (3 hidden layers in network):
    (a) Distribution of the absolute error of
    the HLConcELM solution. The HLConcELM maximum/rms errors
    versus (b) the number of collocation points per direction,
    and (c) the number of nodes in the first hidden layer ($M$).
    Network architecture: $[2, M, 50, 50, 1]$.
    $M=700$ in (a,b), varied in (c). $Q=35\times 35$ in (a,c), varied in (b).
    $\mbs R=(3.0,0.025,1.6)$ in (a,b,c).
  }
  \label{fg_24}
\end{figure}

Figure~\ref{fg_24} illustrates the HLConcELM solutions
obtained on a deeper neural network with $3$ hidden layers.
The neural network architecture is given by $\mbs M=[2, M, 50, 50, 1]$,
where $M$ is either fixed at $M=700$ or varied systematically.
A hidden magnitude vector $\mbs R=(3.0,0.026,1.6)$
is employed in all the simulations.
Figure~\ref{fg_24}(a) shows the absolute-error distribution on $\Omega$,
which corresponds to $M=700$ and a set of $Q=35\times 35$
collocation points.
The HLConcELM result is highly accurate, with the maximum error
on the order of $10^{-7}$ on $\Omega$.
Figures~\ref{fg_24}(b) and (c) depict the
maximum/rms errors of HLConcELM as a function of
the number of collocation points and as a function of $M$,
respectively. The data again signify the exponential
(or near exponential)
convergence of the HLConcELM errors.

%% file: Summary.tex
\section{Concluding Remarks}
\label{sec:summary}


The extreme learning machine (ELM) method can yield
highly accurate solutions to linear and nonlinear PDEs.
In terms of the accuracy
and computational cost,
it exhibits a notable advantage
over the existing DNN-based PDE solvers and the
the traditional numerical methods such as the second-order
and high-order finite element methods (see~\cite{DongL2021,DongY2021}).
To achieve a high accuracy, the existing
ELM method~\cite{DongL2021,DongL2021bip,DongY2021}
requires the last hidden layer of the
 neural network to be wide.
So the ELM neural network typically contains
a large number of nodes in the last hidden layer,
irrespective of the rest of network configuration.
If the last hidden layer is narrow,
the ELM accuracy will suffer and tend to be poor,
even though the neural network may contain a large number
of the nodes in the other hidden layers.

In the current paper we have presented a method to overcome
the above drawback of the existing (conventional) ELM method.
The new method, termed HLConcELM (hidden-layer concatenated ELM),
can produce highly accurate solutions to linear/nonlinear
PDEs when the last hidden layer is wide, and when the last hidden
layer is narrow, in which case the conventional ELM completely losses accuracy.

The new method relies on a type of modified feedforward neural networks (FNN),
which exposes the hidden nodes in all the hidden layers to
the output nodes by incorporating a {\em logical} concatenation
of the hidden layers into the network.
We refer to this modified network as HLConcFNN (hidden-layer concatenated FNN)
in this paper.
In HLConcFNN every hidden node in the network directly influences
the nodes in the output layer, while in conventional FNN only
the hidden nodes in the last hidden layer directly influence
the output nodes.
HLConcFNNs have the interesting property that, if new hidden layers
are appended to an existing network architecture or new nodes are added
to an existing hidden layer, the representation capacity of the
resultant network architecture is guaranteed to be not smaller than that of
the original one.

The HLConcELM method adopts the HLConcFNN as its neural network.
It assigns random values to (and fixes) the weight/bias coefficients
in all the hidden layers of the neural network, while
the coefficients between the output nodes and all the hidden nodes of the network
are trained/computed by a linear or nonlinear least squares method.
Note that in HLConcELM every hidden node in the network is connected
to the output nodes because of the logical hidden-layer concatenation.
HLConcELMs have the property that, as new hidden layers
are appended to an existing network architecture, the representation
capacity of the HLConcELM associated with the resultant architecture
is guaranteed to be not smaller than that associated with the original one,
provided that the random coefficients in the resultant architecture
are assigned in an appropriate fashion.

In essence, when solving PDEs or approximating functions,
the ELM method performs an expansion
of the unknown field solution in terms of a set of random basis functions.
With conventional ELM, the random bases consist of the output fields of
the last hidden layer of the neural network.
With HLConcELM, on the other hand, the random bases consist of
the output fields of the hidden nodes in all the hidden layers
of the neural network. HLConcELM is able to harvest the degrees of freedom
provided by all the hidden nodes in the network, not limited to those
from the last hidden layer.
This is the essential difference
between the HLConcELM method and the conventional ELM method.


We have tested the current HLConcELM method on boundary value problems
and initial/boundary value problems involving a number of linear and nonlinear
PDEs.  In particular, we have compared
HLConcELM and the conventional ELM on network architectures
whose last hidden layer is narrow or wide. The numerical results
demonstrate that the current HLConcELM method yields highly accurate 
results on network architectures with both narrow and wide last hidden layers.
In contrast, the conventional ELM only achieves accurate results on
architectures with a wide last hidden layer, and with a narrow last hidden layer
 it exhibits poor or no accuracy.
The HLConcELM method displays an exponentially convergent behavior for smooth field
solutions, reminiscent of the traditional high-order
techniques~\cite{KarniadakisS2005,ZhengD2011,YangD2020,YangLD2019,YangD2019}.
Its numerical errors decrease exponentially or nearly exponentially
as the number of collocation points or the number of trainable parameters
increases.

%

%% file: Append.tex
\section*{Appendix A: Proofs of Theorems from Section~\ref{sec:method}}

\paragraph{Proof of Theorem~\ref{thm_1}:} 
\begin{proof}
  Consider an arbitrary $u(\bm\theta,\bm\beta,\mbs x)\in U(\Omega,\mbs M_1,\sigma)$,
  where $\bm\theta\in\mbb R^{N_{h1}}$ 
  and $\bm\beta\in\mbb R^{N_{c1}}$,
  with $N_{h1}=\sum_{i=1}^{L-1}(m_{i-1}+1)m_i$ and
  $N_{c1}=\sum_{i=1}^{L-1}m_i$.
  Let $w_{kj}^{(i)}$ ($1\leqslant i\leqslant L-1$, $1\leqslant k\leqslant m_{i-1}$,
  $1\leqslant j\leqslant m_i$) and $b^{(i)}_j$ ($1\leqslant i\leqslant L-1$,
  $1\leqslant j\leqslant m_i$) denote
  the hidden-layer weight/bias coefficients
  of the associated HLConcFNN($\mbs M_1,\sigma$), and
  let $\beta_{ij}$ ($1\leqslant i\leqslant L-1$, $1\leqslant j\leqslant m_i$) denote
   the output-layer coefficients of HLConcFNN($\mbs M_1,\sigma$).
  $u(\bm\theta,\bm\beta,\mbs x)$ is given by~\eqref{eq_7}.

  Consider a function $v(\bm\vartheta,\bm\alpha,\mbs x)\in U(\Omega,\mbs M_2,\sigma)$
  with $\bm\vartheta\in\mbb R^{N_{h2}}$ and $\bm \alpha\in\mbb R^{N_{c2}}$,
  where $N_{c2}=N_{c1}+n$, and $N_{h2}=N_{h1}+(m_{L-1}+1)n$.
  We will choose $\bm\vartheta$ and $\bm\alpha$ 
  such that $v(\bm\vartheta,\bm\alpha,\mbs x) = u(\bm\theta,\bm\beta,\mbs x)$.
  We construct $\bm\vartheta$ and $\bm\alpha$ by setting the hidden-layer and
  the output-layer coefficients of HLConcFNN($\mbs M_2,\sigma$) as follows.

  The HLConcFNN($\mbs M_2,\sigma$) has $L$ hidden layers. We set the weight/bias
  coefficients in its last hidden layer (with $n$ nodes) to arbitrary values.
  We set those coefficients that connect the output node and the $n$ nodes in
  the last hidden layer  to all zeros.
  For the rest of the hidden-layer coefficients and the output-layer coefficients
  in HLConcFNN($\mbs M_2,\sigma$), we use those corresponding coefficient values
  from the network HLConcFNN($\mbs M_1,\sigma$).

  More specifically, let $\xi_{kj}^{(i)}$ and $\eta_j^{(i)}$
  denote the weight/bias coefficients in
  the hidden layers, and $\alpha_{ij}$ denote the output-layer coefficients,
  of HLConcFNN($\mbs M_2,\sigma$) associated with the function
  $v(\bm\vartheta,\bm\alpha,\mbs x)$.
  We set these coefficients by,
  \begin{equation}
    \xi_{kj}^{(i)}=\left\{
    \begin{array}{ll}
      w_{kj}^{(i)}, & \text{for}\ 1\leqslant i\leqslant L-1,\
      1\leqslant k\leqslant m_{i-1},\
      1\leqslant j\leqslant m_i; \\
      \text{arbitrary value}, & \text{for}\ i=L,\
      1\leqslant k\leqslant m_{L-1},\ 1\leqslant j\leqslant n; 
    \end{array}
    \right.
  \end{equation}
  \begin{equation}
    \eta_j^{(i)}=\left\{
    \begin{array}{ll}
      b_j^{(i)}, & \text{for all}\ 1\leqslant i\leqslant L-1,\
      1\leqslant j\leqslant m_i; \\
      \text{arbitrary value}, & \text{for}\ i=L,\
      1\leqslant j\leqslant n; 
    \end{array}
    \right.
  \end{equation}
  \begin{equation}
    \alpha_{ij} = \left\{
    \begin{array}{ll}
      \beta_{ij},& \text{for}\ 1\leqslant i\leqslant L-1,\
      1\leqslant j\leqslant m_i; \\
      0, & \text{for}\ i=L,\
      1\leqslant j\leqslant n.
    \end{array}
    \right.
  \end{equation}

  With the above coefficients,
  the last hidden layer of the network HLConcFNN($\mbs M_2,\sigma$) may output
  arbitrary fields, which however have no effect on the output field of
  HLConcFNN($\mbs M_2,\sigma$) because $\alpha_{Lj}=0$ ($1\leqslant j\leqslant n$).
  The rest of the hidden nodes in HLConcFNN($\mbs M_2,\sigma$) and the
  output node of HLConcFNN($\mbs M_2,\sigma$) produce fields that are identical to
  those of the corresponding nodes in the network HLConcFNN($\mbs M_1,\sigma$).
  We thus conclude that $u(\bm\theta,\bm\beta,\mbs x)=v(\bm\vartheta,\bm\alpha,\mbs x)$.
  So $u(\bm\theta,\bm\beta,\mbs x)\in U(\Omega,\mbs M_2,\sigma)$,
  and the relation~\eqref{eq_9} holds.  
\end{proof}

\paragraph{Proof of Theorem~\ref{thm_2}:}
\begin{proof}
  We use the same strategy as that in the proof of Theorem~\ref{thm_1}.
  Consider an arbitrary $u(\bm\theta,\bm\beta,\mbs x)\in U(\Omega,\mbs M_1,\sigma)$,
  where $\bm\theta\in\mbb R^{N_{h1}}$ 
  and $\bm\beta\in\mbb R^{N_{c1}}$, with $N_{h1}=\sum_{i=1}^{L-1}(m_{i-1}+1)m_i$ and
  $N_{c1}=\sum_{i=1}^{L-1}m_i$.
  The hidden-layer coefficients
  of the associated HLConcFNN($\mbs M_1,\sigma$) are denoted
  by $w_{kj}^{(i)}$ ($1\leqslant i\leqslant L-1$, $1\leqslant k\leqslant m_{i-1}$,
  $1\leqslant j\leqslant m_i$) and $b^{(i)}_j$ ($1\leqslant i\leqslant L-1$,
  $1\leqslant j\leqslant m_i$), and the output-layer coefficients
  are denoted by $\beta_{ij}$ ($1\leqslant i\leqslant L-1$, $1\leqslant j\leqslant m_i$).
  $u(\bm\theta,\bm\beta,\mbs x)$ is given by~\eqref{eq_7}.

  Consider a function $v(\bm\vartheta,\bm\alpha,\mbs x)\in U(\Omega,\mbs M_2,\sigma)$
  with $\bm\vartheta\in\mbb R^{N_{h2}}$ and $\bm \alpha\in\mbb R^{N_{c2}}$,
  where $N_{c2}=N_{c1}+1$, and $N_{h2}=N_{h1}+(m_{s-1}+1)+m_{s+1}$ if $1\leqslant s\leqslant L-2$
  and $N_{h2}=N_{h1}+(m_{s-1}+1)$ if $s=L-1$.
  We construct $\bm\vartheta$ and $\bm\alpha$ by setting the hidden-layer and
  the output-layer coefficients of HLConcFNN($\mbs M_2,\sigma$) as follows.

  In  HLConcFNN($\mbs M_2,\sigma$) we set
  the weight coefficients that connect the extra node of layer $s$
  to those nodes in layer $(s+1)$ to all zeros,
  and we also set the weight coefficient that connects the
  extra node of layer $s$ with the output node to zero.
  We set the weight coefficients that connect the nodes
  of layer $(s-1)$ to the extra node of layer $s$ to arbitrary values,
  and also set the bias coefficient corresponding to the extra node
  of layer $s$ to an arbitrary value.
  For the rest of the hidden-layer and output-layer
  coefficients of HLConcFNN($\mbs M_2,\sigma$),
  we use those corresponding coefficient values
  from the network HLConcFNN($\mbs M_1,\sigma$).

  Specifically, let $\xi_{kj}^{(i)}$ and $\eta_j^{(i)}$ denote the weight/bias coefficients in
  the hidden layers, and $\alpha_{ij}$ denote the output-layer coefficients,
  of the HLConcFNN($\mbs M_2,\sigma$) associated with $v(\bm\vartheta,\bm\alpha,\mbs x)$.
  We set these coefficients by,
  \begin{equation}
    \xi_{kj}^{(i)}=\left\{
    \begin{array}{ll}
      w_{kj}^{(i)}, & \text{for all}\ (1\leqslant i\leqslant s-1,\ \text{or}\
      s+2\leqslant i\leqslant L-1),\
      1\leqslant k\leqslant m_{i-1},\
      1\leqslant j\leqslant m_i; \\
      w_{kj}^{(s)}, & \text{for}\ i=s,\
      1\leqslant k\leqslant m_{s-1},\
      1\leqslant j\leqslant m_s; \\
      \text{arbitrary value}, & \text{for}\ i=s,\
      1\leqslant k\leqslant m_{s-1},\ j=m_{s}+1; \\
      w_{kj}^{(s+1)}, & \text{for}\ i=s+1,\
      1\leqslant k\leqslant m_{s},\
      1\leqslant j\leqslant m_{s+1}; \\
      0, & \text{for}\ i=s+1,\
      k=m_s+1,\
      1\leqslant j\leqslant m_{s+1};
    \end{array}
    \right.
  \end{equation}
  \begin{equation}
    \eta_j^{(i)}=\left\{
    \begin{array}{ll}
      b_j^{(i)}, & \text{for all}\ 1\leqslant i\leqslant L-1,\ i\neq s,\
      1\leqslant j\leqslant m_i; \\
      b_j^{(s)}, & \text{for}\ i=s,\
      1\leqslant j\leqslant m_s; \\
      \text{arbitrary value}, & \text{for}\ i=s,\
      j=m_s+1;
    \end{array}
    \right.
  \end{equation}
  \begin{equation}
    \alpha_{ij} = \left\{
    \begin{array}{ll}
      \beta_{ij},& \text{for all}\ 1\leqslant i\leqslant L-1,\ i\neq s,\
      1\leqslant j\leqslant m_i; \\
      \beta_{sj}, & \text{for}\ i=s,\
      1\leqslant j\leqslant m_s; \\
      0, & \text{for}\ i=s,\ j=m_s+1.
    \end{array}
    \right.
  \end{equation}

  With the above coefficients,
  the extra node in layer $s$ of the network HLConcFNN($\mbs M_2,\sigma$) may output
  an arbitrary field, which however has no contribution to the output field of
  HLConcFNN($\mbs M_2,\sigma$).
  The rest of the hidden nodes and the
  output node of HLConcFNN($\mbs M_2,\sigma$) produce identical
  fields as the corresponding nodes in the network HLConcFNN($\mbs M_1,\sigma$).
  We thus conclude that $u(\bm\theta,\bm\beta,\mbs x)=v(\bm\vartheta,\bm\alpha,\mbs x)$.
  So $u(\bm\theta,\bm\beta,\mbs x)\in U(\Omega,\mbs M_2,\sigma)$
  and the relation~\eqref{eq_10} holds.
\end{proof}